\documentclass[12pt]{amsart}
\usepackage[T1]{fontenc}
\usepackage{amsmath}
\usepackage{amsthm}
\usepackage{amssymb}
\usepackage{amsfonts}
\usepackage{outlines}
\usepackage{dsfont}
\usepackage{stackengine}
\usepackage{comment}
\usepackage{cancel}
\usepackage{graphicx}
\usepackage{enumerate}
\usepackage{mathrsfs}
\usepackage{tikz-cd}
\usepackage{hyperref}
\usepackage[numbers]{natbib}
\usepackage{enumitem}
\usepackage{lipsum}
\usepackage[left = 3.5cm, right = 3.5cm, top = 4.2 cm, bottom = 4.2 cm]{geometry}

\newtheorem{thm}{Theorem} [section]

\newtheorem{lemma}[thm]{Lemma}
\newtheorem{cor}[thm]{Corollary}

\newtheorem{prop}[thm]{Proposition}

\title[Cylindrical tangent flows in mean curvature flow]{Cylindrical tangent flows in mean curvature flow}
\author{Sourav Ghosh}
\date{}

\address{Department of Mathematics, University of Notre Dame, Notre Dame, IN 46556}
\email{sourav.ghosh.1@icloud.com}

\begin{document}

\maketitle

\begin{abstract}
The only non-compact linearly stable singularity models for mean curvature flow are cylindrical, as shown by Colding-Minicozzi \cite{colding2012generic}. The uniqueness of blowups at singularities modeled on cylinders was established by Colding-Minicozzi \cite{colding2015uniqueness}. They also proved rigidity results for cylindrical singularities in their earlier work \cite{colding2015rigidity}. In this paper, we develop a different approach inspired by Székelyhidi \cite{szekelyhidi2020uniqueness} to study cylindrical singularities and prove uniqueness and rigidity results. 
\end{abstract}

\tableofcontents

\section{Introduction}\label{section 1}

A smooth family of hypersurfaces $\{M_t\}_{t \in [0,T)}$ in $\mathbb{R}^{n+1}$ is called a mean curvature flow if
\[
\frac{\partial}{\partial t} x = - H_{M_t}(x),
\]
where $H_{M_t}(x)$ denotes the mean curvature of $M_t$ at the point $x$. Mean curvature flow is the negative gradient flow of the area functional. If the initial hypersurface is smooth and compact, then the flow is known to develop singularities in finite time.

By work of Huisken \cite{huisken1990asymptotic}, White \cite{white1994partial}, and Ilmanen \cite{ilmanen1995singularities}, the singularities of mean curvature flow can be modeled by self-similar shrinking solutions to the flow. To examine a singularity, one blows up around the point where singular behavior occurs. For this we need the notion of a tangent flow \cite{ilmanen1995singularities}, which generalizes the construction of tangent cone from minimal surfaces.

By definition, a tangent flow is the limit of a sequence of rescalings
at a singularity. For instance, a tangent flow to $M_t$ at the origin in space-time is the limit of a sequence of rescaled flows 
\[
\frac{1}{\epsilon_i} M_{\epsilon_i^2 t}, \quad \text{where}\; \epsilon_i \to 0.
\]
A priori, diﬀerent sequences $\epsilon_i$ could give diﬀerent tangent flows.

Uniqueness of blowups is one of the most fundamental questions concerning
singularities in mean curvature flow, as it is essential for understanding the
behavior of the flow both before and after the singular time. There are many
uniqueness results in mean curvature flow under the assumption that the tangent
flow has multiplicity one. We briefly recall some related results.

When a tangent flow is compact, Schulze \cite{schulze2014uniqueness} proved that
it is unique. Colding-Minicozzi \cite{colding2015uniqueness} established
uniqueness of tangent flows when the tangent flow is a generalized cylinder.
Chodosh-Schulze \cite{chodosh2021uniqueness} proved uniqueness when the tangent
flow is smooth and asymptotically conical. The first uniqueness result for
tangent flows modeled on non-smooth shrinkers was obtained by
Lotay-Schulze-Sz\'ekelyhidi \cite{lotay2022neck} in the setting of Lagrangian
mean curvature flow. Stolarski \cite{stolarski2023structure} showed that for
flows with bounded mean curvature, if a tangent flow is a regular multiplicity
one stationary cone, then it is the unique tangent flow. Recently, the author
\cite{ghosh2025lagrangian} proved that for two-dimensional Lagrangian mean
curvature flows with bounded mean curvature, the tangent flow is unique even in
the presence of higher multiplicity.

A well-known conjecture of Huisken asserts that generic mean curvature flows in
\(\mathbb{R}^3\) develop only spherical and cylindrical singularities.
Chodosh-Choi-Mantoulidis-Schulze \cite{chodosh2024mean} proved that the mean
curvature flow starting from a generic surface in \(\mathbb{R}^3\) can be
continued uniquely through spherical and cylindrical singularities, provided
that these singularities have multiplicity one. The assumption that all
singularities have multiplicity one was recently established by
Bamler-Kleiner \cite{bamler2023multiplicity}. Therefore, cylindrical
singularities represent the most important class of non-compact singularities
in mean curvature flow.

The goal of this paper is to reprove the uniqueness result proven in \cite{colding2015uniqueness}, which is the following
\begin{thm}  {\label{thm 1.1}}
    Fix $\mathscr{C}^n = \mathbb{S}^k_{\sqrt{2k}} \times \mathbb{R}^{n-k} \subset \mathbb{R}^{n+1}$ a cylinder. Let $\mathcal{M} = (M_t)_{t \in (-t_1,0)}$ be a mean curvature flow so that the self-similar shrinking multiplicity one flow associated to $\mathscr{C}, \mathcal{M}_\mathscr{C}$, arises as a tangent flow to $\mathcal{M}$ at the space-time point $(0,0)$. Then $\mathcal{M}_\mathscr{C}$ is the unique tangent flow to $\mathcal{M}$ at $(0,0)$.
\end{thm}

We will use the idea introduced by Sz\'ekelyhidi \cite{szekelyhidi2020uniqueness}
to prove the uniqueness of the cylindrical tangent cone \(C\times\mathbb{R}\) for an area-minimizing hypersurface, where \(C\) is the Simons cone
\(C_S=C(S^3\times S^3)\). Colding-Minicozzi
[\citealp{colding2015uniqueness}, \citealp{colding2019regularity}] proved
Łojasiewicz inequalities by a direct geometric approach that relies on the
specific geometry of cylinders. In this paper, instead of exploiting the special
geometry of the cylinder, we perturb the cylinder and construct flows near it. We
will show that the key estimates obtained from the Łojasiewicz inequality are
implicitly encoded in the properties of the constructed rescaled flows. Since
our approach does not depend on the specific geometry of the cylinder, it can be
applied to establish uniqueness results for a broader class of shrinkers.

Let \(y_i\) be coordinates on the \(\mathbb{R}^{n-k}\) factor and let
\(\theta\) be a point in the \(\mathbb{S}^k_{\sqrt{2k}}\) factor. There are two
types of Jacobi fields on \(\mathscr{C}\)
[\citealp{colding2015uniqueness}, \citealp{sun2022generic}]. The first type is
\(y_i f_i(\theta)\), where \(f_i\) is an eigenfunction on
\(\mathbb{S}^k_{\sqrt{2k}}\) with eigenvalue \(\frac{1}{2}\). These Jacobi fields
correspond to rotations. Geometrically, \(f_i\) represents translations in the
spherical directions, while \(y_i\) represents translations in the axis
directions. The second type is \(y_i y_j-2\delta_{ij}\). However, it is known
[\citealp{colding2015rigidity}, \citealp{zhu2020ojasiewicz}] that the quadratic
polynomials \(y_i y_j-2\delta_{ij}\) in the kernel are not generated by
one-parameter families of shrinkers. Geometrically, these Jacobi fields
correspond to non-degenerate neck pinching in the corresponding axis directions
and rotations of these directions. In particular, the kernel contains
non-integrable Jacobi fields.

Denoting the dimension of the space of non-integrable Jacobi fields by \(N\), we
associate each such Jacobi field with an element of \(\mathbb{R}^N\). For
\(\alpha\in\mathbb{R}^N\) with sufficiently small \(\ell^\infty\)-norm, we
construct deformations \(\mathscr{T}^{\alpha,\mathscr{C}}\) of the cylinder
\(\mathscr{C}\). These flows \(\mathscr{T}^{\alpha,\mathscr{C}}\), modeled on
non-integrable Jacobi fields, are perturbations of cylinders defined on
space-time parabolic regions depending on \(\alpha\). Moreover, their maximal
existence times increase as the \(\ell^\infty\)-norm of \(\alpha\) decreases.
Similar constructions have appeared in the literature, including logarithmic
perturbations of cones constructed by Sz\'ekelyhidi \cite{szekelyhidi2020uniqueness}
and logarithmically decaying surfaces constructed by Adams-Simon
\cite{adams1988rates}. We then establish a Łojasiewicz-type inequality
[Lemma~\ref{lemma 4.8}] for our constructed flows.

During a time interval in which the given flow \(M_\tau\) is sufficiently close
to \(\mathscr{C}\), we regard it as a perturbation of the best approximation
within the family of flows \(\mathscr{T}^{\alpha,\mathscr{C}}\). Using a
non-concentration result [Corollary~\ref{cor 5.3}], we show that, relative to
\(\mathscr{T}^{\alpha,\mathscr{C}}\), the flow \(M_\tau\) is modeled by a Jacobi
field over \(\mathscr{C}\) with no degree-zero component. Therefore, we can apply
the three-annulus lemma [Lemma~\ref{lemma 5.6}] to conclude that, as time
progresses, \(M_\tau\) either diverges significantly from
\(\mathscr{T}^{\alpha,\mathscr{C}}\) or converges toward
\(\mathscr{T}^{\alpha,\mathscr{C}}\).

In Proposition~\ref{prop 6.8}, we will show that if \(M_\tau\) diverges from
\(\mathscr{T}^{\alpha,\mathscr{C}}\) over time, then it cannot remain close to any
shrinking solution, which leads to a contradiction. Therefore, we may assume that
\(M_\tau\) approaches \(\mathscr{T}^{\alpha,\mathscr{C}}\). However, the constructed flow \(\mathscr{T}^{\alpha,\mathscr{C}}\) is defined only
for a finite time interval, whose length depends on \(\alpha\). Consequently, the
three-annulus lemma [Lemma~\ref{lemma 5.6}] can be applied only finitely many times. Moreover, since \(\mathscr{T}^{\alpha,\mathscr{C}}\) is not a self-similar shrinking solution, the fact that \(M_\tau\) approaches \(\mathscr{T}^{\alpha,\mathscr{C}}\) does not necessarily imply that \(M_\tau\) becomes sufficiently cylindrical at a suitable rate. Hence, this alone is not enough to conclude the uniqueness of the flow.

Nevertheless, we will show that, in this case, the given flow is very close to
\(\mathscr{T}^{\alpha,\mathscr{C}}\) relative to \(\alpha\). Combining this
estimate with the Łojasiewicz-type inequality [Lemma~\ref{lemma 5.6}] and Huisken's
monotonicity formula, we obtain a quantitative decay rate for the Gaussian area of
the rescaled flow [Proposition~\ref{prop 6.4}]. Once this decay estimate is
established, we can ultimately conclude the uniqueness of the tangent flow.

The structure of the paper is summarized as follows. In Section~\ref{section 2},
we state preliminary results in mean curvature flow and describe the space of
Jacobi fields on \(\mathscr{C}\). In Section~\ref{section 3}, we discuss the
graphical rescaled mean curvature flow equation. In Section~\ref{section 4}, we
construct the rescaled flow \(\mathscr{T}^{\alpha,\mathscr{C}}\), which is
defined on a space-time parabolic region depending on \(\alpha\) and modeled on
the non-integrable Jacobi field \(u_\alpha\) constructed in that section. In Section~\ref{section 5}, we define the distance function $D_{\mathscr{T}^{\alpha,\mathscr{C}}_\tau}(M_{\tau_0}),$
which measures the distance between the given rescaled flow at time \(\tau_0\)
and the constructed flow \(\mathscr{T}^{\alpha,\mathscr{C}}\) at time \(\tau\).
We then prove the non-concentration result [Proposition~\ref{prop 5.2} and Corollary~\ref{cor 5.3}]. Using this non-concentration result, we establish a three-annulus-type lemma [Lemma~\ref{lemma 5.6}] for
\(D_{\mathscr{T}^{\alpha,\mathscr{C}}_\tau}(M_{\tau_0})\). In Section~\ref{section 6}, we apply this three-annulus lemma to prove the decay estimate [Proposition~\ref{prop 6.8}] for the distance function. Finally, given
the decay estimate in Proposition~\ref{prop 6.8}, the uniqueness result
[Theorem~\ref{thm 1.1}] follows by standard arguments.

\subsection{Acknowledgments}
I would like to express my gratitude to my advisor, Professor G\'abor Sz\'ekelyhidi, for introducing me to this circle of ideas and for his unwavering support and encouragement throughout the project. I am deeply thankful to Professor Nicholas Edelen for many insightful discussions at various stages. I am also grateful to Professor William Minicozzi for his valuable feedback and helpful comments.
\section{Preliminaries}\label{section 2} 

In this section, we briefly discuss the basic definitions in mean curvature flow (MCF) that will be useful throughout the article.

We consider smooth, properly immersed submanifolds $\Sigma^n \subset \mathbb{R}^N$. We denote by $x$ the position vector in $\mathbb{R}^N$. For a vector $V$, we write $V^T$ for its projection onto the tangent bundle $T\Sigma$, and $V^\perp=\Pi(V)$ for its projection onto the normal bundle $N\Sigma$.

The second fundamental form is the symmetric $2$-tensor with values in the normal bundle defined by
\[
A(Y,Z)=(\nabla_Y Z)^\perp,
\]
and the mean curvature vector is given by
\[
H=\operatorname{tr}_{\Sigma}A.
\]
We define the shrinker mean curvature by
\[
\phi=-H+\frac12 x^\perp.
\]
A submanifold is called a self-shrinker if $\phi\equiv0$ on $\Sigma$.

For $\lambda>0$, we define the parabolic rescaling
\begin{align*}
\mathcal{D}_\lambda:\mathbb{R}^n\times\mathbb{R}
&\to
\mathbb{R}^n\times\mathbb{R},\\
(x,t)
&\mapsto
(\lambda x,\lambda^2 t).
\end{align*}
If $\mathcal{M}$ denotes the space-time track of a MCF, then $\mathcal{D}_\lambda\mathcal{M}$ is again a MCF.

Suppose that a solution to mean curvature flow develops a singularity at the space-time point $(x_0,T)$. For any sequence $\lambda_i\to\infty$, the rescaled flows
\[
\mathcal{D}_{\lambda_i}\bigl(\mathcal{M}-(x_0,T)\bigr)
\]
admit a subsequence converging as Brakke flows to a limit $\{\nu_t\}_{t\in(-\infty,0)}$ \cite[Lemma ~7.1]{ilmanen1994elliptic}. Any such limit is called a tangent flow at $(x_0,T)$. In general, the tangent flow may depend on the sequence $\lambda_i$. Moreover, one can prove that every tangent flow is a self-shrinking Brakke flow \cite[Lemma ~8]{ilmanen1995singularities}.

We recall Huisken's monotonicity formula \cite{huisken1990asymptotic}:
\[
\begin{aligned}
\frac{d}{dt}\int_{M_t} f\,\rho_{x_0,t_0}\,d\mathcal{H}^n
&=
\int_{M_t}(\partial_t f-\Delta f)\rho_{x_0,t_0}\,d\mathcal{H}^n\\
&\hspace{0.6 cm}-
\int_{M_t}
f\left|
-H+\frac{(x-x_0)^\perp}{2(t-t_0)}
\right|^2
\rho_{x_0,t_0}\,d\mathcal{H}^n .
\end{aligned}
\]
for $t<t_0$, where $f$ is a function on $M_t$ with polynomial growth (locally uniformly in $t$), and
\[
\rho_{x_0,t_0}
=
\frac{1}{(4\pi(t_0-t))^{n/2}}
\exp\left(
-\frac{|x-x_0|^2}{4(t_0-t)}
\right)
\]
is the backward heat kernel.

Given $x_0\in\mathbb{R}^{n+1}$ and $t_0>0$, the functional
\[
F_{x_0,t_0}(M)
=
\frac{1}{(4\pi t_0)^{n/2}}
\int_M
e^{-\frac{|x-x_0|^2}{4t_0}}
\,d\mathcal{H}^n
\]
introduced in \cite{colding2012generic} plays a central role. Note that $M$ is a critical point of $F_{x_0,t_0}$ if and only if it is the time $t=-t_0$ slice of a self-shrinking solution to mean curvature flow that becomes extinct at $x=x_0$ and $t=0$.

The entropy of $M$ is defined by
\[
\lambda(M)
=
\sup_{x_0,t_0}
F_{x_0,t_0}(M).
\]
By Huisken's monotonicity formula, the quantity $t\mapsto\lambda(M_t)$ is decreasing along any mean curvature flow.

We say that $M_t$ has uniformly bounded area ratios if there exists a constant $C_1>0$ such that
\[
\sup_{x,t}
\mathcal{H}^n(M_t\cap B_r(x))
\leq
C_1 r^n
\]
for all $r>0$, where $B_r(x)$ denotes the Euclidean ball of radius $r$ centered at $x\in\mathbb{R}^{n+1}$.

Note that this volume growth condition is equivalent to uniformly bounded entropy.

The Gaussian weight is
\[
\rho
=
\frac{1}{(4\pi)^{n/2}}
e^{-\frac{|x|^2}{4}}.
\]
By $L^2$ and $W^{k,2}$ we denote the weighted Sobolev spaces with respect to $\rho$, and we write $H^k$ in place of $W^{k,2}$.

The Gaussian area functional is defined by
\[
F(M)
=
\int_M \rho\,d\mathcal{H}^n.
\]
We define the excess of $M$ by
\[
\mathcal{A}(M)
=
F(M)-F(\mathscr{C}).
\]
Let $(M_t)_{0\leq t<T}$ be a solution to mean curvature flow and let
$x_0\in\mathbb{R}^{n+1}$. We consider the rescaled flow centered at
$(x_0,T)$, defined by
\[
M_\tau
=
\frac{1}{\sqrt{T-t}}
\bigl(M_t-x_0\bigr),
\qquad
\tau=-\log(T-t).
\]
A straightforward computation shows that the rescaled flow evolves with normal velocity
\[
-H+\frac12 x^\perp.
\]
In terms of the rescaled flow, tangent flows can be studied by taking limits of sequences $M_{\tau_i}$ as $\tau_i\to\infty$. Moreover, Huisken's monotonicity formula becomes
\[
\frac{d}{d\tau}
\int_{M_\tau}
\rho\,d\mathcal{H}^n
=
-
\int_{M_\tau}
\left|
-H+\frac12 x^\perp
\right|^2
\rho\,d\mathcal{H}^n.
\]
Consequently, $\mathcal{A}$ is decreasing along the rescaled flow.

We denote by $A$ the second fundamental form. Consider the elliptic operator
\[
Lu
=
\Delta u
-\frac12 x\cdot\nabla u
+
|A|^2u
+
\frac12 u.
\]
Since on $\mathscr{C}$ we have
\[
|A|^2=\frac12,
\]
the operator reduces to
\[
Lu
=
\Delta u
-\frac12 x\cdot\nabla u
+
u.
\]
We denote by $K$ the $H^2$-kernel of $L$. The linearized operator $L$ is self-adjoint with respect to the weighted $L^2$ inner product
\[
\langle u,v\rangle
=
\int_{\mathscr{C}}
u(x)v(x)
e^{-\frac{|x|^2}{4}}
\,d\mathcal{H}^n.
\]
We say that $\lambda$ is an eigenvalue of $L$ if the equation
\[
Lf=\lambda f
\]
admits a nontrivial $L^2$ solution. In Sun-Wang-Xue [\citealp{sun2020multiplicity}, Section ~5.2], it was shown that the spectrum of the elliptic operator $L$ is given by
\[
\sigma(L)
=
\left\{
-\mu_i-\frac{j}{2}+1
\right\}_{i,j=0}^{\infty},
\]
where
\[
\mu_i
=
\frac{i(i-1+k)}{2k},
\]
and the corresponding eigenfunctions are spanned by
\[
\{\phi_i(\theta)h_j(y)\}_{i,j=0}^{\infty}.
\]
Here $\phi_i$ are eigenfunctions of
$\Delta_{\mathbb{S}^k_{\sqrt{2k}}}$,
and $h_j$ are Hermite polynomials of degree $j$ on
$\mathbb{R}^{n-k}$.

In particular, the first four distinct eigenvalues of $L$ are
\[
\lambda_1=1,
\qquad
\lambda_2=\frac12,
\qquad
\lambda_3=0,
\qquad
\lambda_4
=
\max\left\{
-\frac12,
-\frac1k
\right\}.
\]
We now describe the space of Jacobi fields on the cylinder. Let $K$ denote the space of Jacobi fields on $\mathscr{C}$. The following explicit description is given in \cite[Lemma ~3.25]{colding2015uniqueness} and \cite[Proposition ~3.1, Corollary ~3.6]{colding2019regularity}.

\begin{lemma}\label{lemma 2.1}
We have
\[
K=K_1\oplus K_2.
\]
An orthogonal basis for $K_1$ is given by
\[
\left\{
y_i f :
f \text{ is an eigenfunction on }
\mathbb{S}^k
\text{ with eigenvalue }
\frac12
\right\},
\]
and an orthogonal basis for $K_2$ is given by
\[
\{y_i y_j-2\delta_{ij}: i\leq j\}.
\]
\end{lemma}

One can verify that $K_1$ corresponds to infinitesimal rotations. To make this precise, consider the natural action of $SO(n+1)$ on the ambient space. This induces a map from the Lie algebra $\mathfrak{so}(n+1)$ to the space of normal vector fields along the submanifold, given by
\[
\theta\mapsto J_\theta,
\]
where
\[
J_\theta(p)
=
\left.
\frac{d}{ds}
\right|_{s=0}
\bigl(\exp(s\theta)(p)\bigr)^\perp.
\]
The image of this map is precisely $K_1$. Let
\[
N:=\dim K_2.
\]
Fix an orthonormal ordered basis
\[
\{U_1,U_2,\ldots,U_N\}
\]
of $K_2$. Then every $u\in K_2$ can be written uniquely as
\[
u
=
\sum_{j=1}^{N}
\alpha_jU_j,
\qquad
\alpha_j\in\mathbb{R}.
\]
Thus, we identify $u$ with its coordinate vector $\alpha = (\alpha_1,\alpha_2,\ldots,\alpha_N)
\in \mathbb{R}^N,$ and write $u=u_\alpha$. We define
\[
|\alpha|
=
N\max_{1\leq j\leq N}|\alpha_j|.
\]
Conversely, given any $\alpha = (\alpha_1,\alpha_2,\ldots,\alpha_N) \in \mathbb{R}^N,$ we associate an element $u_\alpha\in K_2$ by
\[
u_\alpha
=
\sum_{j=1}^{N}
\alpha_jU_j.
\]
\section{Graphical Representation}\label{section 3}

Let $\Sigma$ be a hypersurface in $\mathbb{R}^{n+1}$, and let $n$ denote the unit normal to $\Sigma$. Let $u$ be a function on $\Sigma$, and define the graph $\Sigma_u$ by
\[
\Sigma_u = \{p + u(p)n(p) : p \in \Sigma\}.
\]

If $u$ is sufficiently small, then $\Sigma_u$ is contained in a tubular neighborhood of $\Sigma$ on which the normal exponential map is invertible. Let $e_{n+1}$ denote the gradient of the signed distance function to $\Sigma$.

We denote by $\nu_u(p)$ the relative area element, given by
\[
\nu_u(p) = \frac{\sqrt{\det g^u_{ij}(p)}}{\sqrt{\det g_{ij}(p)}},
\]
where $g_{ij}$ is the metric on $\Sigma$ and $g^u_{ij}$ is the pull-back metric on $\Sigma_u$. The relative area element was used in \cite{colding2015uniqueness} to compute the mean curvature and relate the gradient of the $F$-functional to the geometry of the graph.

We also denote by $n_u(p)$ the unit normal to $\Sigma_u$ at $p + u(p)n(p)$, by $H_u(p)$ its mean curvature, by
\[
\eta_u(p) = \langle p + u(p)n(p), n_u(p) \rangle
\]
the support function, and by
\[
w_u(p) = \langle e_{n+1}, n_u(p) \rangle^{-1}
\]
the corresponding weight function.

We have the following lemma (see \citealp{colding2015uniqueness}, Lemma ~A.44, and also \cite{colding2019dynamics}):

\begin{lemma}\label{lemma 3.1}
For sufficiently small $u$, the graphs $\Sigma_u$ over $\Sigma$ evolve by the rescaled mean curvature flow if and only if $u$ satisfies
\[
\partial_t u(p,t) = w_u(p)\Big(-H_u(p) + \tfrac{1}{2}\eta_u(p)\Big) =: M(u).
\]
\end{lemma}

The following lemma gives the decomposition of $M$ into its linearization at $u=0$ and a nonlinear remainder.

\begin{lemma} \label{lemma 3.2}
Let $M$ be the nonlinear operator appearing in Lemma ~\ref{lemma 3.1}. There exists $c >0$ such that for every $u \in C^{2,\alpha}(\Sigma),$ with $\|u\|_{C^{2, \alpha}(\Sigma)} \leq c$ the operator $M$ admits the decomposition
\[
M(u) = Lu + Q(u).
\]
where $L = DM|_0$ is the linearization of $M$ at $u = 0,$ 
\[
Lu = \Delta u - \frac{1}{2}x \cdot \nabla u + |A|^2 u + \frac{1}{2}u.
\]
and $Q$ is a nonlinear map satisfying
$$Q(0) = 0,  \qquad DQ(0) = 0.$$
In particular, $Q$ vanishes to second order at the origin.
\end{lemma}

When the shrinker is $\mathscr{C}$, the nonlinear term $Q$ satisfies the estimate
\[
|Q(u)| \le C\Big(
|\nabla u|^4
+ |\nabla u|^2 |\nabla^2 u|
+ |\nabla u|^2
+ u^2
+ |u|\,|\nabla^2 u|
\Big),
\]
for all $u$ with $\|u\|_{C^{2,\alpha}(\mathscr{C})} \le c$, where $C>0$ depends only on $\mathscr{C}$.

We now extend this formulation to the case where the reference surface is itself evolving under a rescaled mean curvature flow. Let $\{M_t\}_{t \in [0,T]}$ be a rescaled MCF. For each $t$, let $u(\cdot,t)$ be a function defined on $M_t$, and define $M_{u,t}$ as the normal graph of $u(\cdot,t)$ over $M_t$. We denote by $w(t,\cdot)$, $\nu(t,\cdot)$, and $\eta(t,\cdot)$ the corresponding geometric quantities on each hypersurface $M_t$.

We consider $u$ as a function on the space-time manifold $\bigcup_{t \in [0,T]} M_t \times \{t\}$. Then $M_{u,t}$ evolves by rescaled mean curvature flow if and only if $u$ satisfies a quasilinear parabolic equation generalizing the stationary case. This formulation was derived in Sun--Xue \citealp{sun2021initial}, Lemma ~A.3.

\begin{lemma}\label{lemma 3.3}
For sufficiently small $u$, the family of graphs $M_{u,t}$ evolves by rescaled mean curvature flow if and only if $u$ satisfies
\[
\partial_t u(p,t)
=
w\big(p,t,u(p,t),\nabla u(p,t)\big)
\Big(
- H_u
+ \tfrac{1}{2}\eta\big(p,t,u(p,t),\nabla u(p,t)\big)
+ F_t(u)
\Big),
\]
where
\[
F_t(u)
=
H_0
- \tfrac{1}{2}\langle p,n\rangle
- u \left\langle \nabla \left(H_0 - \tfrac{1}{2}\langle p,n\rangle\right), n_u \right\rangle.
\]
\end{lemma}

\section{The comparison flows}\label{section 4}

In this section, we construct the rescaled MCF near a cylinder. Our approach follows the strategy of [\citealp{szekelyhidi2020uniqueness}, Section ~4]. We begin with a lemma analogous to Lemma ~3.22 of \cite{colding2015uniqueness}. Roughly speaking, it states that any $C^{k,\alpha}$ function on the cylinder grows at most quadratically in directions orthogonal to $\ker L$.

\begin{lemma} {\label{lemma 4.1}}
    Let $R \geq 1$ and let $u \in C^{k,\alpha}(\mathscr{C})$. Then,
    $$\|P_{(\ker L)^\perp} u\|_{C^{k,\alpha}(\mathscr{C} \cap B_R(0))} \leq CR^2 \|u\|_{C^{k,\alpha} (\mathscr{C})}.$$
\end{lemma}

\begin{proof}
    Assume that
\[
\|u\|_{C^0(\mathscr C)} = 1.
\]
Then, in particular,
\[
\|u\|_{C^{k,\alpha}(\mathscr C)} \geq 1.
\]
Since \(\mathscr C\) has finite Gaussian volume, there exists a constant \(C\) such that
\[
\|u\|_{L^2(\mathscr C)} \leq C.
\]
Since \(P_{\ker L}\) is the orthogonal projection onto \(\ker L\), it follows that
\[
\|P_{\ker L}u\|_{L^2(\mathscr C)} \leq C.
\]
Applying Lemma~\ref{lemma 2.1}, we obtain
\[
\|P_{\ker L}u\|_{C^{k,\alpha}(\mathscr C\cap B_R(0))}
\leq CR^2.
\]
Using the fact that \(\|u\|_{C^{k,\alpha}(\mathscr C)} \geq 1\), we further deduce that
\[
\|P_{\ker L}u\|_{C^{k,\alpha}(\mathscr C\cap B_R(0))}
\leq CR^2 \|u\|_{C^{k,\alpha}(\mathscr C)}.
\]
Therefore,
\begin{align*}
\|P_{(\ker L)^\perp}u\|_{C^{k,\alpha}(\mathscr C\cap B_R(0))}
&=
\|u-P_{\ker L}u\|_{C^{k,\alpha}(\mathscr C\cap B_R(0))} \\
&\leq
\|u\|_{C^{k,\alpha}(\mathscr C\cap B_R(0))}
+
\|P_{\ker L}u\|_{C^{k,\alpha}(\mathscr C\cap B_R(0))} \\
&\leq
(1+CR^2)\|u\|_{C^{k,\alpha}(\mathscr C)}.
\end{align*}
Since $R\geq 1$, it follows that
\[
(1+CR^2)\|u\|_{C^{k,\alpha}(\mathscr C)}
\leq CR^2\|u\|_{C^{k,\alpha}(\mathscr C)},
\]
after enlarging the constant $C$ if necessary.
\end{proof}

We now prove a lemma analogous to \cite[Lemma ~3.15]{chodosh2021uniqueness}. The argument is essentially a maximum principle argument. It shows that if \(Lu=f\) in the Gaussian space and \(f\) has finite \(C^0\)-norm, then \(u\) satisfies a growth estimate in terms of \(\|f\|_{C^0}\) and the \(C^0\)-norm of \(u\) on a fixed ball.

In our setting, however, we must work with functions lying in \((\ker L)^\perp\), obtained as projections of functions defined on the entire cylinder with finite \(C^0\)-norm. Accordingly, we adapt the argument to this projected setting. 

\begin{lemma} {\label{lemma 4.2}}
    For $f \in L^2(\mathscr{C}) \cap C^0(\mathscr{C}),$ if $Lu =f - P_{\ker L} f$ for $u \in W^{1,2}(\mathscr{C}),$ then there exists a $R_0 > 0$ such that,
    $$\sup\limits_{x\in \mathscr{C}} \big( 1 + |x|^4 \big)^{-1} |u(x)| \leq C\Big( \sup\limits_{x\in \mathscr{C}} |f| + \sup\limits_{\mathscr{C} \cap B_{R_0 (0)}}  |u| \Big)$$
\end{lemma}

\begin{proof}
First note that, if we write $P_{\ker L}f = \sum a_i U_i, U_i$ orthonormal basis of $\ker L,$ then 
$$|a_i| \leq \|f\|_{L^2} \leq C \sup _{\mathscr{C}} |f|.$$
Therefore, 
$$|P_{\ker L}f| \leq C(n) (1+|y|^2) \sup \limits_{\mathscr{C}} |f|.$$
Let $\phi: \mathbb{R}^{n+1} \to \mathbb{R},$ we have
\begin{align*}
    L\phi &= \Delta_\mathscr{C} \phi - \frac{1}{2} x^T. \nabla \phi + \phi  \\
    &= \Delta_{\mathbb{R}^{n+1}} \phi - D^2 \phi (\nu_\mathscr{C}, \nu_\mathscr{C}) - \frac{1}{2} x. \nabla_{\mathbb{R}^{n+1}} \phi + \phi.
\end{align*}
Let 
$$\phi = \alpha |y|^2 - \beta + \gamma |y|^4,$$ 
with constants $\alpha, \beta, \gamma >0$ to be chosen. We compute
$$L\phi = 2\alpha (n-k) - \beta + 4\gamma (n-k+2) |y|^2 - \gamma |y|^4.$$
Fix $R > 0$ be chosen later. Set, 
$$\alpha = \sup \limits_\mathscr{C} |f| + \sup \limits_{\mathscr{C} \cap \partial B_R(0)} |u|.$$
$$\beta = 4 \alpha (n-k).$$
$$\gamma = C(n)\sup \limits_\mathscr{C} |f|.$$
Set 
$$v = u -\phi.$$ 
Then
\begin{align*}
    Lv &= Lu - L\phi \\
    &= f - P_{\ker L} f - 2 \alpha (n-k) + \beta - 4\gamma (n-k+2) |y|^2 + \gamma |y|^4 \\
    &= f - P_{\ker L} f  + 2 \alpha (n-k) - 4\gamma (n-k+2) |y|^2 + \gamma |y|^4 \\
    &\geq f - C(n)(1+|y|^2) \sup_{\mathscr{C}} |f| + 2 (n-k) \Big(\sup_\mathscr{C} |f| + \sup \limits_{\mathscr{C} \cap \partial B_R(0)} |u| \Big) \\
    &\hspace{0.5 cm} - C(n)\sup _\mathscr{C} |f| (4n - 4k + 8) |y|^2 + C(n) \sup _\mathscr{C} |f| |y|^4    \\
    &> 0,
\end{align*}
provided $|y|^4 - (4n - 4k + 9)|y|^2-1 \geq 0$, i.e. $|x|^2 = |y|^2 + 2k \geq C_1(n)$ for some dimensional constant $n.$ Now on $\partial B_R(0),$
\begin{align*}
    \sup_{\mathscr{C} \cap \partial B_R(0)}  v 
    &\leq \sup_{\mathscr{C} \cap \partial B_R(0)}u 
    - (R^2-2k) \Big(\sup \limits_\mathscr{C} |f| + \sup \limits_{\mathscr{C} \cap \partial B_R(0)} |u| \Big) \\
    &\hspace{0.5 cm} + 4(n-k) \Big(\sup_\mathscr{C} |f| + \sup_{\mathscr{C} \cap \partial B_R(0)} |u| \Big) 
    - (R^2-2k)^2 C(n) \sup_\mathscr{C} |f|\\
    &\leq \sup \limits_\mathscr{C} |f| 
    \Big(4n - 2k - R^2 - C(n)(R^2-2k)^2\Big)\\
    &\hspace{0.5 cm} + \sup_{\mathscr{C} \cap \partial B_R(0)} |u|
    \Big(1 + 4n - 2k - R^2\Big) \\
    &< 0.
\end{align*}
provided $R^2 > 4n + 1$. Define 
$$v^+:= \max \{v,0 \}.$$
Now if we choose any $R$ such that $R^2> \max\{4n+1, C_1(n)\}$, then
\begin{align*}
    0 &\leq \int_{\mathscr{C} \setminus B_R(0)} v^+ Lv \;e^{-|x|^2/4} \\
    &= \int_{\mathscr{C} \setminus B_R(0)} v^+ Lv^+ \;e^{-|x|^2/4}\\ 
    &= \int_{\mathscr{C} \setminus B_R(0)} \big(-|\nabla v^+|^2 + |v^+|^2 \big) \;e^{-|x|^2/4},
\end{align*}
where in the last line we have used integral by parts and the fact that $v^+ = 0$ near $\mathscr{C} \cap \partial B_R(0)$. Now,
$$R^2 \int_{\mathscr{C} \setminus B_R(0)} |v^+|^2 \;e^{-|x|^2/4} \leq \int_{\mathscr{C} \setminus B_R(0)} |x|^2|v^+|^2 \;e^{-|x|^2/4}.$$    
Also using Gaussian Poincar\'e type inequality [See for instance \cite{colding2015uniqueness}, Lemma ~3.4], we have
$$\int_{\mathscr{C} \setminus B_R(0)} |x|^2|v^+|^2 \;e^{-|x|^2/4} \leq C_2(n) \int_{\mathscr{C} \setminus B_R(0)} \big(|\nabla v^+|^2 + |v^+|^2 \big) \;e^{-|x|^2/4}.$$
Choose any $R_0$ such that 
$$R_0^2 > \max \{4n+1, C_1(n), C_2(n)\}.$$
As, 
$$R_0^2 \int_{\mathscr{C} \setminus B_{R_0}(0)} |v^+|^2 \;e^{-|x|^2/4} \leq C_2(n) \int_{\mathscr{C} \setminus B_{R_0}(0)} |v^+|^2 \;e^{-|x|^2/4}.$$
it will imply, $v^+ = 0$ as $R_0^2 > C_2(n)$. So, $u \leq \phi$ on $\mathscr{C} \setminus B_{R_0}(0)$. Now we can apply the same reasoning to $-u$ to finish the proof. 
\end{proof}

We now prove Schauder-type estimates for the operator \(L\) on the truncated cylinder \(\mathscr{C}\cap B_R(0)\). These estimates generalize Proposition~8.8 of \cite{kapouleas2018mean} to our setting; see also Proposition~3.5 of \cite{chodosh2021uniqueness}.

Let \(x^T\) denote the tangential component of the position vector along \(\mathscr{C}\), and let \(\varphi_t:\mathscr{C}\to\mathscr{C}\) be the one-parameter family of diffeomorphisms generated by the vector field
\[
-\frac{1}{2t}x^T,
\]
defined for \(t\in[-1,-\tfrac12]\), with \(\varphi_{-1}\) equal to the identity map. For the cylinder \(\mathscr{C}\), the map \(\varphi_t\) can be computed explicitly; see \cite{ecker2012regularity}.

Suppose $F:\mathscr{C}\times[-1,0)\to\mathbb{R}^N$ is a family of embeddings
satisfying
\[
\frac{\partial F}{\partial t}=H.
\]
Define the tangentially reparametrized family
\[
\tilde{F}:=F\circ\varphi_t^{-1}.
\]
Then $\tilde{F}$ satisfies
\[
\left(\frac{\partial \tilde{F}}{\partial t}\right)^{\perp}=H,
\]
where \((\cdot)^{\perp}\) denotes projection onto the normal bundle. Writing $x_1,...,x_{k+1}$ for coordinates on $\mathbb{S}^k$ and $y_1,...,y_{n-k}$ for coordinates on $\mathbb{R}^{n-k}$ the shrinking cylinder may be parametrized by
$$F(x_1,...x_{k+1},y_1,...y_{n-k},t) = (\sqrt{-t} x_1,...\sqrt{-t}x_{k+1},y_1,..., y_{n-k}),$$
and
\begin{align*}
\tilde{F}(x_1,...x_{k+1},y_1,...y_{n-k},t) &= \sqrt{-t} \tilde{F}(x_1,...x_{k+1},y_1,...y_{n-k},-1) \\
&= \sqrt{-t} (x_1,...x_{k+1},y_1,...y_{n-k},).
\end{align*}
Here, \(\varphi_t^{-1}\) acts as the identity on \(\mathbb{S}^k\) and scales the \(\mathbb{R}^{n-k}\)-factor by \(\sqrt{-t}\). We denote the metric
\[
g_t := (-t)\,\varphi_t^{*}g_{\mathscr{C}}.
\]
Observe that \(g_t\) is uniformly equivalent to \(g_{-1}\) for \(t \in [-2,-\tfrac{1}{2}]\). Combined with the fact that the coefficients of the transformed equation are uniformly controlled, standard interior parabolic Schauder estimates apply. This yields the following Schauder estimates on the truncated cylinder.

\begin{lemma} {\label{lemma 4.3}}
    Let $\ell \geq 0.$ There exists a constant $C_\ell$ depending on $\ell$ such that whenever $f \in C^{\ell,\alpha}(\mathscr{C})$ and $u \in C^{2+\ell,\alpha}_{loc}(\mathscr{C})$ satisfies $Lu = f - P_{\ker L} f$, then 
    $$\|u\|_{C^{2+\ell,\alpha}(\mathscr{C}\cap B_R(0))} \leq C_\ell \big(R^2\|f\|_{C^{\ell,\alpha}(\mathscr{C}) } + \|u\|_{C^0(\mathscr{C} \cap B_{4R}(0))}\big).$$
\end{lemma}

\begin{proof}
Denote $g$ by $f - P_{\ker L} f.$ Define the new functions as follows,   
$$\Tilde{u}(x,t) := \sqrt{-t} \big(u(\varphi_t(x))\big), \hspace{0.5 cm} \Tilde{g}(x,t) := -\frac{1}{\sqrt{-t}} \big(g(\varphi_t(x)) \big).$$
Then
$$\partial_t \Tilde{u} = \frac{1}{2\sqrt{-t}} (x. \nabla u) (\varphi_t) - \frac{1}{2\sqrt{-t}} u(\varphi_t).$$
This implies,
\begin{align*}
    \partial_t \Tilde{u} - \Delta_{g_t} \Tilde{u} &=  \frac{1}{2\sqrt{-t}} \big( (x. \nabla_{\mathscr{C}} u) (\varphi_t) - u(\varphi_t) \big) - \frac{1}{\sqrt{-t}} \Delta_{\mathscr{C}} u (\varphi_t) \\
    &= - \frac{1}{\sqrt{-t}} \big(\Delta_{\mathscr{C}} u (\varphi_t) - \frac{1}{2} x. \nabla_{\mathscr{C}} u (\phi_t) +  u(\phi_t) \big) + \frac{1}{2\sqrt{-t}} u(\phi(t))\\
    &= \Tilde{g} - \frac{1}{2t} \Tilde{u}.
\end{align*}
Now, using interior Parabolic Schauder type estimates (see, for instance, Brandt \cite{brandt1969interior}) for a fixed $r > \sqrt{2k}$ and any point $\mathscr{C}$ of the form $x = (0,...,0,y_1,...,y_{n-k}),$ we obtain
$$\sup \limits_{t \in [-1,-\frac{3}{4}]} \|\Tilde{u} \|_{C^{2+\ell,\alpha} (\mathscr{C} \cap B_r(x))} \leq C \sup \limits_{t \in [-\frac{3}{2},-\frac{1}{2}]} \big( \|\Tilde{g}\|_{C^{\ell,\alpha}(\mathscr{C} \cap B_{\frac{3}{2}r}(x))} + \|\Tilde{u}\|_{C^0(\mathscr{C} \cap B_{\frac{3}{2}r}(x))}\big).$$ 
Note that the constant $C_\ell$ depends on $\ell$ but is independent of the point $x$ because of the same geometry. Now at $t= -1$, 
$$\|\tilde{u} (\cdot,-1)\|_{C^{2+\ell,\alpha} (\mathscr{C} \cap B_r(x))} \geq C_\ell^{-1} \|u \|_{C^{2+\ell,\alpha} (\mathscr{C} \cap B_r(x))}.$$ 
Also, by the description of $\phi_t$ for each $t \in [-\frac{3}{2},-\frac{1}{2}]$, we have 
$$\|\Tilde{g}(\cdot, t)\|_{C^{\ell,\alpha}(\mathscr{C} \cap B_{\frac{3}{2}r}(x))} \leq C \|g\|_{C^{\ell,\alpha}(\mathscr{C} \cap \varphi_t (B_{\frac{3}{2}r}(x)))},$$ 
and
$$\|\Tilde{u}(\cdot, t)\|_{C^\ell(\mathscr{C} \cap B_{\frac{3}{2}r}(x))} \leq C_\ell \|u\|_{C^\ell(\mathscr{C} \cap \varphi_t (B_{\frac{3}{2}r}(x)))}.$$ 
For $R\geq r,$ with $|x| \leq R$, one can check that the diffeomorphism $\varphi$ satisfies 
$$\phi_t(B_{\frac{3}{2}r}(x)) \subseteq B_{3R}(x).$$ 
Therefore,
$$\|u\|_{C^{2+\ell,\alpha}(\mathscr{C} \cap B_r(x))} \leq C_\ell \big( \|g\|_{C^{\ell,\alpha}(\mathscr{C} \cap B_{3R}(x))} + \|u\|_{C^0(\mathscr{C} \cap B_{3R}(x))}\big).$$
Here, $C_\ell$ is independent of the choice of $x$. Indeed, all such balls $B_r(x)$, with $x$ on the axis of the cylinder, have identical geometry due to the translation invariance of $\mathscr{C}$ in the $\mathbb{R}^{n-k}$ directions. Hence, the interior Schauder constants are uniform over this family. 

We cover $\mathscr{C} \cap B_R(0)$ by balls $B_r(x),$ where the centers are chosen on the axis of the cylinder and $|x| \leq R$. Since $r > \sqrt{2k}$ every point $p \in \mathscr{C} \cap B_R(0)$ is contained in one of the balls. For such a center $x$, we have $|x| \leq R,$ and hence 
$$B_{3R}(x) \subseteq B_{4R}(0).$$ 
Consequently,
$$\|u\|_{C^{2+\ell,\alpha}(\mathscr{C}\cap B_r(x))} \leq C_\ell \big(\|g\|_{C^{\ell,\alpha}(\mathscr{C} \cap B_{4R}(0))) } + \|u\|_{C^0(\mathscr{C} \cap B_{4R}(0))}\big).$$
Now, let $p,q \in \mathscr{C} \cap B_R(0)$. Clearly $p \in B_r(x)$ for some $x = (0,...,0,y_1,...,y_{n-k})$. The local estimate obtained above therefore implies that, for every derivative  $\nabla^j u$ with $0 \leq j \leq \ell+2,$
$$|\nabla^j u(p)| \leq C_\ell \big(\|g\|_{C^{\ell,\alpha}(\mathscr{C} \cap B_{4R}(0))) } + \|u\|_{C^0(\mathscr{C} \cap B_{4R}(0))}\big).$$
If $p$  and $q$ belong to a common ball $B_r(x)$ then the Holder estimate from the local $C^{2+\ell, \alpha}$ gives
$$\frac{|\nabla^{\ell+2} u(p) - \nabla^{\ell+2} u(q)|}{|p-q|^{\alpha}} \leq C_\ell \big(\|g\|_{C^{\ell,\alpha}(\mathscr{C} \cap B_{4R}(0))) } + \|u\|_{C^0(\mathscr{C} \cap B_{4R}(0))}\big).$$
If $p$ and $q$ do not belong to a common ball $B_r(x)$. Since $r$ is fixed and the covering by balls $B_r(x)$ has uniform geometry, there exists $c_r >0$ such that any two points with distance less than $c_r$ lie in a common ball $B_r(x).$ Therefore $d(p,q) \geq c_r$. Hence,
\begin{align*}
    \frac{|\nabla^{\ell+2} u(p) - \nabla^{\ell+2} u(q)|}{|p-q|^{\alpha}} &\leq c_r^{-\alpha}|\nabla^{\ell+2} u(p) - \nabla^{\ell+2} u(q)| \\
    &\leq c_r^{-\alpha} (|\nabla^{\ell+2} u(p)| + |\nabla^{\ell+2} u(q)|) \\
    &\leq C_\ell \big(\|g\|_{C^{\ell,\alpha}(\mathscr{C} \cap B_{4R}(0))) } + \|u\|_{C^0(\mathscr{C} \cap B_{4R}(0))}\big).
\end{align*}
after enlarging $C_\ell$ if necessary. Finally, Lemma ~\ref{lemma 4.1} implies
$$\|g\|_{C^{\ell,\alpha}(\mathscr{C} \cap B_{4R}(0))} \leq CR^2 \|f\|_{C^{\ell,\alpha}(\mathscr{C})}.$$ 
Therefore,
$$\|u\|_{C^{2+\ell,\alpha}(\mathscr{C}\cap B_R(0))} \leq C_\ell \big(R^2\|f\|_{C^{\ell,\alpha}(\mathscr{C}) } + \|u\|_{C^0(\mathscr{C} \cap B_{4R}(0))}\big),$$
which proves the lemma.
\end{proof}
We now introduce the following spaces on the truncated cylinder, $\mathscr{C} \cap B_R(0)$. Define
\[
Z_{R,\ell}
:=
\left\{
f \in C^{\ell,\alpha}(\mathscr{C} \cap B_R(0))
:
E_R(f) \in C^{\ell,\alpha}(\mathscr{C})
\right\},
\]
where $E_R(f)$ denotes the extension of $f$ by zero outside $\Omega_R$, i.e.
\[
E_R(f)(x) =
\begin{cases}
f(x), & x \in \mathscr{C} \cap B_R(0),\\
0, & x \in \mathscr{C} \setminus B_R(0).
\end{cases}
\]
We equip $Z_{R,\ell}$ with the norm
\[
\|f\|_{Z_{R,\ell}} := \|f\|_{C^{\ell,\alpha}(\mathscr{C} \cap B_R(0))}.
\]
By Lemma~3.11 of \cite{colding2015uniqueness}, the operator
\[
L : (\ker L)^\perp \cap W^{2,2}(\mathscr{C})
\longrightarrow
(\ker L)^\perp \cap L^2(\mathscr{C})
\]
is an isomorphism. Thus, we may define $L^{-1}$ to be its inverse.

For $f \in Z_{R,\ell}$, we define $L_R^{-1} f$ as follows, extend $f$ by zero to $\mathscr{C}$, project the extension onto $(\ker L)^\perp$, solve the equation $L u = P_{(\ker L)^\perp}(E_R(f))$ on $\mathscr{C}$ by setting
\[
u = L^{-1}(P_{(\ker L)^\perp}(E_R(f))),
\]
and then define
\[
L_R^{-1} f := u|_{\mathscr{C} \cap B_R(0)}.
\]
For $R > 0,$ we define 
$$L_R^{-1}(u) := (L^{-1} u)\big|_{\mathscr{C} \cap B_R(0)}.$$
Let \(f \in Z_{R,\ell}\) and $u = L_R^{-1} f.$ If $R \geq R_0$, where $R_0$ is the constant from Lemma~\ref{lemma 4.2}, then Lemma~\ref{lemma 4.2} and Lemma ~\ref{lemma 4.3} imply that
\begin{align*}
\|u\|_{C^{2+\ell,\alpha}(\mathscr{C} \cap B_R(0))}
&\le C_\ell \left(R^2\|f\|_{C^{\ell,\alpha}(\mathscr{C})} + \|u\|_{C^0(\mathscr{C} \cap B_{4R}(0))} \right) \\
&\le C_\ell R^4\left(\|f\|_{C^{\ell,\alpha}(\mathscr{C})} + \|u\|_{C^0(\mathscr{C} \cap B_{R_0}(0))} \right).
\end{align*}
We now prove an implicit function theorem in the spirit of \cite{adams1988rates}, which will be crucial for the construction of our comparison surfaces.

We begin with a cutoff function. Let \(\chi : \mathbb{R} \to [0,1]\) be a smooth function such that
\[
\chi(t) = 1 \quad \text{for } t \le 0,
\qquad
\chi(t) = 0 \quad \text{for } t \ge 1.
\]
For \(R>0\), define the radial cutoff function \(\chi_R : \mathbb{R}^{n+1} \to [0,1]\) by
\[
\chi_R(x) := \chi(|x| - (R-1)).
\]

\begin{prop}\label{prop 4.4}
There exists a constant $R_0>0$ such that for every $\ell\ge0$ there
exists a constant $C_\ell>0$ with the following property. Whenever
$R\ge R_0$ and $\delta>0$ satisfy $\delta R^6 \leq C_\ell^{-1}$, 
there exists a map
\[
N_R:\ker L \cap \{v:\|v\|_{L^2}\le \delta\}
\longrightarrow
C^{2 + \ell,\alpha}(\mathscr C\cap B_R(0))
\]
satisfying the following properties:
\begin{enumerate}[label=\normalfont(\roman*)]
    \item We have
    \[
    \|N_R(u)\|_{C^{2+\ell,\alpha}(\mathscr C\cap B_R(0))}
    \le C_\ell R^6 \|u\|_{L^2}^2.
    \]

    \item For every $u\in \ker L\cap\{v:\|v\|_{L^2}\le \delta\}$,
    \[
    P_{(\ker L)^\perp}
    Q\bigl(\chi_R(u+N_R(u))\bigr)= L(N_R(u))
    \]
    on $\mathscr{C} \cap B_R(0).$

    \item The map $N_R$ is locally Lipschitz. Moreover, for every $u_1,u_2$ in its domain,
    \[
    \|N_R(u_1)-N_R(u_2)\|_{C^{2+\ell,\alpha}(\mathscr C\cap B_R(0))}
    \le
    C_\ell R^6
    \bigl(\|u_1\|_{L^2}+\|u_2\|_{L^2}\bigr)
    \|u_1-u_2\|_{L^2}.
    \]

    \item Let $R_2 \geq R_1 \geq R_0,$ and $\delta R_2^6 \leq C_\ell^{-1}.$ Then  
    \[
    \begin{aligned}
    \|N_{R_1}(u_1)-N_{R_2}(u_2)\|_{H^2}
    \le\;&
    C_\ell R_1^2
    \bigl(\|u_1\|_{L^2}+\|u_2\|_{L^2}\bigr)
    \|u_1-u_2\|_{L^2}
    \\
    &+
    C_\ell
    \bigl\|
    (\chi_{R_1}-\chi_{R_2})
    \bigl(u_2+N_{R_2}(u_2)\bigr)
    \bigr\|_{L^2}.
    \end{aligned}
    \]
\end{enumerate}
\end{prop}

\begin{proof}
Recall the definition of $Z_{R,\ell}$. Define
$$X: = \{f \in \ker L: \|f\|_{L^2} \leq \delta \},$$
and
$$Y_{R,\ell}:= \big \{g \in C^{2+\ell,\alpha}(\mathscr{C} \cap B_R(0)) :\; \|g\|_{C^{2+\ell,\alpha}(\mathscr{C} \cap B_R(0))} \leq \delta\big \}.$$
Consider the map
\begin{align*}
    F_R: X \times Y_{R,\ell} &\longrightarrow Z_{R,\ell} \\
    (u,v) &\longmapsto Q \big(\chi_R (u +v)\big)|_{B_R(0)} 
\end{align*}
We have
\begin{align*}
    F_R(u,v_1) - F_R(u,v_2) &= Q (\chi_R(u+v_1)) - Q(\chi_R(u+v_2)) \\
    &= \int_0^1 DQ\!\left(\chi_R\big(u+t v_1+(1-t)v_2\big)\right) \big(\chi_R(v_1-v_2)\big)\,dt.
\end{align*}
It therefore follows from \cite[Proposition~A.1]{sun2022generic} (see also \cite[Lemma~4.10]{colding2015uniqueness}) that
$$\|F_R(u,v_1) - F_R(u,v_2) \|_{C^{\ell,\alpha} (\mathscr{C})} \leq C_\ell \delta R^2 \| v_1 - v_2\|_{C^{2,\alpha}}.$$
Let $G$ be defined by 
$$G(u,v) := L_R^{-1} F_R(u,v).$$ 
We have shown that
\begin{align*}
    \|G(u,v_1) - G(u,v_2)\|_{C^{2+\ell,\alpha} (\mathscr{C} \cap B_R(0))} \leq C_\ell R^4 \big(&\;\|F_R(u,v_1) - F_R(u,v_2)\|_{C^{\ell,\alpha}(\mathscr{C})} \\
    &+ \|G(u,v_1) - G(u,v_2)\|_{C^0(\mathscr{C} \cap B_{R_0}(0))}\big). 
\end{align*}
Observe that on $B_{R_0}(0)$, 
$$L (G(u,v_1) - G(u,v_2)) = P_{(\ker L)^\perp} \big(F_R(u,v_1) - F_R(u,v_2) \big).$$
So by the standard elliptic Schauder estimate and Lemma ~\ref{lemma 5.1},
\begin{align*}
    \|G(u,v_1) - G(u,v_2)\|_{C^{2+\ell,\alpha}(\mathscr{C} \cap B_{R_0}(0))} \leq& C_\ell \big(R_0^2 \;\|F_R(u,v_1) - F_R(u,v_2)\|_{C^{\ell,\alpha}(\mathscr{C} \cap B_{R_0}(0))} \\
    &+\|G(u,v_1) - G(u,v_2)\|_{\widehat{L}^2(\mathscr{C} \cap B_{R_0}(0))}  \big).
\end{align*}
where in the above inequality $\widehat{L}^2$ is the standard $L^2$ norm, but
$$\|G(u,v_1) - G(u,v_2)\|_{\widehat{L}^2(\mathscr{C} \cap B_{R_0}(0))} \leq e^{\frac{R_0^2}{4}} \|G(u,v_1) - G(u,v_2)\|_{L^2(\mathscr{C} \cap B_{R_0}(0))}.$$
And,
\begin{align*}
     \|G(u,v_1) - G(u,v_2)\|_{L^2(\mathscr{C} \cap B_{R_0}(0))} &\leq \|L^{-1} \big(P_{(Ker L)^\perp} \big( F_R(u,v) - F_R(u,v_2) \big) \big)\|_{L^2(\mathscr{C} )} \\
     &\leq C \|\big(P_{(Ker L)^\perp} \big(F_R(u,v) - F_R(u,v_2) \big) \big)\|_{L^2(\mathscr{C})} \\
     &\leq C \|F_R(u,v_1) - F_R(u,v_2)\|_{C^0(\mathscr{C} \cap B_R(0))}.
\end{align*}
Hence, we can conclude
$$\|G(u,v_1) - G(u,v_2)\|_{C^{2+\ell,\alpha} (\mathscr{C} \cap B_R(0))} \leq C_\ell \delta R^6 \|v_1 - v_2 \|_{C^{2+\ell,\alpha} (\mathscr{C} \cap B_R(0))}.$$
Now if $C_\ell\delta R^6 \leq \frac{1}{2}$, then this shows that $G(u,v)$ is contracting for $\|v\|_{C^{2+\ell,\alpha}(\mathscr{C} \cap B_R(0))} \leq \delta$ for each $\|u\|_{L^2} \leq \delta.$ 

Next we show that there exists a constant $C'_\ell >0$ such that $G(u,v)$ maps 
$$Y_{R,\ell} \cap \{g: \|g \|_{C^{2+\ell,\alpha} (\mathscr{C} \cap B_R(0))} \leq C'_\ell R^6 \|u\|^2_{L^2}\}$$ 
into itself when $u$ is restricted to $X.$ First note 
$$\|Q(\chi_Ru)\|_{C^{\ell,\alpha}(\mathscr{C} \cap B_R(0))} \leq C_\ell \|\chi_Ru\|^2_{C^{2+\ell,\alpha}(\mathscr{C} \cap B_R(0))}.$$
Now 
$$\|Q(\chi_Ru)\|_{C^{\ell,\alpha}(\mathscr{C} \cap B_R(0))} \leq C_\ell R^2 \|u\|^2_{L^2}.$$
Moreover,
$$G(u,0)= L_R^{-1} F_R(u,0) = L_R^{-1} \big(L(\chi_Ru) - \phi(\chi_Ru)\big).$$
Therefore,
$$\|G(u,0)\|_{C^{2+\ell,\alpha}(\mathscr{C} \cap B_R(0))} \leq C_\ell R^6 \|u\|^2_{L^2}.$$
Consequently,
\begin{align*}
    \|G(u,v) \|_{C^{2+\ell,\alpha}(\mathscr{C} \cap B_R(0))} &\leq  \|G(u,0) \|_{C^{2+\ell,\alpha}(\mathscr{C} \cap B_R(0))} + \frac{1}{2}  \|v\|_{C^{2+\ell,\alpha}(\mathscr{C} \cap B_R(0))} \\
    &\leq C_\ell R^6 \|u\|^2_{L^2}  +  \frac{1}{2} C'_\ell R^6 \|u\|^2_{L^2} .
\end{align*}
Choosing $C'_\ell = 2C_\ell$, we see that $G(u,.)$ maps 
$$Y_{R,\ell} \;\cap\; \{g: \|g \|_{C^{2+\ell,\alpha}(\mathscr{C} \cap B_R(0))} \leq C'_\ell R^6 \|u\|^2_{L^2}\}$$ 
into itself. Therefore, by the contraction mapping principle, for every $u \in X$ there exists a unique fixed point 
$$N_R(u) \in Y_{R,\ell} \;\cap\; \{g: \|g \|_{C^{2+\ell,\alpha}(\mathscr{C} \cap B_R(0))} \leq C'_\ell R^6 \|u\|^2_{L^2} \}$$ 
satisfying
$$G(u,N_R(u)) = N_R(u).$$
This implies  
$$N_R(u) = L_R^{-1} F_R(u,N_R(u)) = L_R^{-1} Q (\chi_R (u + N_R(u))).$$
Applying $L$ and using the definition of $L_R^{-1}$, we obtain
$$L(N_R(u)) = P_{(\ker L)^\perp} Q (\chi_R (u +N_R(u))).$$
We now prove ~(iii). Let
\[
N_i=N_R(u_i),
\qquad i=1,2.
\]
Since each $N_i$ is a fixed point of $G$, we have
\[
N_i=G(u_i,N_i).
\]
Hence,
\[
N_1-N_2
=
G(u_1,N_1)-G(u_2,N_2).
\]
Therefore,
\begin{align*}
\|N_1-N_2\|_{C^{2+\ell,\alpha} (\mathscr{C} \cap B_R(0))}
&\le
\|G(u_1,N_1)-G(u_1,N_2)\|_{C^{2+\ell,\alpha} (\mathscr{C} \cap B_R(0))} \\
&\hspace{0.6 cm}
+
\|G(u_1,N_2)-G(u_2,N_2)\|_{C^{2+\ell,\alpha} (\mathscr{C} \cap B_R(0))}.
\end{align*}
The first term is estimated exactly as before
\[
\|G(u_1,N_1)-G(u_1,N_2)\|_{C^{2+\ell,\alpha} (\mathscr{C} \cap B_R(0))} 
\le
\frac12
\|N_1-N_2\|_{C^{2+\ell,\alpha} (\mathscr{C} \cap B_R(0))} .
\]
For the second term recall,
\[
G(u,v) = L_R^{-1} F_R(u,v).
\]
Again,
\begin{align*}
F_R(u_1, N_2) - F_R(u_2, N_2) &= Q(\chi_R(u_1+N_2)) - Q(\chi_R(u_2+N_2)) \\
&= \int_0^1 DQ \big(\chi_R (t u_1 + (1-t)u_2 + N_2) \big) \big(\chi_R(u_1 - u_2) \big) dt.
\end{align*}
We have shown 
\[
\|N_2\|_{C^{2+\ell,\alpha} (\mathscr{C} \cap B_R(0))}
\le
C_\ell R^6 \|u_2\|_{L^2}^2.
\]
Using
\[
\delta R^6\le C_\ell^{-1},
\]
we obtain
\[
\|F_R(u_1,N_2)-F_R(u_2,N_2)\|_{C^{\ell,\alpha}}
\le
C_\ell R^2
\bigl(\|u_1\|_{L^2}+\|u_2\|_{L^2}\bigr)
\|u_1-u_2\|_{L^2}.
\]
Applying the estimate for $L_R^{-1}$ gives
\[
\|G(u_1,N_2)-G(u_2,N_2)\|_{C^{2+\ell,\alpha}(\mathscr C\cap B_R(0))}
\le
C_\ell R^6
\bigl(\|u_1\|_{L^2}+\|u_2\|_{L^2}\bigr)
\|u_1-u_2\|_{L^2}.
\]
Combining the above inequalities yields
\[
\|N_1-N_2\|_{C^{2+\ell,\alpha}(\mathscr C\cap B_R(0))}
\le
C_\ell R^6
\bigl(\|u_1\|_{L^2}+\|u_2\|_{L^2}\bigr)
\|u_1-u_2\|_{L^2}.
\]
This completes the proof of ~(iii). It remains to prove ~(iv).
Write
\[
N_i=N_{R_i}(u_i), \qquad i=1,2.
\]
By the fixed point equation,
\[
N_i
=
L^{-1}
P_{(\ker L)^\perp}
E_{R_i}
\!\left(
Q(\chi_{R_i}(u_i+N_i))
\right).
\]
Hence,
\[
N_1-N_2
=
L^{-1}P_{(\ker L)^\perp}
\big(
E_{R_1}Q(\chi_{R_1}(u_1+N_1)) -
E_{R_2}Q(\chi_{R_2}(u_2+N_2))
\big).
\]
Since $L^{-1}$ is bounded on $(\ker L)^\perp$, it follows that
\[
\begin{aligned}
\|N_1-N_2\|_{H^2}
&\leq
C
\big\|
E_{R_1}Q(\chi_{R_1}(u_1+N_1))
-
E_{R_2}Q(\chi_{R_2}(u_2+N_2))
\big\|_{L^2}
\\
&\leq
C(\|I_1\|_{L^2}+\|I_2\|_{L^2}),
\end{aligned}
\]
where
\[
\begin{aligned}
I_1& :=
E_{R_1}
\left(
Q(\chi_{R_1}(u_1+N_1))
-
Q(\chi_{R_1}(u_2+N_2))
\right),
\\
I_2& :=
E_{R_1}
Q(\chi_{R_1}(u_2+N_2))
-
E_{R_2}
Q(\chi_{R_2}(u_2+N_2)).
\end{aligned}
\]
By the same argument as in the proof of part ~(iii), we obtain
\[
\|I_1\|_{L^2}
\le
C_\ell R_1^2 \bigl(\|u_1\|_{L^2}+\|u_2\|_{L^2}\bigr)
\|u_1-u_2\|_{L^2}.
\]
For the second term, write
\[
w=u_2+N_2,
\]
we get
\[
Q(\chi_{R_1}w)-Q(\chi_{R_2}w)
=
\int_0^1
DQ(t\chi_{R_1} w + (1-t)\chi_{R_2}w)
(\chi_{R_1} w -\chi_{R_2}w )\,dt,
\]
Therefore,
\[
\|I_2\|_{L^2}
\le
C_\ell
\bigl\|
(\chi_{R_1}-\chi_{R_2})
(u_2+N_2)
\bigr\|_{L^2}.
\]
Combining the above estimates completes the proof.
\end{proof}

Let $\ell \geq 0.$ Recall the definition of $u_\alpha$ from Section \ref{section 2} for $\alpha \in \mathbb{R}^{n+1}$. Fix an initial condition $\alpha(0)$, and choose $R$ so that $$2|\alpha(0)| R^6 \leq C_\ell^{-1},$$ 
where $C_\ell$ is the constant from Proposition ~\ref{prop 4.4}. Applying Proposition ~\ref {prop 4.4} with $\delta = |\alpha(0)|$, we then obtain a map $N_R$ defined for all $u_\alpha$ with 
$$|\alpha| \leq 2|\alpha(0)|.$$ 
For every such $\alpha$ and every $r \leq R$, 
$$\|N_R(u_{\alpha})\|_{C^{2+\ell,\alpha}(\mathscr{C} \cap B_R(0))} \leq C_\ell |\alpha|^2R^6,$$ 
and 
$$P_{(\ker L)^\perp} Q \big(\chi_R (u_\alpha + N_R(u_\alpha))) = L(N_R(u_\alpha)).$$ 
We now seek a solution of the equation
$$u_{\alpha'(\tau)} = P_{\ker L} Q \big(\chi_R (u_{\alpha(\tau)} + N_R(u_{\alpha(\tau)}))).$$
Using the identification of $\ker L$ with $\mathbb{R}^{n+1}$ given by $\alpha \mapsto u_\alpha$, this is equivalent to the finite-dimensional ordinary differential equation
$$\alpha'(\tau) = \mathcal{Q}(\alpha(\tau)),$$
where $\mathcal{Q}$ denotes the vector field induced by 
$$P_{\ker L} Q \big(\chi_R (u_\alpha + N_R(u_\alpha))).$$
The vector field $\mathcal{Q}$ is defined on the ball 
$$\{|\alpha| < 2|\alpha(0)|\}.$$ 
Since $\alpha(0)$ lies in the interior of this ball, standard ODE theory yields a solution on some maximal interval $[0,T)$. Furthermore,
$$|\mathcal{Q}(\alpha(\tau))| \leq C|\alpha(\tau)|^2,$$
and hence
$$|\alpha'(\tau)| \leq C|\alpha(\tau)|^2.$$
It follows that there exists a constant $c_0 >0$, independent of $\alpha(0)$, such that
$$|\alpha(\tau)| \leq 2|\alpha(0)|$$
for all
$$\tau \in [0, c_0 |\alpha(0)|^{-1}].$$
In particular, $\alpha(\tau)$ remains in the domain of definition of $\mathcal{Q}$ throughout this interval. Therefore, by the continuation theorem for ordinary differential equations, the solution extends throughout the interval.

We have therefore shown that, for every $\alpha$ with sufficiently small $|\alpha|$, if $R^6 \leq (2C_\ell|\alpha|)^{-1}$, then there exists a solution $\alpha(\tau)$ with $\alpha(0) = \alpha$ on the interval $\tau \in [0, c_1|\alpha|^{-1}]$  
such that
$$\partial_\tau (u_{\alpha(\tau)}) +  L(N_R (u_{\alpha(\tau)})) = Q \big(\chi_R (u_{\alpha(\tau)} + N_R(u_{\alpha(\tau)}))).$$
Fix a sufficiently small constant $c_1>0$. For each $\alpha$, define $R_\alpha$ by
\[
|\alpha|^{\frac94} e^\frac{R_\alpha^2}{8} = c_1.
\]
Then as $|\alpha|\to0$, 
\[
|\alpha|R_\alpha^6\to0.
\]
Hence, for all sufficiently small $|\alpha|$,
\[
R_\alpha^6\leq (2C_\ell|\alpha|)^{-1},
\]
so, the preceding construction applies with $R=R_\alpha$. Throughout the remainder of the paper, we suppress the dependence on $R_\alpha$ and simply write $N.$  

Our next goal is to construct a rescaled mean curvature flow equation which is graphical over the cylinder on $B_{R_{\alpha}-1}(0)$. The function $u_{\alpha(\tau)} + N u_{\alpha(\tau)}$ will serve as the leading-order approximation to the solution.

We now introduce the function spaces that will be used in the sequel. Let $T> 0$. Throughout we will take $T \leq c_0 |\alpha|^{-1}$, where $c_0 >0$ is the constant from the previous construction. For any integer $\ell >0,$ we use $\mathcal{X}_\ell$ to denote the space of functions $u:[0, T] \to H^{\ell+1}(\mathscr{C})$ such that
$$\|u\|_{\mathcal{X}_\ell} := \bigg(\int_0^T \|u(\tau)\|^2_{H^{\ell+1}(\mathscr{C})} d\tau \bigg)^{1/2} + \sup_{\tau \in [0,T]} \|u(\tau)\|_{H^\ell(\mathscr{C})} < \infty.$$
Let $\mathcal{P}: L^2(\mathscr{C}) \to L^2(\mathscr{C})$ denote the $L^2$-orthogonal projection onto the subspace spanned by eigenfunctions corresponding to eigenvalues $\leq \lambda_4 = \max \{-\frac{1}{2}, -\frac{1}{k }\}.$ We will use the following consequence of Corollary $3.3$ in Strehlke \cite{strehlke2020asymptotics}; see also Lemma ~A.3 in Sun-Wang-Xue \cite{sun2025regularity}. 

\begin{lemma} {\label{lemma 4.5}}
Let $v:[0,T]\to H^\ell$ be a continuously differentiable path satisfying
\[
\mathcal{P}v(s)=v(s)
\]
for all $s\in[0,T]$. Then for $\ell\ge 1$, there exists a constant $C_\ell>0,$ such that
$$\sup_{\tau \in [0,T]}\|v(\tau)\|^2_{H^\ell} + \int_0^T \|v(\tau)\|^2_{H^{\ell+1}}d\tau \leq C_\ell \bigg(\|v(0)\|^2_{H^\ell} + \int_0^T \|(\partial_s - L) v(s)\|^2_{H^{\ell - 1}}\; ds \bigg).$$
\end{lemma}
The proof is identical to that of Corollary $3.3$ in \cite{strehlke2020asymptotics}, so we omit the details.

For  $v \in \mathcal{X}_{\ell}$ with sufficintly small $C^{2, \alpha}(\mathscr{C} \cap B_{R_\alpha}(0))$ norm, we define $U_{\alpha}(v)$ to be the unique solution to the initial value problem 
\begin{equation*}
\begin{cases}
(\partial_\tau - L) U_{\alpha}(v) = F_{\alpha}(v) , \\
U_{\alpha}(v)(0) = -\int_0^T e^{-sL} (1- \mathcal{P}) F_{\alpha}(v)(s) ds. &
\end{cases}
\end{equation*}
where
$$F_{\alpha}(v)(\tau) := Q(\chi_{\alpha} (u_{\alpha(\tau)
} + N(u_{\alpha(\tau)}) + v)) - Q(\chi_{\alpha} (u_{\alpha(\tau)
} + N(u_{\alpha(\tau)}))) - \partial_\tau (\chi_{\alpha} N(u_{\alpha(\tau)})).$$
Here $\alpha(\tau)$ denotes the solution of the ODE constructed previously with initial condition $\alpha(0) = \alpha, \chi_{\alpha} = \chi_{R_{\alpha}}$, and $e^{\tau L}$ denotes the heat semigroup for the operator $L$.

The existence of $U_\alpha(v)$ in $\mathcal{X}_\ell$ follows directly from the semigroup argument. Moreover, $U_\alpha$ admits the explicit representation
$$U_{\alpha}(v)(\tau) = \int_0^{\tau} e^{(\tau-s)L} \mathcal{P} F_{\alpha}(v)(s) ds -\int_\tau^T e^{(\tau-s)L} (1- \mathcal{P}) F_{\alpha}(v)(s) ds.$$
We now show that the solution belongs to $\mathcal{X}_\ell.$ The proof is similar to the proof of Theorem ~3.7 in \cite{strehlke2020asymptotics} (see also Lemma ~A.4 in \cite{sun2025regularity}).

\begin{lemma} {\label{lemma 4.6}}
    Let $\ell >  \left\lfloor \frac{n}{2} \right\rfloor +1$, then there exist $C_\ell, r_\ell> 0$ such that when $v_1, v_2 \in \mathcal{X}_\ell$ and $|\alpha| R_\alpha^2 + |\alpha|^2 R_\alpha^6+ e^{\frac{R_\alpha^2}{8}} (\|v_1\|_{\mathcal{X}_\ell} + \|v_2\|_{\mathcal{X}_\ell}) \leq r_\ell$, we have
    $$\|U_\alpha(v_1) - U_\alpha(v_2) \|_{\mathcal{X}_\ell} \leq C_\ell T (|\alpha| R_\alpha^2 + |\alpha|^2 R_\alpha^6 + e^{\frac{R_\alpha^2}{8}} \|v_1\|_{\mathcal{X}_\ell} + e^{\frac{R_\alpha^2}{8}} \|v_2\|_{\mathcal{X}_\ell}) \|v_1 - v_2\|_{\mathcal{X}_\ell}.$$
\end{lemma}

\begin{proof}
    Let 
    $$D(\tau) = U_\alpha(v_1)(\tau) -  U_\alpha(v_2)(\tau),$$ then $D$ satisfies the equations 
    \begin{equation}
    \begin{cases}
    (\partial_\tau - L) D = F_\alpha(v_1) - F_\alpha(v_2)\\
    D(.,0) = -\int_0^T e^{-sL} (1- \mathcal{P}) (F_\alpha(v_1) - F_\alpha(v_2)) ds. &
    \end{cases}
    \end{equation}
    We have
    $$D(\tau) = \int_0^\tau e^{(\tau-s)L} \mathcal{P} (F_\alpha(v_1) - F_\alpha(v_2)) ds - \int_\tau^T e^{(\tau-s)L} (1-\mathcal{P}) (F_\alpha(v_1) - F_\alpha(v_2)) ds.$$
    When $\tau \geq s, \|e^{(\tau-s)L}\| \leq e^{\lambda_4(\tau-s)}$ on $\mathcal{P}(H^\ell),$ and when $\tau \leq s, \|e^{(\tau-s)L}\| \leq 1$ on $(1-\mathcal{P})(H^\ell).$ Therefore,
    $$(1-\mathcal{P}) D(\tau) = -\int_\tau^T e^{(\tau-s)L} (1-\mathcal{P}) (F_\alpha(v_1) - F_\alpha(v_2)) ds.$$
    Since the $H^{\ell-1}$ and $H^{\ell}$ norms are equivalent on the range of $1 - \mathcal{P}$, we can conclude that
    $$\|(1-\mathcal{P}) D(\tau) \|_{H^\ell} \leq C_\ell \int_\tau^T \|(1-\mathcal{P}) (F_\alpha(v_1) - F_\alpha(v_2))\|_{H^{\ell-1}} ds.$$
    As $1 - \mathcal{P}$ is a finite rank operator, so
    $$\|(1-\mathcal{P}) D(\tau) \|_{H^\ell} \leq C_\ell \int_\tau^T \| (F_\alpha(v_1) - F_\alpha(v_2))\|_{H^{\ell-1}} ds.$$
    Now, by the definition of $F_\alpha,$
    $$F_\alpha(v_1) - F_\alpha(v_2) = Q(\chi_\alpha (u_{\alpha(\tau)} + N(u_{\alpha(\tau)}) + v_1)) - Q(\chi_\alpha (u_{\alpha(\tau)} + N(u_{\alpha(\tau)}) + v_2)).$$
    Define 
    $$A := |\alpha| R_\alpha^2 + |\alpha|^2 R_\alpha^6 + e^{\frac{R_\alpha^2}{8}} (\|v_1\|_{\mathcal{X}_\ell} + \|v_2\|_{\mathcal{X}_\ell}).$$
    Now we can use \cite{strehlke2020asymptotics}, Lemma ~3.5 (see also \cite{sun2025regularity}, Lemma ~A.2), Proposition ~\ref{prop 4.4} and Gaussian Sobolev embedding theorem to conclude that
\begin{multline*}
\left\|
(F_\alpha(v_1)-F_\alpha(v_2))
\right\|_{H^{\ell-1}}
\\
\leq
C_\ell
\Big(
\big(
|\alpha|R_\alpha^2
+
|\alpha|^2R_\alpha^6
+
e^{\frac{R_\alpha^2}{8}}\|v_1\|_{H^{\ell+1}}
+
e^{\frac{R_\alpha^2}{8}}\|v_2\|_{H^{\ell+1}}
\big)
\|v_1-v_2\|_{H^\ell}
\\
+
\big(
|\alpha|R_\alpha^2
+
|\alpha|^2R_\alpha^6
+
e^{\frac{R_\alpha^2}{8}}\|v_1\|_{H^\ell}
+
e^{\frac{R_\alpha^2}{8}}\|v_2\|_{H^\ell}
\big)
\|v_1-v_2\|_{H^{\ell+1}}
\Big).
\end{multline*}
Therefore, it follows for all $\tau \in [0,T],$ we have
$$\|(1-\mathcal{P})D(\tau)\|_{H^\ell} \leq C_\ell A \|v_1 - v_2\|_{\mathcal{X}_\ell}.$$
Now we need to bound $\int_0^T \|(1-\mathcal P)D(\tau)\|_{H^{\ell+1}}^2 d\tau.$ Because $(1- \mathcal{P})$ has finite rank, so we have 
$$\|(1-\mathcal{P})D(\tau)\|_{H^{\ell+1}}^2 \leq C_\ell \|(1-\mathcal{P})D(\tau)\|_{H^{\ell}}^2 $$
for every $\tau \in [0,T].$
Now we can integrate from $0$ to $T$ to get
$$\int_0^T \|(1-\mathcal{P})D(\tau)\|_{H^{\ell+1}}^2  \leq C_\ell T \sup_{\tau \in [0,T]} \|(1-\mathcal{P})D(\tau)\|_{H^{\ell}}^2.$$
Then we can use the pointwise bound $\|(1-\mathcal{P})D(\tau)\|_{H^{\ell}}^2$ to conclude
$$\int_0^T \|(1-\mathcal{P})D(\tau)\|_{H^{\ell+1}}^2  \leq C_\ell T A \|v_1 - v_2\|_{\mathcal{X}_\ell}.$$
Combining both we get 
$$\|(1-\mathcal{P})D\|_{\mathcal{X}_\ell} \leq C_\ell A \|v_1- v_2\|_{\mathcal{\chi}_\ell}.$$
We will now bound $\|.\|_{\mathcal{X}_\ell}$ norm of $\mathcal{P} D(\tau).$ Note that $\mathcal{P} D(0) = 0.$ So using Lemma ~\ref{lemma 4.5},
\begin{align*}
    \sup_{\tau \in [0,T]}\|\mathcal{P}D(\tau)\|^2_{H^\ell} + \int_0^T \|\mathcal{P} D(\tau)\|^2_{H^{\ell+1}}d\tau &\leq C_\ell \int_0^T \|\mathcal{P}(F_\alpha(v_1) - F_\alpha(v_2))\|^2_{H^{\ell - 1}}\; ds. \\
    &\leq C_\ell \int_0^T \|F_\alpha(v_1) - F_\alpha(v_2)\|^2_{H^{\ell - 1}}\; ds.
\end{align*} 
Arguing as above, we conclude that
$$\|\mathcal{P}D\|_{\mathcal{X}_\ell} \leq C_\ell A \|v_1- v_2\|_{\mathcal{\chi}_\ell}.$$
Combining the estimates for the $\mathcal{P}$ and $1- \mathcal{P}$ components yields the desired bound in the statement.
\end{proof}

We now construct the rescaled mean curvature flow.. Let 
$$T = T_\alpha := |\alpha|^{-\frac{1}{8}}.$$  
Note that for all sufficiently small $|\alpha|,$
$$T_\alpha \leq c_1 |\alpha|^{-1}.$$
Moreover, from the differential inequality
\[
|\alpha'(\tau)| \le C|\alpha(\tau)|^2,
\]
we obtain
\[
|\alpha(\tau) - \alpha(0)| \le C|\alpha(0)|^2 \tau.
\]
Since \(T_\alpha = |\alpha|^{-\frac{1}{8}}\), it follows that for all sufficiently small \(|\alpha|\) and all \(\tau \in [0,T_\alpha]\),
\[
|\alpha(\tau) - \alpha(0)| \le C|\alpha|^{\frac{15}{8}}.
\]
Hence,
\[
|\alpha(\tau)| \ge \frac{1}{2}|\alpha(0)|
\]
for all \(\tau \in [0,T_\alpha]\), provided \(|\alpha|\) is sufficiently small.
Fix an integer $\ell = \ell(n)$ sufficiently large so that the Sobolev embedding 
$$\|\chi_R v\|_{C^{2,\alpha}} \leq C_\ell e^{\frac{R^2}{8}} \|v\|_{H^\ell} $$
holds for every $v \in \mathcal{X}_\ell$. Define
$$C(r) := \{ v\in \mathcal{X}_\ell: \|v\|_{\mathcal{X}_\ell} \leq r\}.$$
Define the map 
$$\mathcal{F}_\alpha: C(|\alpha|^\frac{5}{2}) \to \mathcal{X}_\ell$$ 
by 
$$\mathcal{F}_\alpha(v) = U_{\alpha}(v),$$
where $\alpha(0) = \alpha$ and $U_{\alpha}(v)$ solves the initial value problem
\begin{equation*}
\begin{cases}
(\partial_\tau - L) U_{\alpha}(v) = F_{\alpha}(v) , \\
U_{\alpha}(v)(0) = -\int_0^{T_\alpha} e^{-sL} (1- \mathcal{P}) F_{\alpha}(v)(s) ds. &
\end{cases}
\end{equation*}
By definition,
    $$F_\alpha(v)(\tau) = Q(\chi_\alpha (u_{\alpha(\tau)} + N(u_{\alpha(\tau)}) + v)) - Q(\chi_\alpha (u_{\alpha(\tau)} + N(u_{\alpha(\tau)}))) - \partial_\tau(\chi_\alpha N(u_{\alpha(\tau)}).$$
Since, 
$$\|\chi_{R_\alpha} v\|_{C^{2,\alpha}} \leq C_\ell e^{\frac{R_\alpha^2}{8}} \|v\|_{H^\ell},$$ 
for every $v \in C(|\alpha|^{\frac{5}{2}}),$
$$\|\chi_{R_\alpha} v\|_{C^{2,\alpha}} \leq C_\ell c_1 |\alpha|^{\frac{1}{4}}.$$
Hence, after choosing $|\alpha|$ sufficiently small, $F_\alpha(v)$ is well defined. Setting $v =0,$ the two nonlinear terms cancel, and hence 
$$F_\alpha(0)(\tau) = -\partial_\tau(\chi_\alpha N(u_{\alpha(\tau)})).$$
Since $\chi_\alpha$ depends only on $\alpha,$ it is independent of $\tau.$ Therefore, 
$$F_\alpha(0)(\tau) = -\chi_\alpha \partial_\tau(N(u_{\alpha(\tau)})) = -\chi_\alpha DN(u_{\alpha(\tau)})(u_{\alpha'(\tau)})$$
Applying Lemma ~3.5 of \cite{strehlke2020asymptotics} (see also \cite{sun2025regularity}, Lemma ~A.2) together with Proposition ~\ref {prop 4.4}, we obtain
$$\|\mathcal{F}_\alpha(0)\|_{\mathcal{X}_\ell} \leq C_\ell T_\alpha |\alpha|^3  R_\alpha^6 \leq C_\ell |\alpha|^{\frac{23}{8}} R_\alpha^6.$$
By the triangle inequality,
$$\|\mathcal{F}_{\alpha}(v)\|_{\mathcal{X}_\ell} \leq \|\mathcal{F}_\alpha(0)\|_{\mathcal{X}_\ell} + \|\mathcal{F}_{\alpha}(v) - \mathcal{F}_\alpha(0)\|_{\mathcal{X}_\ell} $$
Applying Lemma ~\ref{lemma 4.6}, we obtain
$$\|\mathcal{F}_\alpha(v) - \mathcal{F}_\alpha(0)\|_{\mathcal{X}_\ell} \leq C_\ell T_\alpha (|\alpha| R_\alpha^2 + |\alpha|^2 R_\alpha^6 + e^{\frac{R_\alpha^2}{8}} \|v\|_{\mathcal{X}_\ell}) \|v\|_{\mathcal{X}_\ell}.$$
Since $v \in C(|\alpha|^\frac{5}{2})$,
$$\|\mathcal{F}_{\alpha}(v)\|_{\mathcal{X}_\ell} \leq C_\ell |\alpha|^{\frac{23}{8}} R_\alpha^8 + C_\ell |\alpha|^{\frac{27}{8}}R_\alpha^2 + C_\ell|\alpha|^\frac{35}{8} R_\alpha^6 + C_\ell |\alpha|^{-\frac{1}{8}} e^{\frac{R_\alpha^2}{8}} |\alpha|^{5}.$$
Using the identity
$$e^{\frac{R_\alpha^2}{8}} = c_1 |\alpha|^{-\frac{9}{4}},$$
the last term becomes $C_\ell c_1|\alpha|^{\frac{21}{8}}$. Therefore,
$$\|\mathcal{F}_{\alpha}(v)\|_{\mathcal{X}_\ell} \leq C_\ell |\alpha|^{\frac{23}{8}} R_\alpha^8 + C_\ell |\alpha|^{\frac{27}{8}}R_\alpha^2 + C_\ell|\alpha|^\frac{35}{8} R_\alpha^6 + C_\ell c_1 |\alpha|^{\frac{21}{8}}.$$
Since each exponent is strictly larger than $\frac{5}{2}$, it follows that
$$\|\mathcal{F}_{\alpha}(v)\|_{\mathcal{X}_\ell} \leq |\alpha|^{\frac{5}{2}}$$
for all sufficiently small $|\alpha|.$ Thus,
$$\mathcal{F}_\alpha(C(|\alpha|^\frac{5}{2})) \subseteq C(|\alpha|^\frac{5}{2}).$$
Next, applying Lemma ~\ref{lemma 4.6} once more, we obtain
\begin{multline*}
\|\mathcal{F}_\alpha(v_1)-\mathcal{F}_\alpha(v_2)\|_{\mathcal{X}_\ell}
\\
\leq
C_\ell T_\alpha
\Big(
|\alpha|R_\alpha^2
+
|\alpha|^2R_\alpha^6
+
e^{\frac{R_\alpha^2}{8}}\|v_1\|_{\mathcal{X}_\ell}
+
e^{\frac{R_\alpha^2}{8}}\|v_2\|_{\mathcal{X}_\ell}
\Big)
\|v_1-v_2\|_{\mathcal{X}_\ell}.
\end{multline*}
Since $v_1, v_2 \in C(|\alpha|^\frac{5}{2}),$
$$T_\alpha e^{\frac{R_\alpha^2}{8}}\|v_i\|_{\mathcal{X}_\ell}
\leq
|\alpha|^{-\frac{1}{8}} c_1 |\alpha|^{-\frac{9}{4}} |\alpha|^{\frac{5}{2}}
=
c_1 |\alpha|^{\frac{1}{8}}.$$
for $i = 1,2.$ Moreover,
$$T_\alpha (|\alpha| R_\alpha^2 + |\alpha|^2 R_\alpha^6 )
\leq
|\alpha|^{\frac{7}{8}} R_\alpha^2
+
|\alpha|^{\frac{15}{8}} R_\alpha^6.$$
Hence, for all sufficiently small $|\alpha|,$
$$C_\ell T_\alpha (|\alpha| R_\alpha^2 + |\alpha|^2 R_\alpha^6 + e^{\frac{R_\alpha^2}{8}} \|v_1\|_{\mathcal{X}_\ell} + e^{\frac{R_\alpha^2}{8}} \|v_2\|_{\mathcal{X}_\ell}) \leq \frac{1}{2}.$$
Therefore, for all sufficiently small $|\alpha|, \mathcal{F}_\alpha$ is a contraction on $C(|\alpha|^\frac{5}{2})$. By the Banach fixed-point theorem there exists a unique fixed point $v_\alpha = v_\alpha(\tau),$ such that 
$$\mathcal{F}_{\alpha}(v_\alpha) = v_\alpha.$$
Consequently,
$$(\partial_\tau - L) v_\alpha(\tau) = Q(\chi_\alpha(W_\alpha(\tau) + v_\alpha(\tau))) - Q(\chi_\alpha W_\alpha(\tau)) - \chi_\alpha \partial_\tau(N(u_{\alpha(\tau)})),$$
where 
$$W_\alpha(\tau) = u_{\alpha(\tau)} + N(u_{\alpha(\tau)}).$$
Restricting to $\mathscr{C} \cap B_{R_{\alpha}-1}(0),$ where $\chi_\alpha \equiv 1,$ we obtain
$$(\partial_\tau - L) v_\alpha(\tau) = Q(W_\alpha(\tau) + v_\alpha(\tau)) - Q(W_\alpha(\tau)) - \partial_\tau(N(u_{\alpha(\tau)})).$$
On the other hand, restricting the equation satisfied by $\alpha(\tau)$ to $\mathscr{C} \cap B_{R_{\alpha}-1}(0),$ where $\chi_R \equiv 1,$ yields 
$$(\partial_\tau -L) W_\alpha(\tau) = Q(W_\alpha(\tau)) - \partial_\tau(N(u_{\alpha(\tau)})),$$
Adding these two equations gives 
$$(\partial_\tau -L) (W_\alpha(\tau) + v_\alpha(\tau)) = Q(W_\alpha(\tau) + v_\alpha(\tau)).$$
Therefore, the graph of
\[
u_{\alpha(\tau)} + N(u_{\alpha(\tau)}) + v_\alpha(\tau), \qquad \alpha(0)=\alpha,
\]
over the truncated cylinder \(\mathscr{C}\cap B_{R_\alpha-1}(0)\), where \(R_\alpha\) is defined by
\[
|\alpha|^{\frac94}e^{\frac{R_\alpha^2}{8}}=c_1,
\]
evolves by the rescaled mean curvature flow on the time interval
\[
[0,T_\alpha], \qquad T_\alpha=|\alpha|^{-\frac18}.
\]
We denote this flow by $\mathscr{T}^{\alpha,\mathscr{C}}.$
We also denote by $\mathcal{W}^{\alpha,\mathscr{C}}$
the graph of \(u_\alpha+N(u_\alpha)\) over the truncated cylinder
\(\mathscr{C}\cap B_{R_\alpha-1}(0)\).

Let \(|\alpha|\le |\widehat{\alpha}|\), and denote
\[
\alpha(0)=\alpha, \qquad \widehat{\alpha}(0)=\widehat{\alpha}.
\]
By Proposition~\ref{prop 4.4}(iv),
\[
\|N(u_{\alpha(\tau)})-N(u_{\widehat{\alpha}(\tau)})\|_{H^2}
\leq
C\left(|\alpha-\widehat{\alpha}|+e^{-\frac{(R_{\alpha}-1)^2}{8p_0}}\right),
\]
where \(p_0 >1 \) is chosen in Proposition~\ref{prop 5.2}. Moreover, arguing exactly as in the proof of Proposition~\ref{prop 4.4}(iv) and applying the argument to the fixed-point equations
\[
F_{\alpha}(v_{\alpha})=v_{\alpha},
\qquad
F_{\widehat{\alpha}}(v_{\widehat{\alpha}})=v_{\widehat{\alpha}},
\]
we obtain
\[
\|v_{\alpha}(\tau)-v_{\widehat{\alpha}}(\tau)\|_{H^2}
\leq
C\left(|\alpha-\widehat{\alpha}|+e^{-\frac{(R_{\alpha}-1)^2}{8p_0}}\right).
\]
Consequently, by the triangle inequality and after interchanging \(\alpha\) and \(\widehat{\alpha}\) if necessary, for any \(\alpha,\widehat{\alpha}\) we have
\[
\begin{aligned}
&\|(u_{\alpha(\tau)}+N(u_{\alpha(\tau)})+v_{\alpha}(\tau))
-(u_{\widehat{\alpha}(\tau)}+N(u_{\widehat{\alpha}(\tau)})+v_{\widehat{\alpha}}(\tau))\|_{H^2}\\
&\qquad\leq
C\left(
|\alpha-\widehat{\alpha}|
+e^{-\frac{(R_{\alpha}-1)^2}{8p_0}}
+e^{-\frac{(R_{\widehat{\alpha}}-1)^2}{8p_0}}
\right).
\end{aligned}
\]

The next proposition shows that the Taylor expansion of $\phi$ exhibits a nontrivial second-order obstruction. The proof of the next proposition follows the argument of \cite[Proposition ~4.8]{zhu2020ojasiewicz} with minor modifications. This estimate will play a key role in proving that the functional $\phi$ is bounded from below along the constructed flow.

\begin{prop} {\label{prop 4.7}}
    There exist constants $\kappa$ and $R_1 >0$ such that, for every $R \geq R_1,$ and every $u \in K_2$,
    $$\|P_{\ker L} D^2\phi(\chi_R u, \chi_R u) \|_{L^2} \geq \kappa \|u\|^2_{L^2}.$$
\end{prop}

\begin{proof}
Since $D^2\phi$ is bilinear, it is enough to prove the estimate under the normalization $\|u\|_{L^2}=1$.

Let $u \in K_2$ satisfies $\|u\|_{L^2}=1$. By \cite[Proposition 4.8]{zhu2020ojasiewicz}, there exists a constant $\kappa>0$ such that
\[
\|P_{\ker L}D^2\phi(u,u)\|_{L^2}\ge \kappa
\]
Since $\chi_R\equiv1$ on $B_{R-1}(0)$, we have
\[
D^2\phi(\chi_Ru,\chi_Ru)=D^2\phi(u,u)
\]
on $\mathscr C\cap B_{R-1}(0)$. Hence the difference
\[
D^2\phi(\chi_Ru,\chi_Ru)-D^2\phi(u,u)
\]
vanishes on $\mathscr C\cap B_{R-1}(0).$ Moreover,
\[
|D^2\phi(\chi_Ru,\chi_Ru)|,\;
|D^2\phi(u,u)|
\le C(|u|+|\nabla u| + |\nabla^2 u|)^2.
\]
Since $u\in K_2$ and $\|u\|_{L^2}=1$, 
\[
|D^2\phi(\chi_Ru,\chi_Ru)|,\;
|D^2\phi(u,u)|
\le C(1+|x|^2)^2.
\]
Therefore, for some $p >1,$ by choosing $R_1$ sufficiently large, we may assume that for every $R\ge R_1$,
\[
\|D^2\phi(\chi_Ru,\chi_Ru)\|_{L^2(\mathscr C\setminus B_R(0))},
\|D^2\phi(u,u)\|_{L^2(\mathscr C\setminus B_R(0))}
\le
e^{-\frac{R^2}{8p}}.
\]
Therefore,
\begin{align*}
\|P_{\ker L}D^2\phi(\chi_Ru,\chi_Ru)\|_{L^2}
&\geq
\|P_{\ker L}D^2\phi(u,u)\|_{L^2}
\\
&\hspace{0.5cm}
-
\|P_{\ker L}
(D^2\phi(\chi_Ru,\chi_Ru)-D^2\phi(u,u))
\|_{L^2}
\\
&\geq
\kappa
-
\|D^2\phi(\chi_Ru,\chi_Ru)-D^2\phi(u,u)\|_{L^2}
\\
&\geq
\kappa
-
\|D^2\phi(\chi_Ru,\chi_Ru)-D^2\phi(u,u)\|_{L^2(\mathscr C\setminus B_R(0))}
\\
&\geq
\kappa
-
\|D^2\phi(\chi_Ru,\chi_Ru)\|_{L^2(\mathscr C\setminus B_R(0))}
\\
&\hspace{0.5cm}
-
\|D^2\phi(u,u)\|_{L^2(\mathscr C\setminus B_R(0))}
\\
&\geq
\kappa-2e^{-\frac{R^2}{8p}}.
\end{align*}
Hence,
\[
\|P_{\ker L}D^2\phi(\chi_Ru,\chi_Ru)\|_{L^2}\ge\frac{\kappa}{2}.
\]
Replacing $\kappa/2$ by $\kappa$, and using the bilinearity of $D^2\phi$, it follows that for arbitrary $u\in K_2$,
\[
\|P_{\ker L}D^2\phi(\chi_Ru,\chi_Ru)\|_{L^2}
\ge
\kappa\|u\|_{L^2}^2.
\]
This completes the proof.
\end{proof}

The next lemma is the main result of this section and is an application of Proposition ~\ref {prop 4.7}. It establishes a lower bound for the functional $\phi$ along the constructed flow.

\begin{lemma} {\label{lemma 4.8}}
    There exists a constant $\kappa_0 >0,$ independent of $\alpha,$ such that for all sufficiently small $|\alpha|$ and all $\tau \in [0, T_\alpha]$, we have
    $$\|\phi \big(\mathscr{T}^{\alpha,\mathscr{C}}_\tau \big)\|_{L^2} \geq \kappa_0 |\alpha|^2.$$
\end{lemma}

\begin{proof}
By triangle inequality, for every normal graph $v$ with sufficiently small $C^{2,\alpha}$-norm, we have
\begin{align*}
    \frac{1}{2} \|P_{\ker L} (D^2\phi(v ,v))\|_{L^2} \leq& \;\|P_{\ker L} \phi(v)\|_{L^2} + \|P_{\ker L} Lv\|_{L^2} +  \\
    &\; \|P_{\ker L} \big( \phi(v) - Lv - \frac{1}{2}D^2\phi(v,v) \big)\|_{L^2}. 
\end{align*}
Since $Lv \in (\ker L)^{\perp},$ we have 
$$P_{\ker L} L(v) = 0.$$ 
Using
    $$|\phi(v) - Lv - \frac{1}{2}D^2\phi(v,v)| \leq C (|v|+|\nabla v| + |\nabla^2 v|)^3,$$
we obtain
$$|\phi  (\chi_{R_\alpha - 1} u_{\alpha(\tau)}) - L(\chi_{R_\alpha - 1} u_{\alpha(\tau)})  - \frac{1}{2}D^2\phi(\chi_{R_\alpha - 1} u_{\alpha(\tau)},\chi_{R_\alpha - 1} u_{\alpha(\tau)}) |_{L^2} \leq C |\alpha|^3.$$
Consequently,
\begin{align*}
     &\|P_{\ker L} \big(\phi  (\chi_{R_\alpha -1} u_{\alpha(\tau)}) - L(\chi_{R_\alpha -1} u_{\alpha(\tau)})  - \frac{1}{2}D^2\phi(\chi_{R_\alpha -1}  u_{\alpha(\tau)},\chi_{R_\alpha -1} u_{\alpha(\tau)}) \big) \|_{L^2} \\
     &\leq \| \phi  (\chi_{R_\alpha -1} u_{\alpha(\tau)}) - L(\chi_{R_\alpha -1} u_{\alpha(\tau)})  - \frac{1}{2}D^2\phi(\chi_{R_\alpha -1}  u_{\alpha(\tau)},\chi_{R_\alpha -1} u_{\alpha(\tau)}) \|_{L^2} \\
     &\leq C|\alpha|^3.
\end{align*}
Let $c_2 = \frac{\kappa}{8N},$ where $\kappa$ is the constant from Proposition ~\ref{prop 4.7}. Then, Proposition ~\ref{prop 4.7} implies that, for $|\alpha|$ sufficiently small,
$$\frac{1}{2}  \|P_{\ker L} (D^2\phi(\chi_{R_\alpha - 1} u_{\alpha(\tau)}, \chi_{R_\alpha - 1} u_{\alpha(\tau)}))\|_{L^2} \geq c_2 |\alpha|^2.$$ 
Here we used $|\alpha(\tau)| \geq \frac{1}{2}|\alpha|$ for all $\tau \in [0, T_\alpha].$ Combining the previous two estimates yields
$$\|P_{\ker L}\phi (\chi_{R_\alpha -1} u_{\alpha(\tau)}) \|_{L^2} \geq c_2 |\alpha|^2 - C |\alpha|^3 .$$
Therefore, after shrinking $|\alpha|$ if necessary,
$$\|P_{\ker L}\phi (\chi_{R_\alpha -1} u_{\alpha(\tau)}) \|_{L^2} \geq \frac{c_2}{2} |\alpha|^2.$$
Next, write 
$$V(\tau) = N(u_{\alpha(\tau)}) + v_\alpha(\tau).$$
For every normal graph $w$, let $L_w$ denote the linearization of $\phi$ at $w$, and let $Q_w$ denote the corresponding quadratic remainder. Then,
\[
\phi(w+v) = \phi(w) + L(v) + (L_w-L)(v) + Q_w(v).
\]
Taking
\[
w=\chi_{R_\alpha-1}u_{\alpha(\tau)},
\qquad
v=\chi_{R_\alpha-1}V(\tau),
\]
we obtain
\[
\begin{aligned}
\phi\bigl(\chi_{R_\alpha-1}(u_{\alpha(\tau)}+V(\tau))\bigr)
&=
\phi(\chi_{R_\alpha-1}u_{\alpha(\tau)})
+
L(\chi_{R_\alpha-1}V(\tau))
\\
&\quad
+
\bigl(L_{\chi_{R_\alpha-1}u_{\alpha(\tau)}}-L\bigr)
(\chi_{R_\alpha-1}V(\tau))
\\
&\quad
+
Q_{\chi_{R_\alpha-1}u_{\alpha(\tau)}}
(\chi_{R_\alpha-1}V(\tau)).
\end{aligned}
\]
Since
\[
P_{\ker L}L(\chi_{R_\alpha-1}V(\tau))=0,
\]
we conclude that
\begin{align*}
\left\|\phi\bigl(\chi_{R_\alpha-1}(u_{\alpha(\tau)}+V(\tau))\bigr)\right\|_{L^2}
&\ge \left\|P_{\ker L} \phi\bigl(\chi_{R_\alpha-1}(u_{\alpha(\tau)}+V(\tau))\bigr) \right\|_{L^2} \\
&\ge
\frac{c_2}{2}|\alpha|^2
-
\left\|
\bigl(L_{\chi_{R_\alpha-1}u_{\alpha(\tau)}}-L\bigr)
(\chi_{R_\alpha-1}V(\tau))
\right\|_{L^2} \\
&\hspace{0.5 cm}
-
\left\|
Q_{\chi_{R_\alpha-1}u_{\alpha(\tau)}}
(\chi_{R_\alpha-1}V(\tau))
\right\|_{L^2},
\end{align*}
Now, if $|\alpha|$ is small enough
\begin{align*}
    \left\|\bigl(L_{\chi_{R_\alpha-1}u_{\alpha(\tau)}}-L\bigr) (\chi_{R_\alpha-1}V(\tau)) \right\|_{L^2} &\leq C \|u_{\alpha(\tau)} \|_{C^2 (\mathscr{C} \cap B_{R_{\alpha}(0)})}\|V(\tau)\|_{H^2} \\
    &\leq C |\alpha|R_\alpha^2 (|\alpha|^2R_{\alpha}^6 + |\alpha|^{\frac{5}{2}}) \\
    &\leq C|\alpha|^{\frac{5}{2}}.
\end{align*}
Similarly, 
\begin{align*}
    \|Q_{\chi_{R_\alpha-1}u_{\alpha(\tau)}} (\chi_{R_\alpha-1}V(\tau))\|_{L^2} &\leq C(1 + \|u_{\alpha(\tau)} \|_{C^2 (\mathscr{C} \cap B_{R_{\alpha}(0)})})\|V(\tau)\|^2_{H^2} \\
    &\leq C|\alpha|^3.
\end{align*}
Since both error terms are of higher order than $|\alpha|^2$, it follows that 
$$\|\phi \big(\chi_{R_\alpha-1} (u_{\alpha(\tau)} + V(\tau) \big) \|_{L^2} \geq  \frac{c_2}{4}|\alpha|^2$$
for all sufficiently small $|\alpha|$. Arguing exactly as in the proof of Proposition ~\ref{prop 4.7}, for small $|\alpha|$, we have
$$\|\phi \big(\chi_{R_\alpha -1} (u_{\alpha(\tau)} + V(\tau) \big) \|_{L^2(\mathscr{C} \setminus B_{R_\alpha-2}(0))} \leq  \frac{c_2}{8}|\alpha|^2 $$
Therefore,
$$\|\phi \big(\chi_{R_\alpha-1} (u_{\alpha(\tau)} + V(\tau) \big) \|_{L^2(\mathscr{C} \cap B_{R_\alpha-2}(0))} \geq \frac{c_2}{8}|\alpha|^2.$$
Finally, since
\[
\phi \big(\mathscr{T}^{\alpha,\mathscr{C}}_\tau \big) = \phi \big(\chi_{R_\alpha-1} (u_{\alpha(\tau)} + V(\tau) \big) 
\]
on $\mathscr C\cap B_{R_\alpha -2}(0)$, we have
\begin{align*}
    \|\phi \big(\mathscr{T}^{\alpha,\mathscr{C}}_\tau \big) \|_{L^2} &\geq \|\phi \big(\mathscr{T}^{\alpha,\mathscr{C}}_\tau \big) \|_{L^2(\mathscr{C} \cap B_{R_\alpha-2}(0))} \\
    &\geq \frac{c_2}{8}|\alpha|^2.
\end{align*}
Renaming $\frac{c_2}{8}$ as $\kappa_0$ completes the proof.

\end{proof}

The same argument applies to $\mathcal{W}^{\alpha,\mathscr{C}}$. Hence $\mathcal{W}^{\alpha,\mathscr{C}}$ is not a shrinker whenever $|\alpha| \neq 0$.

\section{Non-concentration}\label{section 5}

The goal of this section is to prove the key non-concentration result, Proposition ~\ref{prop 5.2} below similar to Li-Sz\'ekelyhidi [\citealp{li2024singularity}, Lemma $34$], and as a consequence to prove a three annulus lemma, Lemma ~\ref{lemma 5.6}. We will need to consider rescaled flow $M_\tau$ very close to one of the comparison flows $\mathscr {T}^{\alpha, \mathscr{C}}$ constructed in the last section.

We say that a hypersurface $M$ is a $c$-graph over a hypersurface $N$ on the ball $B_R(0)$ if we can write $M \cap B_R(0)$ as the normal graph over $N \cap B_R(0)$ of a function $u$ satisfying $|u| + |\nabla u| + |\nabla^2 u| \leq c.$ 

Now we define our Gaussian weighted $L^2$ distance function following \cite{lotay2022neck} (see also \cite{li2024singularity}). Fix a sufficiently small constant $\eta >0$, to be specified in Proposition ~\ref{prop 5.2}. Let $R \leq R_\alpha-1$ denote the largest radius for which $M$ is a $\eta$-graph over $\mathscr {T}^{\alpha, \mathscr{C}}$ on $B_R(0)$ at time $\tau$. We then define $D_{\mathscr{T}^{\alpha, \mathscr{C}}_\tau} (M)$ by 
$$D_{\mathscr{T}^{\alpha, \mathscr{C}}_\tau} (M) = \bigg(\int_{M \cap B_{R}(0)} d_{\mathscr{T}^{\alpha, \mathscr{C}}_\tau}^2 \;e^{-\frac{|x|^2}{4}} d \mathcal{H}^n\bigg)^{\frac{1}{2}} + e^{-\frac{R^2}{8p_0}}.$$
Here $p_0 > 1$ is chosen sufficiently close to $1$ which will be chosen in Proposition ~\ref{prop 5.2}, and $d_{\mathscr {T}^{\alpha, \mathscr{C}}_\tau}$ is the distance function to $\mathscr {T}^{\alpha, \mathscr{C}}$ at time $\tau$.

Similarly, we define the Gaussian weighted $L^2$ distance from $\mathcal{W}^{\alpha, \mathscr{C}}$ to $M$ as follows. Let $R \leq R_\alpha-1$ denote the largest radius for which $M$ is an $\frac{\eta}{2}$-graph over $\mathcal{W}^{\alpha, \mathscr{C}}$ on $B_R(0)$. Then we define
$$D_{\mathcal{W}^{\alpha, \mathscr{C}}}(M) = \bigg(\int_{M \cap B_{R}(0)} d^2_{\mathcal{W}^{\alpha, \mathscr{C}}}\;e^{-\frac{|x|^2}{4}} d \mathcal{H}^n\bigg)^{\frac{1}{2}} +  e^{-\frac{R^2}{8p_0}}.$$
For $\epsilon >0,$ define 
$$R_\epsilon:= (8p_0 |\ln \epsilon|)^{\frac{1}{2}}.$$
We begin with the following elementary lemma. This lemma is similar to [\citealp{colding2015uniqueness}, Lemma~5.39]. Related statements can also be found in [\citealp{chodosh2021uniqueness}, Proposition~5.31 and \citealp{li2024singularity}, Lemma~30]. It shows that the graphicality of the rescaled flow can be extended to a larger region for a short time, provided that the flow is sufficiently close to the model flow at the initial time.

\begin{lemma}\label{lemma 5.1}
Let $\eta_0>0$. Then there exist constants $\eta_1>0$ and $R_1>0,$
depending only on $\eta_0$ with the following property. Suppose that
$\tau+1\le T_\alpha$ and $M_{\tau_0}$ is an $\eta_1$-graph over
$\mathscr{T}^{\alpha,\mathscr C}_{\tau}$ on $B_R(0)$ for some
$R>0$. Then for every $s\in[0,1]$, $M_{\tau_0+s}$ is an
$\eta_0$-graph over $\mathscr{T}^{\alpha,\mathscr C}_{\tau+s}$ on
$B_{\min\{e^{\frac{s}{2}}(R-R_1), R_\alpha-1\}}(0).$
\end{lemma}

\begin{proof}
Define the shifted rescaled flows by
\[
\widetilde{M}_s := M_{\tau_0+s}, \qquad 
\widetilde{T}_s := \mathscr{T}^{\alpha,\mathscr C}_{\tau+s}.
\]
We compare $\widetilde{M}_s$ and $\widetilde{T}_s$ for $s \in [0,1]$. Let $\widehat{M}_t$ and $\widehat{T}_t$ denote the corresponding unrescaled flows under the standard relation $t = -e^{-\tau}$.

We normalize the time so that $s = 0$ corresponds to $t = -1$ in the unrescaled variables. With this normalization, the parameter $s \in [0,1]$ corresponds to $t = -e^{-s}$, and both unrescaled flows are defined on the common time interval $[-1, -e^{-1}]$ for $|\alpha|$ sufficiently small.

Since for every $\tau \in [0,T_\alpha], \mathscr T^{\alpha,\mathscr C}_\tau$ is a sufficiently small $C^{2,\alpha}$-graph over the cylinder on $B_{R_\alpha -1}(0)$, the corresponding unrescaled flow $\widehat{\mathscr T}^{\alpha,\mathscr C}$ has uniformly controlled geometry on $[-1,-e^{-1}]$. Hence, by the pseudolocality theorem of Ilmanen-Neves-Schulze \cite[Theorem~1.5]{ilmanen2019short} (see also White's local Brakke regularity theorem \cite{white2005local}), there exist constants $\eta_1>0$ and $R_1>0$, depending only on $\eta_0$, such that if $\widehat M_{-1}$ is an $\eta_1$-graph over $\widehat{\mathscr T}^{\alpha,\mathscr C}_{-1}$ on $B_R(0)$, then for every $t\in[-1,-e^{-1}], \widehat M_t$ is an $\eta_0$-graph over $\widehat{\mathscr T}^{\alpha,\mathscr C}_t$ on $B_{R-R_1}(0)$.

Rescaling back to the normalized flow at time $s$, where $t=-e^{-s}$, the spatial scales are multiplied by $e^{\frac{s}{2}}$, and the result follows.
\end{proof}

Now we will prove the following non-concentration estimate for the distance function. 

\begin{prop} {\label{prop 5.2}}
     There are constants $C, \epsilon_0, \theta >0$ and $p >1$ such that we have the following. Suppose that $\tau + 1 \leq T_\alpha$, $|\alpha| \leq \epsilon_0$ and $D_{\mathscr {T}^{\alpha, \mathscr{C}}_\tau}(M_{\tau_0}) \leq \epsilon_0.$ Let $\epsilon = D_{\mathscr{T}^{\alpha,\mathscr C}_\tau}(M_{\tau_0})$. Then $M_{\tau_0+1}$ is an $\frac{\eta}{4}$-graph over $\mathscr {T}^{\alpha, \mathscr{C}}_\tau$ on $B_{\min\{(1+\theta)R_\epsilon, R_\alpha-1\}}(0)$, and
     $$\bigg(\int_{M_{\tau_0 + 1} \cap B_{(1+\theta)R_\epsilon}(0)} d_{\mathscr {T}^{\alpha, \mathscr{C}}_{\tau +1}}^{2p} \;e^{-\frac{|x|^2}{4}} d \mathcal{H}^n\bigg)^{\frac{1}{2p}} \leq C \epsilon.$$
\end{prop}

\begin{proof}
    By definition of $D_{\mathscr {T}^{\alpha, \mathscr{C}}_\tau}(M_{\tau_0})$, there exists a $R>0$ such that $M_{\tau_0}$ is $\eta$-graph over $\mathscr {T}^{\alpha, \mathscr{C}}_{\tau}$ on $B_R(0)$ and 
    $$\bigg(\int_{M_{\tau_0} \cap B_{R}(0)} d_{\mathscr{T}^{\alpha, \mathscr{C}}_\tau}^2 \;e^{-\frac{|x|^2}{4}} d \mathcal{H}^n\bigg)^{\frac{1}{2}} + e^{-\frac{R^2}{8p_0}} = \epsilon.$$
    Since $e^{-\frac{R^2}{8p_0}} \leq \epsilon,$ we have $R \geq R_\epsilon.$

    Now fix a sufficiently small constant $\eta' >0$. Fix $\eta = \eta_1(\eta')$ to be the constant given by Lemma ~\ref{lemma 5.1}. Since $M_{\tau_0}$ is an $\eta$-graph over $\mathscr {T}^{\alpha, \mathscr{C}}_\tau$ on $B_{R_\epsilon}(0)$, Lemma ~\ref{lemma 5.1} implies that for every $s \in [0,1]$ $M_{\tau_0 + s}$ can be written as an $\eta'$ graph of a function $v(x,s)$ over $\mathscr{T}^{\alpha, \mathscr{C}}_{\tau +s}$ on $B_{e^{\frac{s}{2}} (R_\epsilon - R_1)}(0)$. In particular, $|v_s| + |\nabla v_s| + |\nabla^2 v_s| \leq \eta'.$ Provided $\eta'$ is sufficiently small, the graph function satisfies the differential inequality
    $$|\partial_s v - \Delta v + \frac{1}{2} x. \nabla v| \leq C \big(|v| + |\nabla v| \big).$$
    for a uniform constant $C$. Arguing exactly as in the proof of [\citealp{li2024singularity}, Lemma $34$], one absorbs the gradient term and computes the evolution of $|v|^{\frac{3}{2}}$ to obtain
    $$\partial_s |v_s|^{\frac{3}{2}} - \Delta |v_s|^{\frac{3}{2}} + \frac{1}{2} x. \nabla |v_s|^{\frac{3}{2}} \leq C |v_s|^{\frac{3}{2}},$$
    possibly after increasing the constant $C$. 
    
    Writing $M_{\tau_0+s}$ as the normal graph $F_\tau(x) = x + v_s(x)n_{\mathscr {T}^{\alpha, \mathscr{C}}_{\tau+s}}(x),$ over $\mathscr {T}^{\alpha, \mathscr{C}}_{\tau+s}$, we have $d_{\mathscr {T}^{\alpha, \mathscr{C}}_{\tau+s}}(F_s(x)) = |v_s|(x)$. Moreover, on the graphical region the parametrization $F_s$ has uniformly bounded Jacobian, so the induced volume forms on $M_{\tau_0+s}$ and $\mathscr{T}_{\tau+s}^{\alpha, \mathscr{C}}$ are uniformly equivalent on the graphical region. 

    Since the cylinder satisfies 
    $$\langle p, n_{\mathscr{C}}(p) \rangle = \sqrt{2k},$$  
    and $\mathscr {T}^{\alpha, \mathscr{C}}_{\tau+s}$ is a $C^{2,\alpha}$-graph over $\mathscr{C}$ on $B_{R_\alpha -1}(0)$ with $C^{2,\alpha}$ norm bounded by $C|\alpha|^{\frac{1}{4}}$, we have
    $$|n_{\mathscr {T}^{\alpha, \mathscr{C}}_{\tau+s}} - n_{\mathscr{C}}| \leq C|\alpha|^{\frac{1}{4}}.$$
    Hence, 
    $$|\langle x, n_{\mathscr {T}^{\alpha, \mathscr{C}}_{\tau+s}}(x) \rangle| \leq \sqrt{2k} + C|\alpha|^{\frac{1}{4}}|x|.$$
    Since $|x| \leq R_\alpha-1$ on the graphical region and $R_\alpha |\alpha|^{\frac{1}{4}}$ is uniformly bounded, it follows that $\langle x, n_{\mathscr {T}^{\alpha, \mathscr{C}}_{\tau+s}}(x)\rangle$ is uniformly bounded. Therefore, 
    $$|F_s(x)|^2 = |x|^2 + 2 v_s(x) \langle x, n_{\mathscr {T}^{\alpha, \mathscr{C}}_{\tau+s}}(x)\rangle + |v_s(x)|^2,$$
    and the Gaussian weights $e^{-\frac{|x|^2}{4}}$ and $e^{-\frac{|F_s(x)|^2}{4}}$ are uniformly equivalent.
    Consequently, the integrals of $d^{2p}_{\mathscr {T}^{\alpha, \mathscr{C}}_{\tau+s}} e^{-|x|^2/4}$ over a region in $M_{\tau_0+s}$ and of $|v_s|^{2p} e^{-|x|^2 /4}$ over the corresponding region on $\mathscr {T}^{\alpha, \mathscr{C}}_{\tau+s}$ are uniformly equivalent.
    
    Define
    $$f(x,s) = e^{-Cs} |v_s|^{\frac{3}{2}}.$$
    Then $f$ is a subsolution of the drift heat equation on $B_R(0)$ on the time interval $[0,1].$ Also, along any rescaled flow 
    $$(\partial_s - \Delta) |x|^2 = -2n,$$ 
    so $e^{-s} |x|^2 - (R_\epsilon - C_1 - 1)^2$ is a subsolution of the drift heat equation along the rescaled flow. Let us define the function
    $$\widehat{f} =
    \begin{cases}
         f, & \text{for}\; |x| \leq  e^{\frac{s}{2}} (R_\epsilon - C_1 - 1), \\
         \max \big\{f, e^{-s} |x|^2 - (R_\epsilon - C_1 - 1)^2\big\}, & \text{for}\; e^{\frac{s}{2}} (R_\epsilon - C_1 - 1) < |x| \leq e^{\frac{s}{2}} (R_\epsilon - C_1), \\
         e^{-s}|x|^2 - (R_\epsilon - C_1 - 1)^2, & \text{for}\; |x| > e^{\frac{s}{2}} (R_\epsilon - C_1).
    \end{cases}$$
    Note that on $B_{R_\epsilon}(0), |v_s| \leq \eta,$ hence we can choose $\epsilon$ small enough (which implies $R_\epsilon$ large) so that when $|x| \to e^{\frac{s}{2}}(R_\epsilon - C_1),$ one has $e^{-s} |x|^2 - (R_\epsilon - C_1 - 1)^2 \geq f = e^{-Cs} |v_s|^{\frac{3}{2}}.$ It follows that $\widehat{f}$ is continuous and is again a subsolution of the drift heat equation along $\mathscr {T}^{\alpha, \mathscr{C}}_{\tau+s}$ on the time interval $s \in [0, 1].$ Now,
    $$\widehat{f}^\frac{4}{3} \leq 
    \begin{cases}
         |v_s|^2, & \text{for}\; |x| \leq e^{\frac{s}{2}} (R_\epsilon - C_1 - 1), \\
         (|v_s|^2 + |x|^{\frac{8}{3}} ), & \text{for}\; e^{\frac{s}{2}} (R_\epsilon - C_1 - 1) < |x| \leq e^{\frac{s}{2}} (R_\epsilon - C_1),\\
         |x|^{\frac{8}{3}}, & \text{for}\; |x| > e^{\frac{s}{2}} (R_\epsilon - C_1).
    \end{cases}$$
    Choose $q > 1$ so that $1< q < p_0.$ Since $R_\epsilon \to \infty$ as $\epsilon \to 0,$ for all sufficiently small $\epsilon > 0,$ we have
    $$\int_{M_{\tau_0} \setminus B_{R_\epsilon - C_1 - 1}(0)} |x|^{\frac{8}{3}} e^{-\frac{|x|^2}{4}} d\mathcal{H}^n \leq e^{-\frac{(R_\epsilon - C_1 - 1)^2}{4q}},$$
    and
    $$e^{-\frac{(R_\epsilon - C_1 - 1)^2}{4q}} \leq e^{-\frac{R_\epsilon^2}{4p_0}}.$$
    As $R \geq R_\epsilon,$
    $$e^{-\frac{R_\epsilon^2}{4p_0}} \leq e^{-\frac{R^2}{4p_0}} \leq \epsilon^2.$$
    This implies
    $$\int_{\mathscr {T}^{\alpha, \mathscr{C}}_{\tau_0}} \widehat{f}^{\frac{4}{3}} e^{-\frac{|x|^2}{4}} d\mathcal{H}^n \leq C\epsilon^2.$$
    Applying Ecker’s log-Sobolev inequality \cite{ecker2000logarithmic}, we obtain a $p > 1$ and a $C >0$, such that
    $$\int_{\mathscr {T}^{\alpha, \mathscr{C}}_{\tau_0+1}} \widehat f^{\frac{4p}{3}} \;e^{-\frac{|x|^2}{4}} d \mathcal{H}^n \leq C \epsilon^{2p}.$$
    We take $p(0) = \frac{4}{3}$ in \cite{ecker2000logarithmic}, then $p(1) = 1 + \frac{e^2}{3} > p(0).$ So here $p = \frac{p(1)}{p(0)} > 1.$ Choose $\theta > 0$ small enough so that 
    $$1+ \theta < \min \{e^{\frac{1}{2}} , \sqrt p \}.$$
    As $e^{\frac{1}{2}} > 1+\theta,$ so for all small $\epsilon > 0$,  
    $$e^{\frac{1}{2}}(R_\epsilon - C_1) \geq (1+ \theta) R_\epsilon.$$
    On $B_{e^{\frac{1}{2}} (R_\epsilon - C_1)}(0)$, we have
    $$f(.,1) = e^{-C} |v(x,1)|^{\frac{3}{2}} \leq \widehat f(.,1).$$
    Therefore,
    $$\int_{M_{\tau_0 + 1} \cap B_{(1+ \theta) R_\epsilon}(0)} d_{\mathscr {T}^{\alpha, \mathscr{C}}_{\tau +1}}^{2p} \;e^{-\frac{|x|^2}{4}} d \mathcal{H}^n \leq C \epsilon^{2p}.$$
    Using the monotonicity formula centered at different points as in the proof of [\citealp{lotay2022neck}, Lemma $3.5 (2)$], this integral estimate implies the pointwise bound
    $$d_{\mathscr {T}^{\alpha, \mathscr{C}}_{\tau + 1}}^2 \leq C e^\frac{|x|^2}{4p} \epsilon^2$$
    on $M_{\tau_0+1} \cap B_{(1+ \theta) R_\epsilon}(0).$
    This together with interior estimates in turn implies the pointwise bounds for 
    $|\nabla v|$ and $ |\nabla^2 v|$ for $|x| \leq (1+ \theta) R_\epsilon.$
    Therefore, $|v(x,1)| + |\nabla v(x,1)| + |\nabla^2 v(x,1)| \leq C e^\frac{|x|^2}{8p} \epsilon.$ Note that for $|x| \leq (1+ \theta) R_\epsilon,$ 
    $$e^\frac{|x|^2}{8p} \epsilon \leq e^\frac{(1 +\theta)^2 R_\epsilon^2}{8p}. e^{-\frac{R_\epsilon^2}{8p_0}} = e^ {\frac{R_\epsilon^2}{8} \Big( \frac{(1+\theta)^2}{p} - \frac{1}{p_0} \Big)}.$$
    Because $1+ \theta < \sqrt p,$ we have $(1+ \theta)^2 < p$. We now fix $p_0 \in (1,\frac{9}{8})$ sufficiently close to $1$ so that 
    $$\frac{(1+\theta)^2}{p} < \frac{1}{p_0}.$$
    Then the exponent is negative. Hence, for all sufficiently small $\epsilon, |v(x,1)| + |\nabla v(x,1)| + |\nabla^2 v(x,1)| \leq \frac{\eta}{4}$ on $B_{(1+ \theta) R_\epsilon}(0).$ This completes the proof.
\end{proof}

We will use the non-concentration estimate in the following form. This follows from Proposition ~\ref{prop 5.2}.

\begin{cor} \label{cor 5.3}
There are $\epsilon_0, C_0 > 0$ satisfying the following. For any $\delta>0$, there is a $R_0$ depending on $\delta$ such that if $\tau + 1 \leq T_\alpha$, $|\alpha| \leq \epsilon_0$ and $D_{\mathscr {T}^{\alpha, \mathscr{C}}_\tau}(M_{\tau_0}) \leq \epsilon_0$ then we obtain the improved estimate
\[
\begin{aligned}
D_{\mathscr {T}^{\alpha, \mathscr{C}}_{\tau+1}}(M_{\tau_0+1})
\leq
\max \bigg\{
&\delta D_{\mathscr {T}^{\alpha, \mathscr{C}}_\tau}(M_{\tau_0})
+
C_0
\bigg(
\int_{M_{\tau_0+1}\cap B_{R_0}(0)}
d_{\mathscr {T}^{\alpha, \mathscr{C}}_{\tau+1}}^2
e^{-\frac{|x|^2}{4}}
d\mathcal{H}^n
\bigg)^{\frac12},
\\
& C_0e^{-\frac{(R_\alpha-1)^2}{8p_0}}
\bigg\}.
\end{aligned}
\]
\end{cor}

\begin{proof}
Assume $ D_{\mathscr {T}^{\alpha, \mathscr{C}}_\tau}(M_{\tau_0}) =\epsilon \leq \epsilon_0,$ where $\epsilon_0$ is from Proposition ~\ref{prop 5.2}. By the proof of Proposition ~\ref{prop 5.2}, $M_{\tau_0+1}$ is $\eta$ graph over $\mathscr{T}_{\tau+1}^{\alpha,\mathscr{C}}$ on $B_{(1+\theta)R_\epsilon}(0)$ and 
$$\bigg(\int_{M_{\tau_0 + 1} \cap B_{(1+\theta)R_\epsilon}(0)} d_{\mathscr {T}^{\alpha, \mathscr{C}}_{\tau +1}}^{2p} \;e^{-\frac{|x|^2}{4}} d \mathcal{H}^n\bigg)^{\frac{1}{2p}} \leq C \epsilon.$$
By Hölder's inequality and the previous estimate,
$$\int_{M_{\tau_0 + 1} \cap B_{(1+\theta)R_\epsilon}(0) \setminus B_{R_0}(0)} d_{\mathscr {T}^{\alpha, \mathscr{C}}_{\tau+1}}^2 \;e^{-\frac{|x|^2}{4}} d\mathcal{H}^n \leq C \epsilon^2 \bigg(\int_{M_{\tau_0+1} \setminus B_{R_0}(0)} \;e^{-\frac{|x|^2}{4}} d\mathcal{H}^n \bigg)^{1- \frac1p}.$$
Since the flow has uniformly bounded area ratios, choosing $R_0$ sufficiently large depending on $\gamma$ gives 
$$\bigg(\int_{M_{\tau_0+1} \cap B_{(1+\theta)R_\epsilon}(0)\setminus B_{R_0}(0)} \;e^{-\frac{|x|^2}{4}} d\mathcal{H}^n \bigg)^{1-\frac1p} \leq \frac{1}{100}C^{-1}\gamma^2.$$
So, combining
$$\bigg(\int_{M_{\tau_0 + 1} \cap B_{(1+\theta)R_\epsilon}(0)\setminus B_{R_0}(0)} d_{\mathscr {T}^{\alpha, \mathscr{C}}_{\tau +1}}^2 \;e^{-\frac{|x|^2}{4}} d \mathcal{H}^n\bigg)^{\frac{1}{2}} \leq \frac{1}{10} \epsilon \gamma.$$
Consequently,
\begin{align*}
    \int_{M_{\tau_0 + 1} \cap B_{(1+\theta)R_\epsilon}(0)} d_{\mathscr {T}^{\alpha, \mathscr{C}}_{\tau+1}}^2 \;e^{-\frac{|x|^2}{4}} d\mathcal{H}^n  &\leq \frac{1}{100} \epsilon^2 \gamma^2 +  \int_{M_{\tau_0 + 1} \cap B_{R_0}(0)} d_{\mathscr {T}^{\alpha, \mathscr{C}}_{\tau+1}}^2 \;e^{-\frac{|x|^2}{4}} d\mathcal{H}^n .
\end{align*}
By definition,
$$D_{\mathscr{T}^{\alpha, \mathscr{C}}_{\tau+1}} (M_{\tau_0+1}) = \bigg(\int_{M_{\tau_0+1} \cap B_{R}(0)} d_{\mathscr{T}^{\alpha, \mathscr{C}}_{\tau+1}}^2 \;e^{-\frac{|x|^2}{4}} d\mathcal{H}^n \bigg)^{\frac{1}{2}} + e^{-\frac{R^2}{8p_0}},$$
where $R \leq R_\alpha-1$ is the largest radius such that $M_{\tau_0+1}$ is an $\eta$-graph over $\mathscr {T}^{\alpha, \mathscr{C}}$ on $B_R(0)$ at time $\tau+1.$

Proposition \ref{prop 5.2} implies that $R \geq (1+\theta)R_\epsilon.$ If $(1+\theta)R_\epsilon \geq R_\alpha-1$, then we are already done. So we can assume that $(1+\theta)R_\epsilon < R_\alpha-1.$ Therefore,
$$D_{\mathscr {T}^{\alpha, \mathscr{C}}_{\tau+1}}(M_{\tau_0+1})^2 \leq \frac{1}{100} \epsilon^2 \gamma^2 +  \int_{M_{\tau_0 + 1} \cap B_{R_0}(0)} d_{\mathscr {T}^{\alpha, \mathscr{C}}_{\tau+1}}^2 \;e^{-\frac{|x|^2}{4}} d\mathcal{H}^n  + 2e^{-\frac{(1+\theta)^2R_\epsilon^2}{8p_0}}.$$
Finally, by the definition of $R_\epsilon,$
$$2e^{-\frac{(1+\theta)^2R_\epsilon^2}{8p_0}} = 2\epsilon^{2(1+\theta)^2} < \frac{1}{100} \gamma^2\epsilon^2.$$
if $\epsilon$ is sufficiently small, depending on $\theta$. Combining the above estimates proves the claimed inequality.
\end{proof}

We now discuss a $3$-annulus lemma for solutions of the drift heat equation. Let $u$ be a solution to the linearized equation 
$$(\partial_s - L)u = 0,$$
and define
$$I^2(u,s) = \int_\mathscr{C} u^2 e^{-\frac{|x|^2}{4}} d\mathcal{H}^n.$$
It is known by Colding-Minicozzi \cite{colding2020parabolic} that $\log I(s)$ is convex, and is linear only when $u = e^{cs} u_0$ for an eigenfunction $u_0$ of $L$ with eigenvalue $c$. They have proved the log convexity argument for the operator of the form $L_\phi + \lambda$ (in our case we can consider $\lambda = |A_\mathscr{C}|^2 + \frac{1}{2}$ = 1), where 
$$L_\phi(u) = e^\phi\; \text{div} (e^{-\phi} \nabla u).$$   
We will require a $3$-annulus type lemma for solutions of our linearized equation. To obtain this, we will need a slightly stronger result than the log-convexity argument of \cite{colding2020parabolic}. Similar results are proven in [\citealp{lotay2022neck}, Lemma $6.1$] and Simon [\citealp{simon1984isolated}, Lemma $3.3$].

\begin{lemma}\label{lemma 5.4}
There exist constants $0<\delta_1<\delta_2<1$ such that the following holds. Let $u$ be a solution of the linearized equation. If
\[
I(u,s+1)\ge e^{\delta_1} I(u,s),
\]
then
\[
I(u,s+2)\ge e^{\delta_2} I(u,s+1).
\]
\end{lemma}

\begin{proof}
Write
\[
u(s)=\sum_j a_j e^{-\lambda_j s}\phi_j,
\]
where $\{\phi_j\}$ is an orthonormal eigenbasis of the linearized operator. Then
\[
I(u,s)^2=\sum_j a_j^2 e^{-2\lambda_j s}.
\]
By time translation, it suffices to consider $s = 0$. Let $\epsilon>0$ sufficiently small so that the only eigenvalue in $(-3\epsilon, 3\epsilon)$ is $0$. Then
\begin{align*}
I(u,0)^2 + e^{-4\epsilon} I(u,2)^2
&= \sum_j a_j^2 + \sum_j a_j^2 e^{-4\epsilon-4\lambda_j} \\
&= \sum_j a_j^2 e^{-2\epsilon-2\lambda_j}
\left(e^{2\epsilon+2\lambda_j}+e^{-2\epsilon-2\lambda_j}\right). 
\end{align*}
We have chosen $\epsilon$ so that if $\lambda_j \neq 0,$ then $|\lambda_j| \geq 3\epsilon.$ Consequently, there exists a $c_\epsilon > 0$ such that for all $j$,
\[
e^{2\epsilon+2\lambda_j}+e^{-2\epsilon-2\lambda_j}\ge 2(1+c_\epsilon),
\]
and hence
\begin{align*}
    I(u,0)^2 + e^{-4\epsilon} I(u,2)^2 &\ge 2(1+c_\epsilon) \sum_j a_j^2 e^{-2\epsilon-2\lambda_j} \\
    &= 2(1+c_\epsilon) e^{-2\epsilon} I(u,1)^2.
\end{align*}
Let $0 < \delta_1 < \epsilon.$ Assume $I(u,1) \geq e^{\delta_1} I(u,0).$ Then 
\begin{align*}
    I(u,2)^2 &\geq e^{4\epsilon} \big(2(1+c_\epsilon) e^{-2\epsilon} I(u,1)^2 - I(u,0)^2 \big) \\
    &\geq e^{4\epsilon} \big(2(1+c_\epsilon) e^{-2\epsilon} I(u,1)^2 - e^{-2\delta_1} I(u,1)^2 \big) \\
    &\geq e^{2\epsilon} \big(2(1+c_\epsilon) - e^{2\epsilon-2\delta_1} \big) I(u,1)^2. 
\end{align*}
Define $\delta_2$ by 
$$e^{2\delta_2} =e^{2\epsilon} \big(2(1+c_\epsilon) - e^{2\epsilon-2\delta_1}\big).$$ 
Since $2(1+c_\epsilon) > 2 $, we can choose $\delta_1 < \epsilon$ sufficiently close to $\epsilon$, so that $2(1+c_\epsilon) - e^{2\epsilon-2\delta_1} >1$. This will imply $\delta_2> \epsilon.$ But  $\epsilon > \delta_1$, hence $\delta_2 > \delta_1,$ so the lemma follows.
\end{proof}

\begin{lemma} {\label{lemma 5.5}}
Let $u$ be a solution to the linearized equation. Suppose that when we decompose $u$ into eigenfunctions of $L$ at some time $s$, it has no component in the zero eigenspace. Then at least one of the following holds:
$$I(u,s+1) \geq e^{\delta_2}  I(u,s),$$ 
or
$$I(u,s-1) \geq e^{\delta_2} I(u,s).$$
where $\delta_2 >0$ is from Lemma ~\ref{lemma 5.4}.
\end{lemma} 
\begin{proof}
By time translation, it suffices to assume that $s = 1.$ Note that in Lemma~\ref{lemma 5.4}, we can choose $\epsilon$ arbitrarily small and take $c_\epsilon = \cosh(2\epsilon)-1.$ So we can take $\delta_2$ sufficiently small and ensure 
$$e^{2\lambda_j} + e^{-2\lambda_j} \geq 2e^{2\delta_2}$$ for all $\lambda_j \neq 0.$ Now as in the proof of Lemma ~\ref{lemma 5.4},
\begin{align*}
    I(u,0)^2 + I(u,2)^2 &= \sum_j a_j^2 + \sum_j a_j^2 e^{-4\lambda_j} \\
    &= \sum_j a_j^2 e^{-2\lambda_j}e^{2\lambda_j} + \sum_j a_j^2 e^{-2\lambda_j}e^{-2\lambda_j} \\
    &= \sum_j a_j^2 e^{-2\lambda_j} \big(e^{2\lambda_j} + e^{-2\lambda_j}\big).  
\end{align*}
Therefore,
\begin{align*}
     I(u,0)^2 + I(u,2)^2 &\geq 2e^{2\delta_2} \sum_j a_j^2 e^{-2\lambda_j} \\
     &= 2e^{2\delta_2}I(u,1)^2.
\end{align*}
Hence at least one of the terms on the left-hand side is at least $e^{2\delta_2} I(u,1)^2.$ This completes the proof.
\end{proof}

Due to technical considerations, it is convenient to introduce the following modified version of our distance function
$$\mathcal {D}_{\mathscr {T}^{\alpha, \mathscr{C}}_\tau}(M_{\tau_0}) := \sup_{s \in [0,2]} D_{\mathscr {T}^{\alpha, \mathscr{C}}_{\tau-s}}(M_{\tau_0-s}).$$
Note that $\mathcal{D}_{\mathscr {T}^{\alpha, \mathscr{C}}_\tau}(M_{\tau_0})$ gives bound for $M_{\tau_0-s}$ for $s \in [0,2].$ 

Using Proposition ~\ref{prop 5.2} and Lemma ~\ref{lemma 5.4} together with the non concentration estimate [Corollary \ref{cor 5.3}\;], we can show a three-annulus lemma for the distance function by a contradiction argument. 

\begin{lemma} \label{lemma 5.6}
Let $\delta_1, \delta_2$ be the constants from Lemma ~\ref{lemma 5.4} and $C_0$ be the constant from Corollary \ref{cor 5.3}. Then there exist $\epsilon_0 > 0$ and a large $L >0$ such that the following holds. Suppose that 
\begin{enumerate}
    \item $\tau + 2L \le T_\alpha,$
    \item $\tilde{\delta_1},\tilde{\delta_2} \in (\delta_1, \delta_2)$ with $\tilde{\delta_1} < \tilde{\delta_2},$
    \item $|\alpha| \leq \epsilon_0,$
    \item $\mathcal{D}_{\mathscr {T}^{\alpha, \mathscr{C}}_\tau}(M_{\tau_0}) \leq \epsilon_0.$
\end{enumerate}
If 
$$\mathcal{D}_{\mathscr {T}^{\alpha, \mathscr{C}}_{\tau+L}}(M_{\tau_0 + L}) > C_0  e^{- \frac{(R_\alpha-1)^2}{8p_0}},$$
and
$$\mathcal{D}_{\mathscr {T}^{\alpha, \mathscr{C}}_{\tau+L}}(M_{\tau_0 + L}) \geq e^{\tilde{\delta_1}L} \mathcal{D}_{\mathscr {T}^{\alpha, \mathscr{C}}_\tau} (M_{\tau_0}),$$
then
$$\mathcal{D}_{\mathscr {T}^{\alpha, \mathscr{C}}_{\tau + 2L}}(M_{\tau_0 + 2L}) \geq e^{\tilde{\delta_2}L} \mathcal{D}_{\mathscr {T}^{\alpha, \mathscr{C}}_{\tau+L}}(M_{\tau_0 + L}).$$
\end{lemma}

\begin{proof}
The proof is by contradiction and follows the arguments in [\citealp{szekelyhidi2020uniqueness}, Proposition $5.12$] and [\citealp{lotay2022neck}, Proposition $6.2$]. Suppose, for contradiction, that there exists a sequence of flows $M^i_{\tau_0 + s}$ such that $\mathcal{D}_{\mathscr {T}^{\alpha^i, \mathscr{C}}_\tau}(M^i_{\tau_0}) \le \epsilon_i,$ where $ |\alpha^i|,\epsilon_i \to 0$ and for a given large integer $L$ the conclusion of the proposition fails.
Let us denote 
$$d_i := \mathcal{D}_{\mathscr {T}^{\alpha^i, \mathscr{C}}_{\tau+L}}(M^i_{\tau_0+L})$$ 
Then 
$$\mathcal{D}_{\mathscr {T}^{\alpha^i, \mathscr{C}}_{\tau}} (M^i_{\tau_0}) \leq e^{-\tilde{\delta_1} L} d_i,$$
and
$$\mathcal{D}_{\mathscr {T}^{\alpha^i, \mathscr{C}}_{\tau+2L}} (M^i_{\tau_0+2L}) \leq e^{\tilde{\delta_2} L} d_i.$$
In particular $d_i >0.$ Let us define 
$$\mathcal{D}_i := \sup\limits_{s \in [0,2L]} \mathcal{D}_{\mathscr {T}^{\alpha^i, \mathscr{C}}_{\tau + s}} (M^i_{\tau_0 + s}).$$
We claim that there exists a constant $C_{L}$, depending on $L$ such that $\mathcal{D}_i \leq C_{L} d_i$. By the assumptions, for every $s \in [0,2L], M^i_{\tau_0 + s}$ can be written as the graph of a function $u_i(s)$ over $\mathscr {T}^{\alpha^i, \mathscr{C}}_{\tau + s}$ on larger and larger subsets of $\mathscr {T}^{\alpha^i, \mathscr{C}}_{\tau + s}$ as $i \to \infty$. After passing to a subsequence, the rescaled functions $d_i^{-1} u_i$ converge locally smoothly to a solution $u$ of the drift heat equation on $\mathscr{C}.$ If no constant $C_L$ exists as claimed, then $d_i^{-1} \mathcal{D}_i \to \infty$, and therefore it follows that $u(s) = 0$ for $s \in [L-1,L]$ which implies that $u$ is identically $0$ on $[0,2L].$ 

First note that from the definition of $\mathcal{D}_{\mathscr {T}^{\alpha^i, \mathscr{C}}_{\tau+s}}(M^i_{\tau_0+s})$, by our assumption and the fact that $d_i^{-1} \mathcal{D}_i \to \infty,$ we can assume that $s_i \in [2,2L].$ We now use Corollary ~\ref{cor 5.3} to show that this contradicts the assumption that $\mathcal{D}_i = \mathcal{D}_{\mathscr {T}^{\alpha^i, \mathscr{C}}_{\tau +s_i}}(M^i_{\tau_0 +s_i})$ for some $s_i \in [2, 2L].$ We have
$$\mathcal{D}_{\mathscr {T}^{\alpha^i, \mathscr{C}}_{\tau +s_i-1}}(M^i_{\tau_0 +s_i-1}) \leq \mathcal{D}_i.$$ 
Since the limit of $\mathcal{D}_i^{-1} u_i$ vanishes identically, it follows that for any $\epsilon, R_0 > 0$, once $i$ is sufficiently large depending on $\epsilon$ and $R_0$, we have
$$\bigg(\int_{M^i_{\tau_0 + s_i} \cap B_{R_0}(0)} d_{\mathscr {T}^{\alpha^i, \mathscr{C}}_{\tau +s_i}}^2 \;e^{-\frac{|x|^2}{4}} d\mathcal{H}^n \bigg)^{\frac{1}{2}} \leq \epsilon \mathcal{D}_i.$$
So applying Corollary \ref{cor 5.3} with $M_{\tau_0 +s_i -1}$, either we have
$$\mathcal{D}_{\mathscr {T}^{\alpha^i, \mathscr{C}}_{\tau +s_i}}(M^i_{\tau_0 +s_i}) \leq C_0 e^{- \frac{(R_{\alpha^i}-1)^2}{8p_0}},$$
but by our assumptions, this cannot happen. Therefore, given any $\delta > 0$, there exists $R_0$, depending on $\delta$, such that for all sufficiently large $i$,
$$\mathcal{D}_i = \mathcal{D}_{\mathscr {T}^{\alpha^i, \mathscr{C}}_{\tau +s_i}}(M^i_{\tau_0 +s_i}) \leq \delta \mathcal{D}_i + C_0 \bigg(\int_{M^i_{\tau_0 + s_i} \cap B_{R_0}(0)} d_{\mathscr {T}^{\alpha^i, \mathscr{C}}_{\tau +s_i }}^2 \;e^{-\frac{|x|^2}{4}} d\mathcal{H}^n \bigg)^{\frac{1}{2}} .$$
Choose $\delta = \frac{1}{4}$. This determines $R_0.$ Next choose $\epsilon$ so that $C_0\epsilon <\frac{1}{4}$, then for larger $i$ we will have $\mathcal{D}_i <\frac{\mathcal{D}_i}{2},$ which is again a contradiction. Therefore, we have 
$\mathcal{D}_i \leq C_L d_i$ for all $i$.

Now we can consider the rescaling $d_i^{-1} u_i$, extract a subsequence which converges to a function say $u$, a solution of the drift heat equation on $\mathscr{C},$ where the convergence is smooth on compact subsets $\mathscr{C}.$ By our assumption,
$$I(u,0) \leq Ce^{-\tilde{\delta_1} L},$$
and
$$I(u,2L) \leq Ce^{\tilde{\delta_2} L},$$
for a constant $C$ independent of $L.$

Let $\epsilon >0$. We claim that if $L$ is sufficiently large depending on $\epsilon,$ then $I(u,L-1) \leq \epsilon.$ To see this, note that we have two possibilities.

Let $I(u,k+1) \leq e^{\delta_1} I(u,k)$ for $k = 0,1,..., L-2$ then we have
$$I(u,L-1) \leq e^{(L-1)\delta_1} I(u,0) \leq Ce^{(L-1)\delta_1 - \tilde{\delta_1} L}.$$
Since $\tilde{\delta_1} > \delta_1$ this implies that $I(u,L-1) \leq \epsilon$ if $L$ is chosen sufficiently large.

So now we can assume that $I(u,k) \geq e^{\delta_1} I(u,k-1)$ for some $k \leq L-1,$ then by Lemma ~\ref{lemma 5.4} we have $I(u,k+1) \geq e^{\delta_2} I(u,k)$ for $k = L-1, L,...,2L-1.$ But this implies
$$I(u, L-1) \leq e^{-(L+1)\delta_2} I(u,2L) \leq Ce^{-(L+1)\delta_2 + \tilde{\delta_2} L}.$$
Since $\delta_2 > \tilde{\delta_2}$ this implies that $I(u,L-1) \leq \epsilon$ if $L$ is chosen sufficiently large.

Since $I(u, L-1) \leq \epsilon$, the same argument as in Proposition ~\ref{prop 5.2} (see also \cite{lotay2022neck}, Lemma $3.5 (2)$ and Proposition $5.3 (2)$), yields the pointwise estimate
$$|u(x,L)|^2 \leq C_1\epsilon^2e^{\frac{|x|^2}{4p}},$$
for some $C_1 > 0$ and $p > 1.$

We now use the local smooth convergence of the $d_i^{-1} u_i$ to $u$. For any fixed $R_0,$ this implies that, as $i\to \infty$, the functions $d_i^{-1} d_{\mathscr{T}^{\alpha^i, \mathscr{C}}_{\tau+s}}$ on $M^i_{\tau_0+s}$ converges to $|u|$ on $\mathscr{C} \cap B_{R_0}(0)$. Hence, for a given $R_0$, the local smooth convergence implies that, for sufficiently large $i$, the function $d_i^{-2} d^2_{\mathscr{T}^{\alpha^i, \mathscr{C}}_{\tau+s}}(M^i_{\tau_0+s})$ differs from $|u(.,s)|^2$ by at most $\epsilon^2$. This implies
\begin{align*}
    \int_{M^i_{\tau_0+L} \cap B_{R_0}(0)} d^2_{\mathscr{T}^{\alpha^i, \mathscr{C}}_{\tau+L}} e^{-\frac{|x|^2}{4}} d\mathcal{H}^n &\leq d_i^2\epsilon^2 + d_i^2 \int_{\mathscr{C} \cap B_{R_0}(0)} |u|^2 e^{-\frac{|x|^2}{4}} d\mathcal{H}^n \\
    &\leq C_2 d_i^2\epsilon^2,
\end{align*}
for some $C_2>0.$ Applting the non-concentration lemma, we estimate $\mathcal{D}_{\mathscr {T}^{\alpha^i, \mathscr{C}}_{\tau + L}}(M^i_{\tau_0 + L}).$ By our assumption and Proposition ~\ref{prop 5.2} we have
$$\mathcal{D}_{\mathscr {T}^{\alpha^i, \mathscr{C}}_{\tau + L-1}}(M^i_{\tau_0 + L-1}) \leq C_L d_i$$
for an $L$ dependent constant $C_L.$ Now let $\delta >0.$ Assuming $d_i > 0$ is sufficiently small,  Corollary \ref{cor 5.3} yields , for every $\delta >0$, there exists a $R_0$  depending on $\delta$ such that 
$$d_i = \mathcal{D}_{\mathscr {T}^{\alpha^i, \mathscr{C}}_{\tau + L}}(M^i_{\tau_0 + L}) \leq \delta C_L d_i + C_0 \bigg( \int_{M^i_{\tau_0+L} \cap B_{R_0}(0)} d^2_{\mathscr{T}^{\alpha^i, \mathscr{C}}_{\tau+L}} e^{-\frac{|x|^2}{4}} d\mathcal{H}^n \bigg)^{\frac{1}{2}}.$$
First choose $\epsilon > 0$ so that $C_0C_2^{\frac12}\epsilon \leq \frac{1}{4}.$ The choice of $\epsilon$ determines an $L$, such that we have the estimate $I(u,L-1) \leq \epsilon.$ Choosing $L$ determines the constant $C_L$. Next choose $\delta > 0$ so that $\delta C_L \leq \frac{1}{4}.$ The choice of $\delta$ determines the $R_0$ and then for sufficiently large $i$ we have 
$$\int_{M^i_{\tau_0+L} \cap B_{R_0}(0)} d^2_{\mathscr{T}^{\alpha^i, \mathscr{C}}_{\tau+L}} \leq C_2d_i^2\epsilon^2.$$ 
For such a large $i$, it implies $d_i = \frac{d_i}{4} + C_0C_2^{\frac12}\epsilon d_i \leq \frac{d_i}{2}$, which is a contradiction.

\end{proof}

\section{The main argument}\label{section 6}

In this section, we give the proof of our main result. Recall that $\mathcal{W}^{\alpha, \mathscr{C}}$ denotes the graph of the function $u_\alpha + N_{R_{\alpha}}(u_\alpha)$ over the truncated cylinder $\mathscr{C} \cap B_{R_\alpha -1}(0)$. Let
$$\mathcal{W} = \Big\{\mathcal{W}^{\alpha, \mathscr{RC}}: |\alpha| \leq \epsilon_1, \;\mathscr{R} \in SO(n+1)\Big\},$$
where $\mathscr{RC}$ denotes the image of $\mathscr{C}$ under rotation $\mathscr{R}.$ Here $\epsilon_1 > 0$ is a fixed constant so that $\mathcal{W}^{\alpha, \mathscr{RC}}, \mathscr{T}^{\alpha,\mathscr{RC}}$ are defined for $|\alpha| \leq \epsilon_1$ and Lemma~\ref{lemma 4.8} holds. 

We define the distance between two hypersurfaces as follows,
$$D(M_1, M_2) = \inf_{W} \Big\{D_W(M_1) + D_W(M_2) :  W \in \mathcal{W} \Big\}.$$
Following \cite{colding2015rigidity}, we define a metric \(d_V\) on the space of Radon measures on \(\mathbb{R}^{n+1}\) with finite Gaussian mass. Let
\(\{f_j\}_{j=1}^{\infty}\) be a countable dense subset of the unit ball of
\(C^0_c(\mathbb{R}^{n+1})\), equipped with the $C^0$ norm. For Radon measures \(\mu_1\) and \(\mu_2\), define
\begin{equation*}
d_V(\mu_1,\mu_2)
=
\sum_{j=1}^{\infty}2^{-j}
\left|
\int_{\mathbb{R}^{n+1}} f_j e^{-\frac{|x|^2}{4}}\,d\mu_1
-
\int_{\mathbb{R}^{n+1}} f_j e^{-\frac{|x|^2}{4}}\,d\mu_2
\right|.
\end{equation*}
It is straightforward to verify that \(d_V\) defines a metric on the space of Radon measures satisfying \(F(\mu)<\infty\). Moreover, the topology induced by \(d_V\) coincides with the standard weak topology on this space. In particular,
\[
\mu_i \rightharpoonup \mu
\quad\Longleftrightarrow\quad
d_V(\mu_i,\mu)\longrightarrow 0.
\]
By Proposition~\ref{prop 7.1}, we may choose the countable dense subset above so that each $f_j$ belongs to $C_c^1(\mathbb{R}^{n+1})$ and satisfies $(\frac{3}{2})^{-j}\|f_j\|_{C^1} \leq 1$. 

Let $p_0 >1$ be the constant introduced in Section \ref{section 5}. Throughout this
section, we fix an constant $1 < q < p_0$. We will show that the distance $D$ can be used to control the distance $d_v$. First, we need the following lemma.

\begin{lemma}\label{lemma 6.1}
Let $N\subset \mathbb{R}^{n+1}$ be a smooth hypersurface satisfying $\lambda(N)\leq\lambda_0$ and $|\langle x,n_N(x)\rangle|\leq C_0$.
Let $M$ be a normal graph over $N$ on $B_R(0)$ given by a function
$u$ with sufficiently small $C^2$ norm. Then, for every
$f\in C^1(\mathbb{R}^{n+1})$, there exists a constant
$C=C(C_0,\lambda_0)>0$ such that
\begin{multline*}
\left|
\int_{M\cap B_R(0)}
f e^{-\frac{|x|^2}{4}}\,d\mathcal H^n
-
\int_{N\cap B_R(0)}
f e^{-\frac{|x|^2}{4}}\,d\mathcal H^n
\right|
\\
\leq
C\|\nabla f\|_{C^0}\|u\|_{L^2}
+
C\|f\|_{C^0}
\left(
\|u\|_{L^2}\|\phi_N\|_{L^2}
+\|u\|_{L^2}\|u\|_{H^2}
+e^{-\frac{(R-1)^2}{4q}}
\right),
\end{multline*}
where all $L^2$ and $H^2$ norms are taken over $N\cap B_R(0)$.
\end{lemma}

\begin{proof}
Consider the one-parameter family of hypersurfaces
\[
M_t=\{x+t\chi_R(x)u(x)n_N(x):x\in N\cap B_R(0)\},
\qquad t\in[0,1].
\]
Define the weighted functional
\[
\mathcal{F}_f(M_t)
=
\int_{M_t} f e^{-\frac{|x|^2}{4}}\,d\mathcal{H}^n .
\]
By the first variation formula for the Gaussian area,
\[
\frac{d}{dt}\mathcal{F}_f(M_t)
=
\int_{M_t}
\chi_Ru
\left(
\partial_{n_{M_t}}f
-
f\phi_{M_t}
\right)
e^{-\frac{|x|^2}{4}}\,d\mathcal{H}^n .
\]
Hence, by the fundamental theorem of calculus,
\begin{multline*}
\int_{M_1} f e^{-\frac{|x|^2}{4}}\,d\mathcal{H}^n
-
\int_{N\cap B_R(0)} f e^{-\frac{|x|^2}{4}}\,d\mathcal{H}^n
\\
=
\int_0^1
\int_{M_t}
\chi_Ru
\left(
\partial_{n_{M_t}}f
-
f\phi_{M_t}
\right)
e^{-\frac{|x|^2}{4}}
\,d\mathcal{H}^n\,dt .
\end{multline*}
Observe that $M_1=M$ on $B_{R-1}(0)$, since
$\chi_R\equiv1$ there. Thus, by the triangle inequality,
\begin{align*}
&\left|
\int_{M\cap B_R(0)}
f e^{-\frac{|x|^2}{4}}\,d\mathcal{H}^n
-
\int_{M_1}
f e^{-\frac{|x|^2}{4}}\,d\mathcal{H}^n
\right| \\
&\leq
\|f\|_{C^0}
\left|
\int_{M\cap(B_R(0)\setminus B_{R-1}(0))}
e^{-\frac{|x|^2}{4}}\,d\mathcal{H}^n
-
\int_{N\cap(B_R(0)\setminus B_{R-1}(0))}
e^{-\frac{|x|^2}{4}}\,d\mathcal{H}^n
\right|
\\
&\quad
+
\|f\|_{C^0}
\left|
\int_{M_1\setminus B_{R-1}(0)}
e^{-\frac{|x|^2}{4}}\,d\mathcal{H}^n
-
\int_{N\cap(B_R(0)\setminus B_{R-1}(0))}
e^{-\frac{|x|^2}{4}}\,d\mathcal{H}^n
\right|.
\end{align*}
Since $u$ has sufficiently small $C^1$ norm, the induced volume forms on
$M_t$ and $N$ are uniformly equivalent. Therefore, arguing as in the proof
of Proposition~\ref{prop 5.2}, there exists a constant $C$ such that
\[
\left|
\int_{M\cap B_R(0)}
f e^{-\frac{|x|^2}{4}}\,d\mathcal{H}^n
-
\int_{M_1}
f e^{-\frac{|x|^2}{4}}\,d\mathcal{H}^n
\right|
\leq
C\|f\|_{C^0}
e^{-\frac{(R-1)^2}{4q}}.
\]
Using the fact that the induced volume forms on $M_t$ and $N$ are uniformly
equivalent, we get
\[
\int_{M_t}
|u\partial_{n_{M_t}}f|
e^{-\frac{|x|^2}{4}}\,d\mathcal{H}^n
\leq
C\|\nabla f\|_{C^0}
\int_{N\cap B_R(0)}
|u|e^{-\frac{|x|^2}{4}}\,d\mathcal{H}^n .
\]
By the Cauchy--Schwarz inequality and the entropy bound $\lambda(N)\leq\lambda_0$, there exists a constant $C$ such that
\[
\int_{M_t}
|u\partial_{n_{M_t}}f|
e^{-\frac{|x|^2}{4}}\,d\mathcal{H}^n
\leq
C\|\nabla f\|_{C^0}\|u\|_{L^2}.
\]
Moreover, using
\[
\|\phi_{M_t}\|_{L^2}
\leq
C\left(
\|\phi_N\|_{L^2}
+
\|u\|_{H^2}
\right),
\]
we obtain
\[
\int_{M_t}
|uf\phi_{M_t}|
e^{-\frac{|x|^2}{4}}\,d\mathcal{H}^n
\leq
C\|f\|_{C^0}
\|u\|_{L^2}
\left(
\|\phi_N\|_{L^2}
+
\|u\|_{H^2}
\right).
\]
Integrating these estimates with respect to $t\in[0,1]$ and combining them
with the estimate comparing $M$ and $M_1$ yields the desired inequality.
\end{proof}

As an application of Lemma~\ref{lemma 6.1}, we obtain the following estimate.

\begin{lemma}\label{lemma 6.2}
There exists a constant $C=C(\lambda_0)>0$ such that for any two
hypersurfaces $M_1$ and $M_2$ satisfying
$\lambda(M_1),\lambda(M_2)\leq\lambda_0$, we have
\[
d_V(\mu_{M_1},\mu_{M_2})
\leq
C\,D(M_1,M_2),
\]
where $\mu_{M_i}=\mathcal H^n\lfloor M_i$ denotes the Radon measure
associated to $M_i$.
\end{lemma}

\begin{proof}
Fix an admissible hypersurface $W\in\mathcal W$. We first establish the
following estimate. For any $f\in C^1(\mathbb{R}^{n+1})$ and any hypersurface $M$, we have
\[
\left|
\int_M f e^{-\frac{|x|^2}{4}}\,d\mathcal H^n
-
\int_W f e^{-\frac{|x|^2}{4}}\,d\mathcal H^n
\right|
\leq
C\|f\|_{C^1}D_W(M).
\]
Let $R>0$ be the largest radius such that $M$ is an $\eta$-graph over
$W$ on $B_R(0)$, given by a function $u$. By the triangle inequality,
\begin{align*}
&\bigg|
\int_M f e^{-\frac{|x|^2}{4}}\,d\mathcal{H}^n
-
\int_W f e^{-\frac{|x|^2}{4}}\,d\mathcal{H}^n
\bigg|
\\
&\leq
\bigg|
\int_{M\cap B_R(0)}
f e^{-\frac{|x|^2}{4}}\,d\mathcal{H}^n
-
\int_{W\cap B_R(0)}
f e^{-\frac{|x|^2}{4}}\,d\mathcal{H}^n
\bigg|
\\
&\quad+
\bigg|
\int_{M\setminus B_R(0)}
f e^{-\frac{|x|^2}{4}}\,d\mathcal{H}^n
-
\int_{W\setminus B_R(0)}
f e^{-\frac{|x|^2}{4}}\,d\mathcal{H}^n
\bigg|.
\end{align*}
By the definition of $D_W(M)$, Proposition~\ref{prop 5.2}, and the choice of
$q$, we have
\[
\begin{gathered}
\|u\|_{L^2}\leq CD_W(M),\\
\|u\|_{H^2}\leq C\eta,\\
\|\phi_W\|_{L^2}\leq C|\alpha|^2.
\end{gathered}
\]
Moreover,
\[
e^{-\frac{(R-1)^2}{4q}}
\leq
Ce^{-\frac{R^2}{4p_0}}
\leq
Ce^{-\frac{R^2}{8p_0}}
\leq
CD_W(M).
\]
Therefore, Lemma~\ref{lemma 6.1} gives
\begin{align*}
&\left|
\int_{M\cap B_R(0)}
f e^{-\frac{|x|^2}{4}}\,d\mathcal H^n
-
\int_{W\cap B_R(0)}
f e^{-\frac{|x|^2}{4}}\,d\mathcal H^n
\right| \\
&\qquad\leq
C\|\nabla f\|_{C^0}D_W(M)
+
C\|f\|_{C^0}
\left(
D_W(M)|\alpha|^2
+
D_W(M)\eta
+
D_W(M)
\right).
\end{align*}
Since $|\alpha|\leq\epsilon_0$ and $\eta$ is fixed, we obtain
\[
\left|
\int_{M\cap B_R(0)}
f e^{-\frac{|x|^2}{4}}\,d\mathcal H^n
-
\int_{W\cap B_R(0)}
f e^{-\frac{|x|^2}{4}}\,d\mathcal H^n
\right|
\leq
C\|f\|_{C^1}D_W(M).
\]
For the contribution outside $B_R(0)$, we have  
\[
\left|
\int_{M\setminus B_R(0)}
f e^{-\frac{|x|^2}{4}}\,d\mathcal H^n
-
\int_{W\setminus B_R(0)}
f e^{-\frac{|x|^2}{4}}\,d\mathcal H^n
\right|
\leq
C\|f\|_{C^0}e^{-\frac{R^2}{4q}}
\leq
C\|f\|_{C^0}D_W(M).
\]
Combining the above inequalities, we obtain
\[
\left|
\int_M f e^{-\frac{|x|^2}{4}}\,d\mathcal H^n
-
\int_W f e^{-\frac{|x|^2}{4}}\,d\mathcal H^n
\right|
\leq
C\|f\|_{C^1}D_W(M).
\]
Applying this estimate to both $M_1$ and $M_2$ and using the triangle inequality, we get
\[
\left|
\int_{M_1} f e^{-\frac{|x|^2}{4}}\,d\mathcal H^n
-
\int_{M_2} f e^{-\frac{|x|^2}{4}}\,d\mathcal H^n
\right|
\leq
C\|f\|_{C^1}
\bigl(D_W(M_1)+D_W(M_2)\bigr).
\]
Therefore,
\begin{align*}
d_V(\mu_{M_1},\mu_{M_2})
&=
\sum_{j=1}^{\infty}2^{-j}
\left|
\int_{\mathbb R^{n+1}}
f_j e^{-\frac{|x|^2}{4}}\,d\mu_{M_1}
-
\int_{\mathbb R^{n+1}}
f_j e^{-\frac{|x|^2}{4}}\,d\mu_{M_2}
\right|\\
&\leq
C\sum_{j=1}^{\infty}2^{-j}\|f_j\|_{C^1}
\bigl(D_W(M_1)+D_W(M_2)\bigr)\\
&\leq
C\bigl(D_W(M_1)+D_W(M_2)\bigr).
\end{align*}
Finally, taking the infimum over all admissible $W\in\mathcal W$ gives
\[
d_V(\mu_{M_1},\mu_{M_2})
\leq
CD(M_1,M_2),
\]
which completes the proof.
\end{proof}

We will also use the estimate for the one-parameter family introduced in
the proof of Lemma~\ref{lemma 6.1}. Namely, for
\[
M_t=\{x+t\chi_R(x)u(x)n_N(x):x\in N\cap B_R(0)\},
\qquad t\in[0,1],
\]
with sufficiently small $C^1$ norm, we have
\[
\|\phi_{M_t}\|_{L^2(M_t)}
\leq
C\left(
\|\phi_N\|_{L^2}
+\|Lu\|_{L^2}
+\|Q(u)\|_{L^2}
\right).
\]
Applying Lemma~\ref{lemma 6.1} with $f\equiv1$ and using the above
estimate, we obtain
\begin{align*}
&\left|
\int_{M\cap B_R(0)}
e^{-\frac{|x|^2}{4}}\,d\mathcal H^n
-
\int_{N\cap B_R(0)}
e^{-\frac{|x|^2}{4}}\,d\mathcal H^n
\right| \\
&\qquad\leq
C\|u\|_{L^2}
\left(
\|\phi_N\|_{L^2}
+\|Lu\|_{L^2}
+\|Q(u)\|_{L^2}
\right)
+
Ce^{-\frac{(R-1)^2}{4q}} .
\end{align*}

\begin{prop}\label{prop 6.7}
There exists $\epsilon>0$ such that if
$|\alpha|\leq\epsilon$, then for every $\tau\in[0,T_\alpha]$,
\[
\left|
F(\mathscr{T}^{\alpha,\mathscr C}_{\tau})
-
F(\mathscr C)
\right|
\leq
|\alpha|^{\frac52}.
\]
\end{prop}

\begin{proof}
Applying the above estimate with
\[
N=\mathscr C,\qquad
R=R_\alpha-1,\qquad
u=u_{\alpha(\tau)}+N(u_{\alpha(\tau)})+v(\tau),
\]
we obtain
\[
\left|
\int_{\mathscr{T}^{\alpha,\mathscr C}_{\tau}}
e^{-\frac{|x|^2}{4}}\,d\mathcal H^n
-
\int_{\mathscr C\cap B_{R_\alpha-1}(0)}
e^{-\frac{|x|^2}{4}}\,d\mathcal H^n
\right|
\leq
\frac12|\alpha|^{\frac52}
+
Ce^{-\frac{(R_\alpha-2)^2}{4q}}.
\]
Furthermore,
\[
\int_{\mathscr C\setminus B_{R_\alpha-1}(0)}
e^{-\frac{|x|^2}{4}}\,d\mathcal H^n
\leq
e^{-\frac{(R_\alpha-1)^2}{4q}}.
\]
Since $R_\alpha$ is chosen such that
\[
|\alpha|^{\frac94}e^{\frac{R_\alpha^2}{8}}=c_1,
\]
it follows that, for sufficiently small $|\alpha|$,
\[
Ce^{-\frac{(R_\alpha-2)^2}{4q}} + e^{-\frac{(R_\alpha-1)^2}{4q}}
\leq
\frac12|\alpha|^{\frac52}.
\]
Combining the above estimates and using the triangle inequality, we obtain
\begin{align*}
\left|
F(\mathscr{T}^{\alpha,\mathscr C}_{\tau})
-
F(\mathscr C)
\right|
&\leq
\left|
\int_{\mathscr{T}^{\alpha,\mathscr C}_{\tau}}
e^{-\frac{|x|^2}{4}}\,d\mathcal H^n
-
\int_{\mathscr C\cap B_{R_\alpha-1}(0)}
e^{-\frac{|x|^2}{4}}\,d\mathcal H^n
\right|
\\
&\qquad
+
\int_{\mathscr C\setminus B_{R_\alpha-1}(0)}
e^{-\frac{|x|^2}{4}}\,d\mathcal H^n
\\
&\leq
|\alpha|^{\frac52}.
\end{align*}
\end{proof}

We have seen in Lemma~\ref{lemma 4.8} that there is a second order obstruction in the Taylor expansion of $\phi$. We will now use this result to show that when the distance from $\mathscr {T}_\tau^{\alpha, \mathscr C}$ to $M_{\tau_0}$ is controlled at the scale $|\alpha|^2$, the decrease in the Gaussian area of the rescaled flow can be estimated by the corresponding decrease for $\mathscr{T}^{\alpha,\mathscr C}$.

For $\gamma > 0$, we define
\[
\mathcal A^\gamma(M)
:=
|\mathcal A(M)|^{\gamma-1}\mathcal A(M).
\]
Note that $\mathcal A^\gamma$ is decreasing along the rescaled flow.

\begin{prop} \label{prop 6.4}
There are $C, \kappa_1, \epsilon_0 > 0$ and $\gamma \in (0,1)$ with the following property. Suppose that $\tau + 1 \leq T_\alpha, |\alpha| \leq \epsilon_0$ and $\mathcal{D}_{\mathscr {T}^{\alpha, \mathscr{C}}_\tau} (M_{\tau_0}) \leq \epsilon_0.$ Suppose that we have the potentially better bound $\mathcal{D}_{\mathscr {T}^{\alpha, \mathscr{C}}_\tau} (M_{\tau_0}) \leq \kappa_1 |\alpha|^2.$ Then we have the following
$$\mathcal{A}^\gamma(M_{\tau_0}) - \mathcal{A}^\gamma(M_{\tau_0 + 1}) \geq C^{-1} |\alpha|^2.$$
\end{prop}

\begin{proof}
Let $\epsilon_0 >0$ be the constant from Proposition~\ref{prop 5.2}. Then by Proposition~\ref{prop 5.2}, if $\mathcal {D}_{\mathscr {T}^{\alpha, \mathscr{C}}_\tau} (M_{\tau_0}) = \epsilon \leq \epsilon_0,$ we can write $M_{\tau_0+s}$ as a graph over $\mathscr {T}^{\alpha, \mathscr{C}}_{\tau+s}$ on $B_{R_\epsilon}(0)$ of a function $u$, for $s \in [0,1]$ satisfying $|u| + |\nabla u| + |\nabla^2 u| \leq Ce^{\frac{|x|^2}{4p}} \epsilon$. Now, by Lemma~\ref{lemma 6.1},
\begin{align*}
\bigg|
&\int_{M_{\tau_0+s}\cap B_{R_\epsilon}(0)}
e^{-\frac{|x|^2}{4}}\,d\mathcal{H}^n
-
\int_{\mathscr{T}^{\alpha,\mathscr C}_{\tau+s}\cap B_{R_\epsilon}(0)}
e^{-\frac{|x|^2}{4}}\,d\mathcal{H}^n
\bigg| \\
&\leq C\left(
\|u\|_{L^2}\|\phi_{\mathscr{T}^{\alpha,\mathscr C}_{\tau+s}}\|_{L^2}
+\|u\|_{L^2}\|u\|_{H^2}
+e^{-\frac{(R_\epsilon-1)^2}{4q}}
\right) \\
&\leq C(\epsilon|\alpha|^2+\epsilon^2).
\end{align*}
Also,
\begin{align*}
\bigg|
&\int_{M_{\tau_0+s}\setminus B_{R_\epsilon}(0)}
e^{-\frac{|x|^2}{4}}\,d\mathcal{H}^n
-
\int_{\mathscr{T}^{\alpha,\mathscr C}_{\tau+s}\setminus B_{R_\epsilon}(0)}
e^{-\frac{|x|^2}{4}}\,d\mathcal{H}^n
\bigg| \\
&\leq
\int_{M_{\tau_0+s}\setminus B_{R_\epsilon}(0)}
e^{-\frac{|x|^2}{4}}\,d\mathcal{H}^n
+
\int_{\mathscr{T}^{\alpha,\mathscr C}_{\tau+s}\setminus B_{R_\epsilon}(0)}
e^{-\frac{|x|^2}{4}}\,d\mathcal{H}^n \\
&\leq Ce^{-\frac{R_\epsilon^2}{4p_0}}
\leq C\epsilon^2.
\end{align*}
By our assumption, $\epsilon \leq \kappa_1|\alpha|^2$. Hence, after enlarging $C$ if necessary,
$$|F(M_{\tau_0+s}) - F(\mathscr {T}^{\alpha, \mathscr{C}}_{\tau+s})| \leq C \kappa_1 |\alpha|^4.$$
By Lemma~\ref{lemma 4.8} and Huisken's monotonicity formula, we have
\begin{align*}
    \mathcal{A}(\mathscr {T}^{\alpha, \mathscr{C}}_{\tau}) - \mathcal{A}(\mathscr {T}^{\alpha, \mathscr{C}}_{\tau +1}) &= \int_0^{1} \int_{\mathscr {T}^{\alpha, \mathscr{C}}_{\tau+s}} \phi^2 e^{-\frac{|x|^2}{4}} d\mathcal{H}^n \\
    &\geq \kappa_0^2 |\alpha|^4.
\end{align*}
Therefore, 
\begin{align*}
    \mathcal{A}(M_{\tau_0}) - \mathcal{A}(M_{\tau_0+1}) &\geq \mathcal{A}(\mathscr {T}^{\alpha, \mathscr{C}}_{\tau}) - \mathcal{A}(\mathscr {T}^{\alpha, \mathscr{C}}_{\tau +1}) - 
    |\mathcal{A}(M_{\tau_0}) - \mathcal{A}(\mathscr {T}^{\alpha, \mathscr{C}}_{\tau})| \\ &\quad- |\mathcal{A}(M_{\tau_0+1}) - \mathcal{A}(\mathscr {T}^{\alpha, \mathscr{C}}_{\tau+1})|\\
    &\geq \kappa_0^2 |\alpha|^4 - 2C \kappa_1 |\alpha|^4 \\
    &\geq \frac{\kappa_0^2}{2} |\alpha|^4.
\end{align*}
where we have chosen $\kappa_1$ small such that $\kappa_1 < \frac{\kappa_0^2}{4C}.$ Moreover, for $s\in[0,1]$, we have
\[
|F(M_{\tau_0+s})-F(\mathscr{T}^{\alpha,\mathscr{C}}_{\tau_0+s})|
\leq C|\alpha|^4,
\]
and by Proposition~\ref{prop 6.7},
\[
|F(\mathscr{T}^{\alpha,\mathscr{C}}_{\tau_0+s})-F(\mathscr{C})|
\leq |\alpha|^{\frac52}.
\]
Therefore, by the triangle inequality, for sufficiently small $|\alpha|$,
\[
|\mathcal{A}(M_{\tau_0+s})|
=
|F(M_{\tau_0+s})-F(\mathscr{C})|
\leq
2|\alpha|^{\frac52}.
\]
Now, by Huisken's monotonicity formula, for $s\in[0,1]$ and
$\gamma\in(0,1)$,
\[
\frac{d}{ds}\mathcal{A}(M_{\tau_0+s})^\gamma
=
\gamma |\mathcal{A}(M_{\tau_0+s})|^{\gamma-1}
\frac{d}{ds}\mathcal{A}(M_{\tau_0+s}).
\]
Integrating with respect to $s\in[0,1]$, we obtain
\begin{align*}
\mathcal{A}(M_{\tau_0})^\gamma
-\mathcal{A}(M_{\tau_0+1})^\gamma
&=
\gamma\int_0^1
|\mathcal{A}(M_{\tau_0+s})|^{\gamma-1}
\left(-\frac{d}{ds}\mathcal{A}(M_{\tau_0+s})\right)\,ds\\
&\geq
\gamma(2|\alpha|^{5/2})^{\gamma-1}
\int_0^1
\left(-\frac{d}{ds}\mathcal{A}(M_{\tau_0+s})\right)\,ds\\
&=
\gamma 2^{\gamma-1}
|\alpha|^{\frac{5\gamma}{2}-\frac52}
\left(
\mathcal{A}(M_{\tau_0})
-\mathcal{A}(M_{\tau_0+1})
\right) \\
&\geq \gamma 2^{\gamma-1} |\alpha|^{\frac{5\gamma}{2}-\frac52} \frac{\kappa_0^2}{2}|\alpha|^4 \\
&= \gamma \kappa_0^2 2^{\gamma -2}|\alpha|^{\frac{5\gamma}{2} - \frac12} |\alpha|^2.
\end{align*}
Choosing $0<\gamma<\frac15,$ we have $|\alpha|^{\frac{5\gamma}{2} - \frac12} \geq \epsilon_0^{\frac{5\gamma}{2} - \frac12}$. Therefore, this factor, together with the remaining constants, can be absorbed into $C^{-1}$, and the required result follows.
\end{proof}

For a function $u\in C^2(M)$ with sufficiently small $C^2$ norm, we
write
\[
M_u=\{x+u(x)\nu_M(x):x\in M\}
\]
for the normal graph of $u$ over $M$. For a fixed reference hypersurface $N$, we define the Gaussian weighted
$L^2$ distance between $M$ and $N$ by
\[
D_N^2(M)
=
\int_M d_N^2 e^{-\frac{|x|^2}{4}}\,d\mathcal{H}^n,
\]
where $d_N$ denotes the distance function to $N$. The following lemma estimates the Gaussian weighted distance to N for two hypersurfaces that are realized as sufficiently small normal graphs over the same reference hypersurface.

\begin{lemma} \label{lemma 6.5}
Let $M\subset \mathbb{R}^{n+1}$ be a smooth hypersurface satisfying $|A_M|\leq C_0,$ and $|\langle x,\nu_M\rangle|\leq C_0 .$ Suppose that $u,v\in C^2(M)$ satisfy $\|u\|_{C^2},\|v\|_{C^2}\leq \epsilon$, where $\epsilon$ sufficiently small. Moreover, assume that $N$ can be written as a normal graph over both $M_u$ and $M_v$. Then there exists a constant $C$ such that
$$D_N(M_u) \leq CD_N(M_v) + C\|u-v\|_{L^2(M)}. $$
\end{lemma}

\begin{proof}
Let
\[
X_w(x)=x+w(x)\nu_M(x)
\]
be the normal graph parametrization of $M_w$ over $M$, and let $J_w$
denote the corresponding Jacobian.

Since the distance function is $1$-Lipschitz, we have
\[
|d_N(X_u)-d_N(X_v)|
\leq
|X_u-X_v|.
\]
Since $X_u-X_v=(u-v)\nu_M,$ we obtain
\[
d_N^2(X_u)
\leq
2d_N^2(X_v)+2|u-v|^2.
\]
Moreover, for any sufficiently small graph function $w$, using $|\langle x,\nu_M\rangle|\leq C_0$ as in Proposition~\ref{prop 5.2} we have
\[
e^{-|X_w(x)|^2/4}
\leq
Ce^{-|x|^2/4}.
\]
Since $J_w$ depends smoothly on $w,\nabla w$ and
$\|u\|_{C^2},\|v\|_{C^2}\leq\epsilon$, the Jacobians are uniformly
comparable. In particular, $J_u\leq CJ_v$ and $J_u\leq C.$ Therefore,
\begin{align*}
D_N^2(M_u)
&=
\int_M
d_N^2(X_u)
e^{-|X_u|^2/4}
J_u\,d\mathcal{H}^n
\\
&\leq
2\int_M
d_N^2(X_v)
e^{-|X_u|^2/4}
J_u\,d\mathcal{H}^n
+
2\int_M
|u-v|^2
e^{-|X_u|^2/4}
J_u\,d\mathcal{H}^n
\\
&\leq
C\int_M
d_N^2(X_v)
e^{-|x|^2/4}
J_v\,d\mathcal{H}^n
+
C\int_M
|u-v|^2
e^{-|x|^2/4}\,d\mathcal{H}^n
\\
&=
CD_N^2(M_v)
+
C\int_M
|u-v|^2
e^{-|x|^2/4}\,d\mathcal{H}^n.
\end{align*}
So, the lemma follows.
\end{proof}

Combining the previous lemma with Lemma~\ref{lemma 7.2}, we obtain the following estimates.

\begin{lemma}\label{lemma 6.6}
There exist constants $C>0$ and $\epsilon>0$, with $C$ depending on the entropy bound, such that the following hold. Suppose that
\[
|\alpha|,\ |\widehat{\alpha}|,\ |\mathscr{R}-Id|,\ 
D_{\mathscr{T}^{\alpha,\mathscr{C}}_\tau}(M_{\tau_0}),\
D_{\mathcal{W}^{\alpha(\tau),\mathscr{C}}}(M_{\tau_0}), 
D_{\mathscr{T}^{\alpha,\mathscr{C}}_{\tau+L}}(M_{\tau_0+L})
\leq \epsilon .
\]
Then
\[
D_{\mathscr{T}^{\widehat{\alpha},\mathscr{RC}}_{\tau+1}}(M_{\tau_0+1})
\leq
C\bigg(
D_{\mathscr{T}^{\alpha,\mathscr{C}}_\tau}(M_{\tau_0})
+|\mathscr{R}-Id|
+|\alpha-\widehat{\alpha}|
+e^{-\frac{(R_{\widehat{\alpha}}-1)^2}{8p_0}}
\bigg).
\]
For relating the distance from $\mathcal{W}^{\alpha,\mathscr{C}}$ to the distance from $\mathscr{T}^{\alpha,\mathscr{C}}$, we have
\[
D_{\mathcal{W}^{\widehat{\alpha}(\tau+1),\mathscr{RC}}}(M_{\tau_0+1})
\leq
C\bigg(
D_{\mathscr{T}^{\alpha,\mathscr{C}}_\tau}(M_{\tau_0})
+|\mathscr{R}-Id|
+|\alpha-\widehat{\alpha}|
+e^{-\frac{(R_{\widehat{\alpha}}-1)^2}{8p_0}}
\bigg).
\]
Also,
\[
D_{\mathscr{T}^{\widehat\alpha,\mathscr{RC}}_\tau}(M_{\tau_0})
\leq
C\bigg(
D_{\mathcal{W}^{\alpha(\tau),\mathscr{C}}}(M_{\tau_0})
+|\mathscr{R}-Id|
+|\alpha-\widehat{\alpha}|
+e^{-\frac{(R_{\widehat{\alpha}}-1)^2}{8p_0}}
\bigg).
\]
In the case $|\widehat{\alpha}|\geq |\alpha|$, the exponential terms in the above estimates are controlled by the corresponding distance terms and therefore can be omitted.

Finally,
\[
D_{\mathcal{W}^{\alpha(\tau+1),\mathscr{C}}}(M_{\tau_0+1})
+
D_{\mathcal{W}^{\alpha(\tau+1),\mathscr{C}}}(M_{\tau_0+L+1})
\leq
C\Big(
D_{\mathscr{T}^{\alpha,\mathscr{C}}_\tau}(M_{\tau_0})
+
D_{\mathscr{T}^{\alpha,\mathscr{C}}_{\tau+L}}(M_{\tau_0+L})
+
|\alpha|^2
\Big).
\]
\end{lemma}

\begin{proof}
By definition of $D_{\mathscr {T}^{\alpha, \mathscr{C}}_\tau}(M_{\tau_0})$, there exists $R \leq R_\alpha-1$ such that $M_{\tau_0}$ is $\eta$-graph over $\mathscr {T}^{\alpha, \mathscr{C}}_{\tau}$ on $B_R(0)$ and 
$$a := D_{\mathscr {T}^{\alpha, \mathscr{C}}_\tau}(M_{\tau_0})= \bigg(\int_{M_{\tau_0} \cap B_{R}(0)} d_{\mathscr{T}^{\alpha, \mathscr{C}}_\tau}^2 \;e^{-\frac{|x|^2}{4}} d \mathcal{H}^n\bigg)^{\frac{1}{2}} + e^{-\frac{R^2}{8p_0}}.$$ 
By Proposition~\ref{prop 5.2}, provided that $|\alpha|$ and $a$ are sufficiently small, $M_{\tau_0+1}$ is an $\frac{\eta}{4}$-graph over $\mathscr {T}^{\alpha, \mathscr{C}}_{\tau+1}$ on $B_{R_a}(0),$ and 
$$\bigg(\int_{M_{\tau_0 + 1} \cap B_{R_a}(0)} d_{\mathscr {T}^{\alpha, \mathscr{C}}_{\tau +1}}^{2p} \;e^{-\frac{|x|^2}{4}} d \mathcal{H}^n\bigg)^{\frac{1}{2p}} \leq C a.$$
By Lemma ~\ref{lemma 7.2}, $M_{\tau_0 +1}$ is an $\frac{\eta}{2}$-graph over $\mathscr{T}^{\widehat{\alpha}, \mathscr{RC}}_{\tau+1}$ on $B_{\min\{R_a, R_{\widehat \alpha}-1\}}(0)$ when $|\widehat{\alpha}|, |\mathscr{R} - Id|$ are sufficiently small. By definition,
$$D_{\mathscr{T}^{\widehat{\alpha}, \mathscr{RC}}_{\tau+1}}(M_{\tau_0+1})= \bigg(\int_{M_{\tau_0+1} \cap B_{\widehat{R}}(0)} d_{\mathscr{T}^{\widehat{\alpha}, \mathscr{RC}}_{\tau+1}}^2 \;e^{-\frac{|x|^2}{4}} d \mathcal{H}^n\bigg)^{\frac{1}{2}} + e^{-\frac{\widehat{R}^2}{8p_0}}$$
where $\widehat{R}$ denotes the largest radius such that $M_{\tau_0+1}$ is an $\eta$-graph over $\mathscr{T}^{\widehat{\alpha}, \mathscr{RC}}_{\tau+1}$. Hence, 
$$\widehat{R} \geq \min\{R_a, R_{\widehat \alpha}-1\}.$$ 
Consequently,
\[
e^{-\frac{\widehat{R}^2}{8p_0}}
\le
\max\left\{
e^{-\frac{R_a^2}{8p_0}},
\,
e^{-\frac{(R_{\widehat{\alpha}}-1)^2}{8p_0}}
\right\}.
\]
By Lemma ~\ref{lemma 6.5}, 
\begin{multline*}
\bigg(
\int_{M_{\tau_0+1} \cap B_{R_a}(0)}
d_{\mathscr{T}^{\widehat{\alpha},\mathscr{RC}}_{\tau+1}}^2
e^{-\frac{|x|^2}{4}}
\,d\mathcal{H}^n
\bigg)^{\frac{1}{2}}
\\
\leq
C
\bigg(
\int_{M_{\tau_0+1}\cap B_{R_a}(0)}
d_{\mathscr{T}^{\alpha,\mathscr{C}}_{\tau+1}}^{2p}
e^{-\frac{|x|^2}{4}}
\,d\mathcal{H}^n
\bigg)^{\frac{1}{2p}}
+
C
\bigg(
|\mathscr{R}-Id|
+
|\alpha-\widehat{\alpha}|
+
e^{-\frac{(R_{\alpha}-1)^2}{8p_0}}
+
e^{-\frac{(R_{\widehat{\alpha}}-1)^2}{8p_0}}
\bigg).
\end{multline*}
By the entropy bound,
$$\bigg(\int_{M_{\tau_0+1} \cap B_{\widehat{R}}(0) \setminus B_{R_a}(0)} d_{\mathscr{T}^{\widehat{\alpha}, \mathscr{RC}}_{\tau+1}}^2 \;e^{-\frac{|x|^2}{4}} d \mathcal{H}^n\bigg)^{\frac{1}{2}} \leq C \eta e^{-\frac{R_a ^2}{8p_0}}.$$
As $R \leq R_\alpha-1,$
$$e^{-\frac{(R_\alpha-1)^2}{8p_0}} \leq e^{-\frac{R^2}{8p_0}} \leq a.$$
Combining these estimates proves the first inequality. 

We can follow the same argument considering the surface $\mathcal{W}^{\widehat{\alpha}(\tau+1),\mathscr{RC}}$ in place of $\mathscr{T}^{\widehat{\alpha},\mathscr{RC}}_{\tau+1}$ to get 
\[
D_{\mathcal{W}^{\widehat{\alpha}(\tau+1),\mathscr{RC}}}(M_{\tau_0+1})
\leq
C\bigg(
D_{\mathscr{T}^{\alpha,\mathscr{C}}_\tau}(M_{\tau_0})
+|\widehat{\alpha}|^{\frac52}
+|\mathscr{R}-Id|
+|\alpha-\widehat{\alpha}|
+e^{-\frac{(R_{\widehat{\alpha}}-1)^2}{8p_0}}
\bigg).
\]
Now
$$e^{-\frac{(R_{\widehat{\alpha}}-1)^2}{8p_0}} \geq e^{-\frac{R_{\widehat{\alpha}}^2}{8p_0}} = c_1^{-\frac{1}{p_0}} |\widehat{\alpha}|^{\frac{9}{4p_0}} \geq c_1^{-\frac{1}{p_0}} |\widehat{\alpha}|^{\frac52}$$
So the second inequality follows.

Now by Lemma ~\ref{lemma 7.2},
\[
D_{\mathscr{T}^{\alpha,\mathscr{C}}_\tau}(M_{\tau_0})
\leq
C\bigg(
D_{\mathcal{W}^{\alpha(\tau),\mathscr{C}}}(M_{\tau_0})
+|\alpha|^{\frac52}
+|\mathscr{R}-Id|
+|\alpha-\widehat{\alpha}|
+e^{-\frac{(R_{\widehat{\alpha}}-1)^2}{8p_0}}
\bigg).
\]
Again using 
$$D_{\mathcal{W}^{\alpha(\tau),\mathscr{C}}}(M_{\tau_0}) \geq e^{-\frac{(R_{{\alpha}}-1)^2}{8p_0}} \geq c_1^{-\frac{1}{p_0}} |{\alpha}|^{\frac52},$$
the third inequality follows.

Now let 
$$b := \max \{D_{\mathscr{T}^{\alpha,\mathscr{C}}_\tau}(M_{\tau_0}), D_{\mathscr{T}^{\alpha,\mathscr{C}}_{\tau+L}}(M_{\tau_0+L})\}.$$ 
Then by definition both $M_{\tau_0}, M_{\tau_0 +L}$ are $\eta$-graph over $\mathscr {T}^{\alpha, \mathscr{C}}_{\tau}, \mathscr {T}^{\alpha, \mathscr{C}}_{\tau + L}$ on $B_{R_b}(0).$ So, as before, using Proposition ~\ref{prop 5.2}, $M_{\tau_0+1}$ and $M_{\tau_0+L+1}$ can be written as $\frac{\eta}{4}$-graphs over $\mathcal{W}^{\alpha(\tau+1),\mathscr{C}}$ and $\mathcal{W}^{\alpha(\tau +L +1),\mathscr{C}}$ respectively on $B_{R_b}(0)$. 

By Lemma~\ref{lemma 7.2}, $M_{\tau_0+ L+1}$ can be written as a $\frac{\eta}{2}$-graphs over $\mathcal{W}^{\alpha(\tau+1),\mathscr{C}}$ on $B_{R_b}(0)$. Now, applying Lemma~\ref{lemma 6.5} together with the estimate $|\alpha(\tau) - \alpha(0)| \leq C|\alpha|^2$, we obtain
$$D_{\mathcal{W}^{\alpha(\tau+1),\mathscr{C}}}(M_{\tau_0+L+1}) \leq C \bigg(D_{\mathcal{W}^{\alpha(\tau+ L+1),\mathscr{C}}}(M_{\tau_0+L+1}) + |\alpha|^2 + e^{-\frac{R_b^2}{8p_0}}\bigg).$$
Using the second inequality we obtain
$$D_{\mathcal{W}^{\alpha(\tau+1),\mathscr{C}}}(M_{\tau_0+1}) \leq C D_{\mathscr{T}^{\alpha,\mathscr{C}}_\tau}(M_{\tau_0}),$$
and
$$D_{\mathcal{W}^{\alpha(\tau+ L+1),\mathscr{C}}}(M_{\tau_0+L+1}) \leq CD_{\mathscr{T}^{\alpha,\mathscr{C}}_{\tau+L}}(M_{\tau_0+L}).$$
Combining the above estimates proves the last inequality.
\end{proof}

Now we will characterize shrinkers near $\mathscr{C}$. We will see that cylinders are rigid in the space of shrinkers. Similar result has been proved in [\citealp{colding2015rigidity}, \citealp{zhu2020ojasiewicz}]. The key point is that $\mathscr{C}$ is not integrable in the sense that there does not exist a family of shrinkers near $\mathscr{C}$ that can be modeled on the Jacobi field in $K_2$. The proof follows a similar approach to [\citealp{szekelyhidi2020uniqueness}, Proposition ~6.5].

\begin{prop} \label{prop 6.7}
    There exists $\epsilon > 0$ such that the following holds. Suppose that $M$ is a shrinker in $\mathbb{R}^{n+1}$, such that $d_V(\mu_{\mathscr{C}}, \mu_M) \leq \epsilon.$ Then $M$ is a rotation of $\mathscr{C}$.
\end{prop}

\begin{proof}
We will argue by contradiction. Let $M^i$ be a sequence of shrinkers with $d_V(\mu_{\mathscr{C}}, \mu_{M^i}) \to 0$ such that the $M^i$ are not rotations of $\mathscr{C}$. For any $i$, if
$$\inf \limits_{\alpha} \Big\{D_{\mathcal{W}^{\alpha, \mathscr{RC}}}(M^i) :  |\alpha| \leq \epsilon_1 \Big \} = 0,$$
then we must have $M^i = \mathcal{W}^{\alpha, \mathscr{RC}}$ for some $\alpha$ and rotation $\mathscr{R}$, but the surfaces $\mathcal{W}^{\alpha, \mathscr{RC}}$ are not shrinkers for $|\alpha| \neq 0$, so $M^i = \mathscr{RC}$. Therefore, we can assume that the infimum above is positive for all $i$. We choose a sequence of rotations $\mathscr{R}^i$ and $\alpha^i$ with $|\alpha^i| \to 0$ such that
$$D_i := D_{\mathcal{W}^{\alpha^i, \mathscr{R}^i \mathscr{C}}}(M^i) < 2 \inf \limits_{\mathscr{R},\alpha} \Big \{D_{\mathcal{W}^{\alpha, \mathscr{RC}}}(M^i): |\alpha| \leq \epsilon_1 \Big\}.$$
By Allard’s regularity theorem ~\cite{allard1972first}, the convergence $d_V(\mu_{\mathscr{C}}, \mu_{M^i}) \to 0$ implies that $M^i$ converges locally smoothly to $\mathscr{C}$. Therefore we have $D_i \to 0$. By Lemma ~\ref{lemma 6.6},
$$D_{\mathscr{T}^{\alpha^i, \mathscr{R}^i \mathscr{C}}_0}(M^i) \leq CD_i.$$
Using the argument in the proof of Proposition~\ref{prop 5.2}, we obtain that for any $\tau \in [1,3]$, 
$$D_{\mathscr{T}^{\alpha^i, \mathscr{R}^i \mathscr{C}}_\tau}(M^i) \leq CD_i.$$
This implies  
$$\mathcal{D}_{\mathscr{T}^{\alpha^i, \mathscr{R}^i \mathscr{C}}_3}(M^i) \leq CD_i.$$
for possibly a larger constant $C.$ We can write $M^i$ as the graph of $u_i$ over larger and larger region of $\mathscr{T}_\tau^{\alpha^i, \mathscr{R}^i \mathscr{C}}$ for $\tau \in [0,5]$. Choosing a subsequence, we can assume that $D^{-1}_i u_i \to u_0$ locally uniformly on $\mathscr{C}$ and $u_0$ is a Jacobi field on $\mathscr{C}.$ So, by Lemma ~\ref{lemma 2.1}, we can write  
$$u_0 = f + \sum_{j=1}^N \beta_j U_j,$$
where $f$ is the Jacobi field corresponding to a rotation. Let $\widehat{\mathscr{R}}^i$ denote the rotation corresponding to the Jacobi field $D_i f$, and let $\widehat \alpha^i = \alpha^i + D_i \beta,$ where $\beta = (\beta_1, \beta_2, ... , \beta_N).$ We can also view $M^i$ as graphs of functions $\widehat u_i$ over larger and larger regions of  $\mathscr{T}_\tau^{\widehat \alpha^i, \widehat {\mathscr{R}}^i \mathscr{R}^i \mathscr{C}}$ for $\tau \in [0,5]$. Denote 
$$\widehat D_i := D_{\mathscr{T}_4^{\widehat \alpha^i, \widehat {\mathscr{R}}^i \mathscr{R}^i \mathscr{C}}}(M^i).$$
So by Lemma~\ref{lemma 6.6}, for large $i$, we have 
$$\widehat D_i \leq C \bigg(D_i + e^{-\frac{(R_{\widehat \alpha^i}-1)^2}{8p_0}}\bigg).$$
We claim that, for sufficiently large $i$, 
$$D_i \geq e^{-\frac{(R_{\widehat{\alpha}^i}-1)^2}{8p_0}}.$$ 
Indeed, if this were false, then
$$\widehat D_i \leq 2Ce^{-\frac{(R_{\widehat \alpha^i}-1)^2}{8p_0}}.$$
Proposition ~\ref{prop 5.2} gives
$$\mathcal{D}_{\mathscr{T}_4^{\widehat \alpha^i, \widehat {\mathscr{R}}^i \mathscr{R}^i \mathscr{C}}}(M^i) \leq 2C e^{-\frac{(R_{\widehat\alpha^i}-1)^2}{8p_0}}.$$
Let $\kappa_1$ be the constant from Proposition ~\ref{prop 6.4}. By the definition of $R_{\widehat\alpha^i}$ $$|\widehat\alpha^i|^{\frac94}e^{\frac{R_{\widehat\alpha^i}^2}{8}}=c_1,$$
Since $p_0 \in (1, \frac{9}{8}),$ this implies that, for sufficiently small $|\widehat\alpha^i|$, $$\mathcal{D}_{\mathscr{T}_4^{\widehat \alpha^i, \widehat {\mathscr{R}}^i \mathscr{R}^i \mathscr{C}}}(M^i)  \leq 2C e^{-\frac{(R_{\widehat\alpha^i}-1)^2}{8p_0}} \leq \kappa_1|\widehat\alpha^i|^2.$$
Since $M^i$ is a shrinker, it gives a stationary rescaled flow. Hence, Proposition~\ref{prop 6.4} yields
$$\mathcal{A}^\gamma(M^i_{\tau_0}) - \mathcal{A}^\gamma(M^i_{\tau_0+1}) \geq C^{-1} |\widehat\alpha^i|^2,$$
This contradicts the fact that $M^i$ is a stationary rescaled flow. Therefore, after enlarging $C$ if necessary,
$$\widehat D_i \leq C D_i.$$ 
Moreover, by construction and Lemma ~\ref{lemma 6.6} we have $D_i^{-1} \widehat u_i \to 0$ for $\tau \in [1,5].$ We can apply the non-concentration estimate Corollary \ref{cor 5.3}. It follows that for any $\delta > 0$ we have 
$$\widehat D_i < \delta D_i$$ 
for sufficiently large $i$. So we can use Lemma~\ref{lemma 6.6} to conclude 
$$D_{\mathcal{W}^{\widehat \alpha^i(4), \widehat {\mathscr{R}}^i \mathscr{R}^i \mathscr{C}}}(M^i) \leq C\delta D_i.$$
This contradicts the definition of $D_i$. 
\end{proof}

We introduce the following notation for the distance between two time slices of the rescaled flow:
$$\mathcal {D} (M_{\tau_1},M_{\tau_2}) := \sup_{s \in [0,1]} D(M_{\tau_1 -s}, M_{\tau_2 -s}).$$
This quantity measures the maximal distance between the two corresponding time intervals of the rescaled flow. We now proceed to prove the main result, which leads to the uniqueness of the tangent flow, analogous to [\citealp{szekelyhidi2020uniqueness}, Proposition ~6.6], as follows.

\begin{prop} \label{prop 6.8}
There are $C, L, \epsilon > 0$ and $\gamma \in (0,1)$ with the following property. Suppose that $\sup\limits_{s\in [0,1]} d_V(\mu_{\mathscr{C}}, \mu_{M_{\tau_0-s}}) \leq \epsilon,$ $\mathcal A^\gamma (M_{\tau_0 - 1}) - \lim\limits_{\tau \to \infty} \mathcal A^\gamma(M_\tau) \leq \epsilon$.
Then one of the following holds: 
\begin{enumerate} [label=\normalfont(\roman*) ]
    \item $\mathcal{D}(M_{\tau_0 + L}, M_{\tau_0 + 2L}) \leq \frac{1}{2} \mathcal{D}(M_{\tau_0}, M_{\tau_0 + L}).$
    \item $\mathcal{D}(M_{\tau_0 + L}, M_{\tau_0 + 2L}) \leq C \Big(\mathcal A^\gamma(M_{\tau_0 + L}) -\mathcal A^\gamma(M_{\tau_0 + 2L}) \Big).$
\end{enumerate}
\end{prop}

\begin{proof}
Suppose, by contradiction, that the conclusion fails. Then there exists a sequence of flows $M^i$ and a sequence of cylinders $\mathscr{C}^i$ such that $\sup\limits_{s\in [0,1]} d_V(\mu_{\mathscr{C}^i}, \mu_{M^i_{\tau_0-s}})$, $\mathcal{A}(M^i_{\tau_0 - 1}) - \lim\limits_{\tau \to \infty} \mathcal{A}(M^i_\tau) \leq \epsilon_i$ where $\epsilon_i \to 0$, but neither of the two conclusions holds. We will derive a contradiction, provided that $L$ is chosen sufficiently large. First note that, up to choosing a subsequence, we can replace the sequence $\mathscr{C}^i$ by a single $\mathscr{C}.$

Let $\tilde{\delta_1},\tilde{\delta_2}$ be the constants from Lemma ~\ref{lemma 5.6}. We first claim that, if $L$ is chosen sufficiently large, then along a subsequence we can find $|\alpha^i|, \epsilon_i \to 0$ and rotations $\mathscr{R}^i$ converging to identity such that one of the following possibilities holds: 
\begin{enumerate} [label=\normalfont(\alph*)]
    \item $\mathcal{D}_{\mathscr{T}^{\alpha^i, \mathscr{R}^i \mathscr{C}}_{1+2L}} (M^i_{\tau_0 + 2L}) \geq e^{\tilde{\delta_1} L} \mathcal{D}_{\mathscr{T}^{\alpha^i, \mathscr{R}^i \mathscr{C}}_{1+L}} (M^i_{\tau_0 + L}),$ \vspace{0.15 cm}
    \item $\mathcal{D}_{\mathscr{T}^{\alpha^i, \mathscr{R}^i \mathscr{C}}_2} (M^i_{\tau_0 + 1}) \geq e^{\tilde{\delta_2} L} \mathcal{D}_{\mathscr{T}^{\alpha^i, \mathscr{R}^i \mathscr{C}}_{1+L}} (M^i_{\tau_0 + L}).$ 
\end{enumerate}
We can also assume that $M^i$ does not satisfy any of the given conditions.
Let us first choose sequences $\mathscr{R}^i, \alpha^i$ such that
$$\mathcal{D}_{\mathscr{T}^{\alpha^i, \mathscr{R}^i \mathscr{C}}_{1+L}} (M^i_{\tau_0 + L}) < 2 \;\text{inf} \Big\{\mathcal{D}_{\mathscr{T}^{\alpha, \mathscr{RC}}_{1+L}} (M^i_{\tau+L})\Big\},$$
where the infimum is over all rotations $\mathscr{R}$ and $\alpha$ for which $\mathscr{T}^{\alpha, \mathscr{R}^\lambda \mathscr{C}}$ is defined on the time interval $[0, T_\alpha]$, i.e. $\mathscr{T}^{\alpha^i, \mathscr{R}^i \mathscr{C}}$ is approximately a best rescaled flow to the rescaled flow $M^i$ among the family of comparison flows on the time interval $[0, T_{\alpha^i}]$. Suppose that both $(\textnormal{a})$ and $(\textnormal{b})$ fail for the sequence $\alpha^i$ and $\mathscr{R}^i$. Letting 
$$\mathscr{C}^i = \mathscr{R}^i \mathscr{C},$$ 
and 
$$d_i = \mathcal{D}_{\mathscr{T}^{\alpha^i, \mathscr{R}^i \mathscr{C}}_{1+L}} (M^i_{\tau_0 + L}),$$ 
we have $d_i \to 0$ and since neither ~(a) nor ~(b) holds, we have
$$\mathcal{D}_{\mathscr{T}^{\alpha^i, \mathscr{C}^i}_{1+2L}} (M^i_{\tau_0 + 2L}) < e^{\tilde{\delta_1} L} d_i,$$
and
$$\mathcal{D}_{\mathscr{T}^{\alpha^i, \mathscr{C}^i}_2} (M^i_{\tau_0+1}) < e^{\tilde{\delta_2} L} d_i,$$
As in the proof of Lemma ~\ref{lemma 5.6}, we can find $M^i_{\tau_0 + s}$ as the graph of $u_i(s)$ over $\mathscr{T}_{1+s}^{\alpha^i, \mathscr{C}^i}$ on larger and larger subsets of $\mathscr{T}_{1+s}^{\alpha^i, \mathscr{C}^i}$ for $s \in [-1, 2L]$. In addition, up to choosing a subsequence, the rescaled functions $d_i^{-1} u_i$ converge locally smoothly to a solution $u$ of the drift heat equation on $\mathscr{C}.$ Note that as in Lemma ~\ref{lemma 5.6} we have,
$$I(u,0) \leq Ce^{\tilde{\delta_2} L},$$
and
$$I(u,2L) \leq Ce^{\tilde{\delta_1} L}.$$
for a constant $C$ independent of $L$.

Let $u_0$ be the degree zero homogeneous Jacobi field, such that $u - u_0$ has no degree zero piece. By Lemma~\ref{lemma 4.1} we can write
$$u_0 = f + \sum_{j=1}^N \beta_j U_j,$$
where $f$ is the Jacobi field corresponding to a rotation. Define $\beta \in \mathbb{R}^N$ as $\beta = (\beta_1, \beta_2, ... , \beta_N).$ Let $\widehat{\alpha}^i \in \mathbb{R}^N$ be defined by $\widehat{\alpha}^i = \alpha^i + d_i\beta$ and let $\widehat{\mathscr{R}}^i$ denote the rotation corresponding to the Jacobi field $d_i f$. By construction, $M^i_{\tau_0+s}$ are graphs of new functions $\widehat u_i$ over $\mathscr{T}^{\widehat \alpha_i, \widehat{\mathscr{R}}^i \mathscr{C}^i}_{1+s}$ on larger and larger subsets of $\mathscr{T}^{\widehat\alpha^i, \widehat{\mathscr{R}}^i \mathscr{C}^i}_{1+s}$ for $s \in [0, 2L]$. By Lemma ~\ref{lemma 6.6}, for large $i$ and for $s \in [1,2L],$ we have
$$\mathcal{D}_{\mathscr{T}^{\widehat \alpha^i, \widehat{\mathscr{R}}^i\mathscr{C}^i}_{1+s}}(M^i_{\tau_0 +s}) \leq C\bigg(d_i + e^{-\frac{(R_{\widehat \alpha^i}-1)^2}{8p_0}}\bigg).$$ 
Note that, as in Proposition ~\ref{prop 6.7}, if $d_i < e^{-\frac{(R_{\widehat\alpha^i}-1)^2}{8p_0}}$ for large $i,$ then for large $i,$ we will have 
$$\mathcal{D}_{\mathscr{T}^{\widehat {\alpha}^i, \widehat{\mathscr{R}}^i \mathscr{C}^i}_{1 + L}}(M^i_{\tau_0 +L}), \mathcal{D}_{\mathscr{T}^{\widehat {\alpha}^i, \widehat{\mathscr{R}}^i \mathscr{C}^i}_{1 + 2L}}(M^i_{\tau_0 + 2L}) \leq \kappa_1 |\widehat {\alpha}^i|^2.$$
where $\kappa_1$ is from Proposition ~\ref{prop 6.4}. Therefore Proposition~\ref{prop 6.4} will imply
$$\mathcal{A}(M^i_{\tau_0 + L })^{\gamma} - \mathcal{A}(M^i_{\tau_0 + 2L})^{\gamma} \geq C^{-1} |\widehat\alpha^i|^2.$$
Let $s \in [0,1].$ By Lemma ~\ref{lemma 6.6}, we have
\begin{align*}
&D_{\mathcal{W}^{\widehat{\alpha}^i(1 + L-s),\widehat{\mathscr{R}}^i\mathscr{C}^i}}
(M^i_{\tau_0+L-s})
+
D_{\mathcal{W}^{\widehat{\alpha}^i(1 +L-s),\widehat{\mathscr{R}}^i\mathscr{C}^i}}
(M^i_{\tau_0+2L-s})
\\
&\leq
C\Big(
D_{\mathscr{T}^{\widehat{\alpha}^i, \widehat{\mathscr{R}}^i\mathscr{C}^i}_{L-s}}
(M^i_{\tau_0+L-s-1})
+
D_{\mathscr{T}^{\widehat{\alpha}^i, \widehat{\mathscr{R}}^i\mathscr{C}^i}_{2L-s}}
(M^i_{\tau_0+2L-s-1})
+
|\widehat{\alpha}^i|^2
\Big)
\\
&\leq
C\Big(
\mathcal{D}_{\mathscr{T}^{\widehat{\alpha}^i, \widehat{\mathscr{R}}^i\mathscr{C}^i}_{1+L}}
(M^i_{\tau_0+L})
+
\mathcal{D}_{\mathscr{T}^{\widehat{\alpha}^i, \widehat{\mathscr{R}}^i\mathscr{C}^i}_{1+2L}}
(M^i_{\tau_0+2L}) +
|\widehat{\alpha}^i|^2
\Big)
\\
&\leq
C\Big(\kappa_1 |\widehat {\alpha}^i|^2 + |\widehat {\alpha}^i|^2  \Big).
\end{align*}
Therefore,
$$D(M^i_{\tau_0 + L-s}, M^i_{\tau_0 + 2L-s}) \leq C \Big(\mathcal{A}(M^i_{\tau_0 + L })^{\gamma} - \mathcal{A}(M^i_{\tau_0 + 2L})^{\gamma} \Big).$$
Taking sup over $s \in [0,1]$, we get
$$\mathcal{D}(M^i_{\tau_0 + L}, M^i_{\tau_0 + 2L}) \leq C \Big(\mathcal{A}(M^i_{\tau_0 + L })^{\gamma} - \mathcal{A}(M^i_{\tau_0 + 2L})^{\gamma} \Big).$$
Hence we can assume that 
$$\mathcal{D}_{\mathscr{T}^{\widehat \alpha^i, \widehat{\mathscr{R}}^i\mathscr{C}^i}_{1 + s}}(M^i_{\tau_0 + s}) \leq C d_i,$$
for possibly a larger $C$ for any $s \in [1,2L].$ 

Then as in Lemma~\ref{lemma 5.6}, we can consider the new rescaling function $d_i^{-1} \widehat u_i$ which will converge to $u^{\perp}$ where, 
$$u^{\perp} = u - u_0.$$
Since $u_0$ has degree zero, the function $u^{\perp}$ also satisfies the inequalities 
$$I(u^\perp,0) \leq Ce^{\tilde{\delta_2} L},$$
and
$$I(u^\perp,2L) \leq Ce^{\tilde{\delta_1} L}.$$
for a possibly larger constant $C$ but independent of $L$. 

Let $\epsilon >0.$ If $L$ is chosen sufficiently large depending on $\epsilon$, then we will show that $I(u^\perp, L-1) \le \epsilon.$ 

Let $I(u^\perp,k+1) \leq e^{-\delta_2} I(u,k)$ for $k = 0,1,...,L-2$ then we have,
$$I(u^\perp,L-1) \leq e^{-(L-1)\delta_2} I(u^\perp ,0) \leq Ce^{\tilde{\delta_2} L} e^{-(L-1)\delta_2}.$$
So, $I(u^\perp, L-1) \le \epsilon$ if $L$ is sufficiently large.

Hence, we can assume that for some $k \leq L-1,$ $I(u^\perp, k) > e^{-\delta_2} I(u^\perp, k-1)$, i.e. $I(u^\perp, k-1) < e^{\delta_2} I(u^\perp, k).$ Then by Lemma ~\ref{lemma 5.5} $I(u^\perp, k+1) \geq e^{\delta_2} I(u^\perp,k)$. But then we can use Lemma ~\ref{lemma 5.4} to conclude $I(u^\perp,k+1) \geq e^{\delta_2} I(u^\perp,k)$ for $k = L, L+1,..., 2L-1$ which implies
$$I(u^{\perp}, L-1) \leq e^{-L\delta_2} I(u^{\perp},2L) \leq Ce^{-(L+1)\delta_2 + \tilde{\delta_1} L}.$$
Therefore, $I(u^\perp, L-1) \le \epsilon$ if $L$ is sufficiently large.

Hence, using Corollary \ref{cor 5.3}, either we have
$$\mathcal{D}_{\mathscr {T}^{\widehat\alpha^i, \widehat{\mathscr{R}}^i\mathscr{C}}_{1+L}}(M^i_{\tau_0 + L}) \leq C_0e^{-\frac{(R_{\widehat\alpha^i}-1)^2}{8p_0}},$$  
which implies for sufficiently large $i,$ 
$$\mathcal{D}_{\mathscr {T}^{\widehat\alpha^i, \widehat{\mathscr{R}}^i\mathscr{C}^i}_{1+L}}(M^i_{\tau_0 + L}), \mathcal{D}_{\mathscr {T}^{\widehat\alpha^i, \widehat{\mathscr{R}}^i\mathscr{C}^i}_{1+2L}}(M^i_{\tau_0 + 2L}) \leq \kappa_1|\widehat\alpha^i|^2.$$ 
This will imply by Proposition~\ref{prop 6.4} and arguing as before that
$$\mathcal{D}(M^i_{\tau_0 + L}, M^i_{\tau_0 + 2L}) \leq C\Big(\mathcal{A}(M^i_{\tau_0 +L })^\gamma- \mathcal{A}(M^i_{\tau_0 + 2L})^\gamma \Big).$$
Therefore for any $\delta$ there exists a constant $C_L$ depending on $L$ such that for sufficiently large $i$ we have a $R_0$ such that 
$$d_i = \mathcal{D}_{\mathscr {T}^{\widehat{\alpha}^i, \widehat{\mathscr{R}}^i \mathscr{C}^i}_{1 + L}}(M^i_{\tau_0 + L}) \leq \delta C_L d_i + C_0 \bigg( \int_{M^i_{\tau_0+L} \cap B_{R_0}(0)} d^2_{\mathscr{T}^{\widehat{\alpha}^i, \widehat{\mathscr{R}}^i \mathscr{C}^i}_{1 + L}} e^{-\frac{|x|^2}{4}} d\mathcal{H}^n \bigg)^{\frac{1}{2}}.$$
Now we will proceed as in the proof of Lemma~\ref{lemma 5.6}. For a given $R_0,$ for large $i$, using the local smooth convergence of the $d_i^{-1}u_i$ to $u$, we have
$$\int_{M^i_{\tau_0+L} \cap B_{R_0}(0)} d^2_{\mathscr{T}^{\widehat{\alpha}^i, \widehat{\mathscr{R}}^i \mathscr{C}^i}_{1+L}} e^{-\frac{|x|^2}{4}} d\mathcal{H}^n \leq Cd_i^2\epsilon^2.$$
We first choose $\epsilon$ such that $C_0C^{1/2}\epsilon < \frac{1}{4}.$ The choice of $\epsilon$ determines an $L$, such that
we have the estimate $I(u,L-1) \leq \epsilon.$ Choosing $L$ determines the constant $C_L$ and now we choose $\delta$ such that
$\delta C_L < \frac{1}{4}.$ The choice of $\delta$ determines the $R_0$ and then for sufficiently large $i$ it implies $d_i \leq \frac{d_i}{2}$, which is a contradiction.

Let us suppose condition $(\textnormal{a})$ holds for all $i,$ so 
$$\mathcal{D}_{\mathscr{T}^{\alpha^i, \mathscr{C}^i}_{1+2L}} (M^i_{\tau_0 + 2L}) \geq e^{\tilde{\delta_1} L} \mathcal{D}_{\mathscr{T}^{\alpha^i, \mathscr{C}^i}_{1+L}} (M^i_{\tau_0 + L}),$$
for $|\alpha^i| \to 0$. We want to apply Lemma~\ref{lemma 5.6} repeatedly to estimate
$D^{M^i}_{\mathscr{T}^{\alpha^i, \mathscr{C}^i}} (\tau _0+ kL)$ for $k = 2, 3, ...$ as long as $1 + kL \leq T_{\alpha^i}.$

There are two cases we need to consider. 

First, Suppose that as long as $1+ kL \leq T_{\alpha^i},$ we have 
$$\mathcal{D}_{\mathscr{T}^{\alpha^i, \mathscr{C}^i}_{1+kL}}(M^i_{\tau_0 + kL}) < \epsilon_0,$$ 
where $\epsilon_0$ is from Lemma~\ref{lemma 5.6}. Then either we can apply the three-annulus lemma $L^{-1}(T_{\alpha^i}-1)$ times to find that
$$\mathcal{D}_{\mathscr{T}^{\alpha^i, \mathscr{C}^i}_{1+L}} (M^i_{\tau_0 + L}), \mathcal{D}_{\mathscr{T}^{\alpha^i, \mathscr{C}^i}_{1 + 2L}} (M^i_{\tau_0 + 2L}) \leq C\epsilon_0 e^{-\tilde{\delta_2} (T_{\alpha^i}-1)}.$$
By definition,
$$ T_\alpha := |\alpha|^{-\frac{1}{8}}.$$
So it follows that for sufficiently large $i$, we have
$$\mathcal{D}_{\mathscr{T}^{\alpha^i, \mathscr{C}^i}_{1+L}} (M^i_{\tau_0 + L}), \mathcal{D}_{\mathscr{T}^{\alpha^i, \mathscr{C}^i}_{1 + 2L}} (M^i_{\tau_0 + 2L}) \leq \kappa_1 |\alpha^i|^2,$$
where $\kappa_1$ is from Proposition~\ref{prop 6.4}. Therefore by Proposition~\ref{prop 6.4},
$$\mathcal{A}(M^i_{\tau_0+L})^{\gamma} - \mathcal{A}(M^i_{\tau_0 + 2L})^{\gamma} \geq \kappa_1 |\alpha^i|^2,$$
which implies arguing as before,
$$\mathcal{D}(M^i_{\tau_0 +L}, M^i_{\tau_0 + 2L}) \leq C \Big(\mathcal{A}(M^i_{\tau_0+L})^{\gamma} - \mathcal{A}(M^i_{\tau_0 + 2L})^{\gamma}  \Big),$$
Next, we consider the remaining case. Namely, suppose that for some $k_i$ with $1+ k_iL \leq T_{\alpha^i}$, we have 
$$\mathcal{D}_{\mathscr{T}^{\alpha^i, \mathscr{C}^i}_{1+k_iL}} (M^i_{\tau_0 + k_iL}) > \epsilon_0,$$ 
and $k_i$ is the smallest such choice for each $i$. We have 
$$\mathcal{D}_{\mathscr{T}^{\alpha^i, \mathscr{C}^i}_{1 + k_iL}} (M^i_{\tau_0 + k_iL}) \geq e^{\tilde{\delta_1} L} \mathcal{D}_{\mathscr{T}^{\alpha^i, \mathscr{C}^i}_{1 + (k_i-1)L}} (M^i_{\tau_0 + (k_i-1)L}),$$ which implies 
\begin{align*}
\mathcal{D}_{\mathscr{T}^{\alpha^i, \mathscr{C}^i}_{1 + k_iL}} (M^i_{\tau_0 + k_iL}) - \mathcal{D}_{\mathscr{T}^{\alpha^i, \mathscr{C}^i}_{1 + (k_i-1)L}} (M^i_{\tau_0 + (k_i-1)L}) &\geq (1 - e^{-\tilde{\delta_1} L}) \epsilon_0 \\
&\geq \frac{1}{2} \epsilon_0,
\end{align*}
by choosing $L$ large. Define the shifted rescaled flows by
\[
\widetilde{M}^i_s := M_{\tau_0 + (k_i-1)L + s}, \qquad 
\widetilde{T}^i_s := \mathscr{T}^{\alpha^i,\mathscr C^i}_{1 + (k_i-1)L + s}.
\]
So we obtain a sequence of flows $\widetilde{M}^i_s$ and $\widetilde{T}^i_s$ defined for $s \in [-2,L]$, which satisfy 
\begin{align*}
    \mathcal{A}^{\gamma} (\widetilde{M}^i_{-2}) - \mathcal{A}^{\gamma} (\widetilde{M}^i_L) &\leq \epsilon_i,\\
    \mathcal{D}_{\widetilde{T}^i_L} ( \widetilde{M}^i_L) - \mathcal{D}_{\widetilde{T}^i_0} ( \widetilde{M}^i_0) &\ge \frac{\varepsilon_{0}}{2},\\
    \mathcal{D}_{\widetilde{T}^i_0} ( \widetilde{M}^i_0) &< \epsilon_0.
\end{align*}
for a sequence $\epsilon_i \to 0.$ First note that the sequence of flow $\widetilde{T}^i_s$ converges to $\mathscr{C}$ as $i \to \infty.$ Since $\epsilon_i\to0$, after passing to a subsequence, the flows
$\widetilde M^i_s$ converge to a static flow. Now for $\epsilon_0$ is sufficiently small, using rigidity of cylinder i.e. Proposition ~\ref{prop 6.7} and the fact  that $\mathcal{D}_{\widetilde{T}^i_0} ( \widetilde{M}^i_0) < \epsilon_0,$ we can conclude that $\widetilde{M}^i_s$ converges to 
some cylinder $\mathscr{RC}$ for a rotation $\mathscr{R}.$ In particular, for any $s \in [-2,L],$ the $\widetilde{M}^i_s$ converge locally smoothly to $\mathscr{RC}$ on compact sets. So for any $R>0,$ we have as $i \to \infty,$
$$\bigg(\int_{\widetilde{M}^i_L \cap B_R(0)} d_{\widetilde{T}^i_L}^2 e^{-\frac{|x|^2}{4}} d\mathcal{H}^n\bigg)^{\frac12} - \bigg(\int_{\widetilde{M}^i_0 \cap B_R(0)} d_{\widetilde{T}^i_0}^2 e^{-\frac{|x|^2}{4}} d\mathcal{H}^n\bigg)^{\frac12} \to 0.$$
By definition,
\begin{align*}
     \mathcal{D}_{\widetilde{T}^i_L}(\widetilde{M}^i_L) -
     \mathcal{D}_{\widetilde{T}^i_0}(\widetilde{M}^i_0)
     &= \left(\int_{\widetilde{M}^i_L\cap B_{R_{i,1}}(0)} d_{\widetilde{T}^i_L}^{\,2} e^{-\frac{|x|^2}{4}} d\mathcal{H}^n\right)^{\frac12} - \left(\int_{\widetilde{M}^i_0\cap B_{R_{i,2}}(0)} d_{\widetilde{T}^i_0}^{\,2} e^{-\frac{|x|^2}{4}} d\mathcal{H}^n\right)^{\frac12} \\ &\quad + e^{-\frac{R_{i,1}^2}{8p_0}} - e^{-\frac{R_{i,2}^2}{8p_0}}.
\end{align*}
where $R_{i,1}$ and $R_{i,2}$ are the largest values for which ${M}^i_L$ and ${M}^i_0$ are $\eta-$graphs over $\widetilde{T}^i_L$ and $\widetilde{T}^i_0,$ respectively. Since $R_{i,1}, R_{i,2} \to \infty$ as $i \to \infty$, it follows that 
$$\mathcal{D}_{\widetilde{T}^i_L} ( \widetilde{M}^i_L) - \mathcal{D}_{\widetilde{T}^i_0} (\widetilde{M}^i_0) \to 0.$$
which contradicts 
$$\mathcal{D}_{\widetilde{T}^i_L} ( \widetilde{M}^i_L) - \mathcal{D}_{\widetilde{T}^i_0} (\widetilde{M}^i_0) \ge \frac{\varepsilon_{0}}{2}.$$

Let us suppose condition $(\textnormal{b})$ holds for all $i,$ while condition $(\textnormal{a})$ fails. Again, there are two cases we need to consider.

First, let $\kappa_1 > 0$ be the constant from Proposition ~\ref{prop 6.4} and suppose that 
$$d_i \leq \kappa_1 e^{-\tilde{\delta_1} L}|\alpha^i|^2.$$ 
Then arguing similarly as before, we have
$$\mathcal{D}(M_{\tau_0 + L}, M_{\tau_0 + 2L}) \leq C\Big(\mathcal{A}(M_{\tau_0 +L})^{\gamma} - \mathcal{A}(M_{\tau_0 + 2L})^{\gamma} \Big).$$
Thus, it remains to consider the case
$$d_i \geq \kappa_1 e^{-\tilde{\delta_1} L}|\alpha_i|^2.$$ 
In addition,
$$\mathcal{D}_{\mathscr{T}^{\alpha^i, \mathscr{C}^i}_2} (M^i_{\tau_0 + 1}) \geq e^{\tilde{\delta_2} L} d_i,$$
and
$$\mathcal{D}_{\mathscr{T}^{\alpha^i, \mathscr{C}^i}_{1+2L}} (M^i_{\tau_0 + 2L}) \leq e^{\tilde{\delta_1} L} d_i.$$
We can use Lemma ~\ref{lemma 6.6} to estimate $\mathcal{D}(M^i_{\tau_0 + L-s}, M^i_{\tau_0 + 2L-s})$ from above. Let $s \in[0,1].$
\begin{align*}
    D(M^i_{\tau_0 + L-s}, M^i_{\tau_0 + 2L-s}) &\leq D_{\mathcal{W}^{\alpha^i(1+L-s),\mathscr{C}^i}} (M^i_{\tau_0 + L-s}) + D_{\mathcal{W}^{\alpha^i(1+L-s),\mathscr{C}^i}} (M^i_{\tau_0 + 2L -s}) \\
    &\leq C \Big(D_{\mathscr{T}^{\alpha^i, \mathscr{C}^i}_{L-s}} (M^i_{\tau_0 + L - s -1}) + D_{\mathscr{T}^{\alpha^i, \mathscr{C}^i}_{2L-s}} (M^i_{\tau_0 + 2L- s -1}) + |\alpha^i|^2\Big)\\
    &\leq C \Big(\mathcal{D}_{\mathscr{T}^{\alpha^i, \mathscr{C}^i}_{1+L}} (M^i_{\tau_0 + L}) + \mathcal{D}_{\mathscr{T}^{\alpha^i, \mathscr{C}^i}_{1+2L}} (M^i_{\tau_0 + 2L}) + d_i\kappa_1^{-1} e^{\tilde{\delta_1} L} \Big)\\
    &\leq Cd_i(1 + e^{\tilde{\delta_1}L}).
    \end{align*}
where the factor ${\kappa_1}^{-1}$ has been absorbed into the constant $C$. Since $s \in [0,1]$ arbitrary, 
$$\mathcal{D}(M^i_{\tau_0}, M^i_{\tau_0 + L}) \leq Cd_i(1 + e^{\tilde{\delta_1}L}).$$
We claim that for sufficiently large $L$, we have 
$$\mathcal{D}(M^i_{\tau_0 + L}, M^i_{\tau_0 + 2L}) \leq \frac{1}{2} \mathcal{D}(M^i_{\tau_0}, M^i_{\tau_0 + L}).$$
once $i$ is large enough. Suppose that this were not the case, i.e.
$$\mathcal{D}(M^i_{\tau_0}, M^i_{\tau_0 + L}) < 2 \mathcal{D}(M^i_{\tau_0+L}, M^i_{\tau_0 + 2L}).$$
Let $s \in [0,1].$ By definition there are rotations $\widehat{\mathscr{R}}^i$, and $\widehat\alpha^i$, such that
$$D_{\mathcal{W}^{\widehat \alpha^i, \widehat{\mathscr{R}}^i\mathscr{C}^i}}(M^i_{\tau_0-s}) + D_{\mathcal{W}^{\widehat \alpha^i, \widehat{\mathscr{R}}^i\mathscr{C}^i}}(M^i_{\tau_0+L-s}) < 2 D(M^i_{\tau_0+L}, M_{\tau_0 + 2L}).$$
Now using 
$$D_{\mathcal{W}^{\widehat \alpha^i, \widehat{\mathscr{R}}^i\mathscr{C}^i}}(M ^i_{\tau_0 + L-s}) \leq 2 D(M^i_{\tau_0+L}, M^i_{\tau_0 + 2L}),$$
we have
$$D_{\mathcal{W}^{\widehat \alpha^i, \widehat{\mathscr{R}}^i\mathscr{C}^i}}(M^i_ {\tau_0+L-s}) \leq Cd_i(1 + e^{\tilde{\delta_1}L}).$$
By Lemma ~\ref{lemma 6.6},
$$D_{\mathscr{T}_0^{\widehat \alpha^i, \widehat{\mathscr{R}}^i\mathscr{C}^i}}(M^i_ {\tau_0+L-s}) \leq Cd_i(1 + e^{\tilde{\delta_1}L}).$$
Using the proof of Proposition ~\ref{prop 5.2},
$$d_{\mathscr{T}_1^{\widehat \alpha^i, \widehat{\mathscr{R}}^i\mathscr{C}^i}} \leq C d_i(1 + e^{\tilde{\delta_1}L})$$
on $M_{\tau_0+L-s+1} \cap B_1(0).$ Applying the same argument to
$\mathscr{T}^{\alpha^i,\mathscr{C}^i}$ and using the definition $d_i = \mathcal{D}_{\mathscr{T}^{\alpha^i, \mathscr{C}^i}_{1+L}} (M^i_{\tau_0 + L})$, we obtain
$$d_{\mathscr{T}^{\alpha^i, \mathscr{C}^i}_{2+L-s}} \leq Cd_i.$$
on $M_{\tau_0+L-s+1} \cap B_1(0).$ Therefore, 
$$|\widehat{\mathscr{R}}^i - Id|, |\alpha^i(2+L-s) - \widehat \alpha^i(1)| \leq Cd_i(1 + e^{\tilde{\delta_1}L}).$$
Using the estimate $|\alpha^i(\tau)-\alpha^i(0)|\leq C|\alpha^i|^2,$
we have
\begin{align*}
    |\alpha^i(1-s) - \widehat{\alpha}^i(0)| &\leq C |\alpha^i(2-s) - \widehat{\alpha}^i(1)| \\
    &\leq C|\alpha^i(2+L-s) - \widehat \alpha^i(1)| + C|\alpha^i(2+L-s) -\alpha^i(2-s)|\\
    &\leq C d_i(1 + e^{\tilde{\delta_1} L}).
\end{align*}
Now using 
$$D_{\mathcal{W}^{\widehat \alpha^i, \widehat{\mathscr{R}}^i\mathscr{C}^i}}(M^i_ {\tau_0-s}) \leq 2 D(M^i_{\tau_0+L}, M^i_{\tau_0 + 2L}),$$
and Lemma ~\ref{lemma 6.6} we have,
$$D_{\mathscr{T}^{\alpha^i, \mathscr{C}^i}_{1-s}} (M^i_{\tau_0 - s}) \leq C d_i(1 + e^{\tilde{\delta_1} L}).$$
Since $s\in[0,1]$ was arbitrary, the estimate holds for all $s\in[0,1]$. Applying Proposition~\ref{prop 5.2}, we obtain
$$\mathcal{D}_{\mathscr{T}^{\alpha^i, \mathscr{C}^i}_2} (M^i_{\tau_0+1}) \leq C d_i(1 + e^{\tilde{\delta_1} L})$$
increasing $C$ further if necessary. At the same time, we have 
$$\mathcal{D}_{\mathscr{T}^{\alpha^i, \mathscr{C}^i}_2} (M^i_{\tau_0+1}) \geq e^{\tilde{\delta_2} L} d_i.$$ 
So we obtain 
$$e^{\tilde{\delta_2} L} d_i \leq C d_i(1 + e^{\tilde{\delta_1} L})$$ 
for sufficiently large $i$. Since $\tilde{\delta_2} > \tilde{\delta_1}$, this is a contradiction if $L$ is sufficiently large.
\end{proof}

We can now prove our main Theorem~\ref{thm 1.1}. The argument follows exactly as in \cite[Theorem~6.7]{szekelyhidi2020uniqueness}. Assume that $\mathcal A^\gamma (M_{\tau_0 - 1}) - \lim\limits_{\tau \to \infty} \mathcal A^\gamma(M_\tau)$ is sufficiently small. Suppose that
\(M_{\tau_0+kL-s}\) is sufficiently close to \(\mathscr C\), in the sense that
\(d(M_{\tau_0+kL-s})\) is sufficiently small for every
\(s\in[0,1]\) and \(k=0,1,\ldots,N\). We apply Proposition~\ref{prop 6.4}
repeatedly. Let \(k_1,\ldots,k_m\) denote the values of \(k\) for which
alternative~(ii) holds, while alternative~(i) holds for the remaining values.
Then
\begin{align*}
\sum_{k=0}^{N}
\mathcal{D}(M_{\tau_0+kL},M_{\tau_0+(k+1)L})
&\leq
2\mathcal{D}(M_{\tau_0},M_{\tau_0+L}) +2\sum_{j=1}^{m}
\mathcal{D}(M_{\tau_0+k_jL},M_{\tau_0+(k_j+1)L})
\\
&\leq
2\mathcal{D}(M_{\tau_0},M_{\tau_0+L})
+2C\sum_{j=1}^{m}
\Big(
\mathcal{A}^{\gamma}(M_{\tau_0+k_jL})
-
\mathcal{A}^{\gamma}(M_{\tau_0+(k_j+1)L})
\Big)
\\
&\leq
2\mathcal{D}(M_{\tau_0},M_{\tau_0+L})
+
2C
\Big(
\mathcal{A}^{\gamma}(M_{\tau_0})
-
\lim_{\tau\to\infty}\mathcal{A}^{\gamma}(M_\tau)
\Big).
\end{align*}
where the last inequality follows from the monotonicity of
\(\mathcal A^\gamma\) along the rescaled flow. Denote \(d_V(\mu_{M_1},\mu_{M_2})\) simply by \(d(M_1,M_2)\).
Using Lemma~\ref{lemma 6.2}, we obtain
\begin{align*}
    \sup_{s\in [0,1]} d(M_{\tau_0 -s}, M_{\tau_0 + (N+1)L -s}) &\leq \sup_{s\in [0,1]} \left( \sum_{k=0}^{N} d(M_{\tau_0 + kL-s}, M_{\tau_0 + (k+1)L-s}) \right) \\
    &\leq \sum_{k=0}^{N} \sup_{s\in [0,1]} d(M_{\tau_0 + kL-s}, M_{\tau_0 + (k+1)L-s}) \\
    &\leq C \sum_{k=0}^{N} \mathcal{D}(M_{\tau_0 + kL}, M_{\tau_0 + (k+1)L}) \\
    &\leq C \mathcal{D}(M_{\tau_0}, M_{\tau_0 + L}) + C\Big(\mathcal{A}^{\gamma}(M_{\tau_0}) - \lim\limits_{\tau \to \infty} \mathcal{A}^\gamma(M_\tau)\Big).
\end{align*}
If \(M_{\tau_0-s}\) is sufficiently close to \(\mathscr C\) for every
\(s\in[0,1]\), then the estimate above implies that
\(M_{\tau_0+(N+1)L-s}\) is also sufficiently close to \(\mathscr C\) for every
\(s\in[0,1]\), provided $\mathcal A^\gamma(M_{\tau_0-1}) - \lim\limits_{\tau \to \infty} \mathcal A^\gamma(M_\tau)$ is sufficiently small. Hence Proposition~\ref{prop 6.4} can be applied
repeatedly, and therefore \(N\) may be taken arbitrarily large.

Since \(\mathscr C\) is a tangent flow, there exists a sequence
\(\tau_i\to\infty\) such that
\[
\sup_{s \in [0,1]}d(M_{\tau_i-s}, \mathscr{C})\to 0.
\]
Moreover, by Huisken's monotonicity formula,
\[
\mathcal A^\gamma(M_{\tau_i})
-
\lim_{\tau\to\infty}
\mathcal A^\gamma(M_\tau)
\to 0.
\]
Applying the preceding estimate with \(\tau_0=\tau_i\) and letting
\(i\to\infty\), we conclude that
\[
\lim_{i\to\infty}
\sup_{\tau \ge\tau_i} \sup_{s \in [0,1]}
d(M_{\tau -s}, \mathscr{C})
=
0.
\]
Therefore, \(\mathscr C\) is the unique tangent flow.

\section{Appendix}\label{section 7}

We begin by recording a simple approximation result used in the proof of Proposition~\ref{lemma 6.2}.
\begin{lemma} \label{prop 7.1}
There exists a countable subset
\[
\{f_j\}_{j=1}^\infty
\subset
\{f\in C_c^1(\mathbb{R}^{n+1}):\|f\|_{C^0}\le1\}
\]
that is dense in the unit ball of $C_c^0(\mathbb{R}^{n+1})$ with respect to the $C^0$ norm and satisfies
\[
\left(\frac32\right)^{-j}\|f_j\|_{C^1}\le1
\]
for every $j \in\mathbb N$.
\end{lemma}

\begin{proof}
Let
\[
B=\{f\in C_c^1(\mathbb{R}^{n+1}):\|f\|_{C^0}\le1\}.
\]
Since $B$ is dense in the unit ball of $C_c^0(\mathbb{R}^{n+1})$ with respect to the $C^0$ norm, there exists a countable dense subset
\[
\{g_i\}_{i=1}^{\infty}\subseteq B.
\]
For each $i\in\mathbb N$, we have $\|g_i\|_{C^1}<\infty$. Since
\[
\left(\frac32\right)^{-n}\to0
\quad\text{as }n\to\infty,
\]
there exists an integer $M_i$ such that
\[
\left(\frac32\right)^{-n}\|g_i\|_{C^1}\le1
\]
for every $n\ge M_i$.

Define a strictly increasing sequence of integers
$\{N_i\}_{i=1}^{\infty}$ recursively by
\[
N_1=M_1,
\]
and
\[
N_i=\max\{M_i,N_{i-1}+1\},
\qquad i\ge2.
\]
Then
\[
N_1<N_2<\cdots,
\]
and
\[
\left(\frac32\right)^{-N_i}\|g_i\|_{C^1}\le1
\]
for every $i\in\mathbb N$.

Define
\[
f_j=
\begin{cases}
g_i,&\text{if } j=N_i\text{ for some }i,\\
0,&\text{otherwise}.
\end{cases}
\]

Since
\[
\{g_i:i\in\mathbb N\}\subseteq\{f_j: j \in\mathbb N\},
\]
the set $\{f_j: j\in\mathbb N\}$ is dense in $B$. As $B$ is dense in the unit ball of
$C_c^0(\mathbb{R}^{n+1})$ with respect to the $C^0$ norm, it follows that
$\{f_j: j \in\mathbb N\}$ is also dense in the unit ball of
$C_c^0(\mathbb{R}^{n+1})$.

By construction, if $j=N_i$ for some $i$, then
\[
\left(\frac32\right)^{-j}\|f_j\|_{C^1}
=
\left(\frac32\right)^{-N_i}\|g_i\|_{C^1}
\le1.
\]
Otherwise, $f_j=0$, and hence
\[
\left(\frac32\right)^{-j}\|f_j\|_{C^1}=0.
\]
Therefore,
\[
\left(\frac32\right)^{-j}\|f_j\|_{C^1}\le1
\]
for every $j \in\mathbb N$.
\end{proof}

The following elementary lemma allows one to compare graphical representations over nearby hypersurfaces.
\begin{lemma}\label{lemma 7.2}
Let $M\subset \mathbb{R}^{n+1}$ be a smooth hypersurface satisfying $|A_M|\leq C_0$ and positive normal injectivity radius.  Suppose that $u,v\in C^2(M)$ satisfy $\|u\|_{C^2(M)}$, $\|v\|_{C^2(M)}\leq \epsilon ,$
where $\epsilon>0$ is sufficiently small. Assume that $N$ is an
$\frac{\eta}{2}$-graph over $M_u$. Then $N$ is an
$\eta$-graph over $M_v$.
\end{lemma}

\begin{proof}
Let
\[
X_u(x)=x+u(x)\nu_M(x),\qquad
X_v(x)=x+v(x)\nu_M(x),
\]
be the normal graph parametrizations of $M_u$ and $M_v$, respectively.

Since $|A_M|\leq C_0$, $M$ has positive normal injectivity radius, and
$\|u\|_{C^2(M)},\|v\|_{C^2(M)}$ are sufficiently small, the hypersurfaces
$M_u$ and $M_v$ have uniformly bounded geometry. Therefore, there exists
$r>0$ such that the normal coordinate maps
\[
G_u,G_v:M\times(-r,r)\to \mathbb{R}^{n+1},
\]
defined by
\[
G_u(x,s)=X_u(x)+s\nu_u(x),\qquad
G_v(x,s)=X_v(x)+s\nu_v(x),
\]
are $C^2$ diffeomorphisms onto tubular neighborhoods of $M_u$ and $M_v$.

Since $N$ is an $\frac{\eta}{2}$-graph over $M_u$, there exists
$u_N\in C^2(M_u)$ such that
\[
N=\{p+u_N(p)\nu_u(p):p\in M_u\},
\]
with
\[
\|u_N\|_{C^2(M_u)}\leq \frac{\eta}{2}.
\]
Pulling back by $X_u$, we can write
\[
N=\{G_u(x,\widetilde u_N(x)):x\in M\},
\]
where
\[
\|\widetilde u_N\|_{C^2(M)}
\leq C\|u_N\|_{C^2(M_u)}.
\]
Moreover, since $X_u$ converges to the identity in $C^2$ as
$\|u\|_{C^2}\to0$, the constant $C$ can be chosen arbitrarily close to
$1$ by taking $\epsilon$ sufficiently small. Hence, after decreasing
$\epsilon$ if necessary,
\[
\|\widetilde u_N\|_{C^2(M)}\leq \frac{3\eta}{4}.
\]
We continue to denote $\widetilde u_N$ by $u_N$.

For $\epsilon$ sufficiently small, $M_u$ and $M_v$ are close enough so that
$N$ lies in the tubular neighborhood of $M_v$. Consider
\[
\Phi=G_v^{-1}\circ G_u.
\]
Write $\Phi=(\Phi_1,\Phi_2)$. Since $u$ and $v$ are sufficiently small in $C^2(M)$, both $G_u$ and $G_v$ are $C^2$-close to the normal coordinate map of $M$. Consequently,
$\Phi$ is $C^2$-close to the identity map.
Hence,
\[
\partial_s\Phi_2(x,s)>0
\]
for all $(x,s)$ in the tubular neighborhood, provided $\epsilon$ is
sufficiently small. Hence the image of the graph of $u_N$ under $\Phi$ is
again the graph of a $C^2$ function over $M$,
\[
\Phi(\operatorname{graph}(u_N))
=
\operatorname{graph}(u_v)
\]
for some $u_v\in C^2(M)$. Equivalently,
\[
N=\{G_v(x,u_v(x)):x\in M\}.
\]
Since $\Phi$ is $C^2$-close to the identity, the graphing
function $u_v$ depends continuously on $(u,v,u_N)$ in the
$C^2$ topology. Consequently, after choosing $\epsilon$
sufficiently small,
\[
\|u_v-u_N\|_{C^2(M)}
<
\frac{\eta}{4}.
\]
Therefore, 
\[
\|u_v\|_{C^2(M)}
\leq \frac{3\eta}{4}+\frac{\eta}{4}
=\eta .
\]
Thus $N$ is an $\eta$-graph over $M_v$.
\end{proof}

\bibliographystyle{plain}
\bibliography{references}

\end{document}